\documentclass[11pt]{amsart}
\usepackage{amsmath,amssymb,amsthm}
\usepackage{graphicx}
\usepackage{booktabs}
\usepackage[margin=1in]{geometry}
\usepackage{hyperref}

\newtheorem{theorem}{Theorem}[section]
\newtheorem{proposition}[theorem]{Proposition}
\newtheorem{lemma}[theorem]{Lemma}

\theoremstyle{definition}
\newtheorem{definition}[theorem]{Definition}
\theoremstyle{remark}
\newtheorem{remark}[theorem]{Remark}

\newcommand{\ve}{\varepsilon}
\newcommand{\R}{\mathbb{R}}

\newcommand{\Lop}{\mathcal{L}}

\title[The high activation energy limit: rate and a convergent scheme]{The high activation energy limit of an equidiffusional\\ flame model: rate of convergence and\\ a convergent front-following scheme}

\author{Juli\'an Fern\'andez Bonder}
\address{Instituto de C\'alculo (UBA--CONICET) and Departamento de Matem\'atica, FCEyN,
Universidad de Buenos Aires, Ciudad Universitaria, (1428) Buenos Aires, Argentina}
\email{jfbonder@dm.uba.ar}

\date{\today}


\begin{document}

\begin{abstract}
We study the high activation energy limit of an equidiffusional model of
premixed flames, whose limit is a free boundary problem for the heat equation
in which the free boundary condition depends on perturbations of order $\ve$ of
the data. We prove a rate of convergence for this singular limit, of order
$\ve$ up to a logarithm, by means of explicit barriers whose error does not
depend on the strictness of the super- and subsolutions that they approximate.
The same barriers, evaluated on a mesh, are discrete barriers for a monotone
scheme: for radial data, and a perturbation that need not be constant, we prove that a
fully discrete scheme, on meshes with
a fine band that follows the front with a prescribed velocity, converges to
the solution of the free boundary problem, with an explicit error bound.
Numerical experiments, including a perturbation that is not constant,
illustrate the rate and the scheme.
\end{abstract}

\maketitle

\section{Introduction}

The system
\begin{equation}\label{eq:system}
T^\ve_t=\Delta T^\ve+Y^\ve f_\ve(T_f-T^\ve),\qquad
Y^\ve_t=\Delta Y^\ve-Y^\ve f_\ve(T_f-T^\ve),\qquad f_\ve(s)=\frac1{\ve^2}f\Big(\frac s\ve\Big),
\end{equation}
is the constant density (thermo-diffusive) model of a premixed flame with one-step
chemistry and equal diffusivities of heat and reactant (Lewis number one)
\cite{BuckmasterLudford1982,MatkowskySivashinsky1979}: $T^\ve$ is the
temperature and $Y^\ve\ge0$ the mass fraction of the reactant, scaled so that the
reaction conserves $T^\ve+Y^\ve$; $T_f$ is a reference temperature; $\ve$ is
essentially the inverse of the activation energy (the inverse of the Zeldovich number),
and the factor $\ve^{-2}$ is the scaling for which the flame speed is of order one. The
function $f$ is Lipschitz, positive on $(-\infty,1)$
and zero on $[1,\infty)$, so that the reaction only takes place where
$T^\ve>T_f-\ve$: the kinetics has a cut-off (an ignition temperature) at $T_f-\ve$. In
terms of $u^\ve=T_f-T^\ve$, \eqref{eq:system} reads
$\Delta u^\ve-u^\ve_t=Y^\ve f_\ve(u^\ve)=\Delta Y^\ve-Y^\ve_t$, so that
$w^\ve=Y^\ve-u^\ve$ solves the heat equation and \eqref{eq:system} reduces to
\begin{equation}\label{eq:Peps}\tag{$P_\ve$}
u_t=\Delta u-(u+w^\ve)f_\ve(u),\qquad w^\ve_t=\Delta w^\ve.
\end{equation}
The function $w^\ve=Y^\ve+T^\ve-T_f$ is the excess of total enthalpy with respect to the
reference: where the reactant is exhausted the temperature is $T_f+w^\ve$, so that $T_f$
is the adiabatic flame temperature exactly when $w^\ve\equiv0$, that is, $u^\ve=Y^\ve$.

The free boundary problem obtained in the limit was introduced by Caffarelli
and V\'azquez \cite{CaffarelliVazquez1995} (see also \cite{Vazquez1996}): in
the case $u^\ve=Y^\ve$ they proved uniform estimates for the singular
perturbation problem, passed to the limit
$\ve\to0$, and constructed explicit solutions, among them the self-similar
solutions that we use in Section~\ref{sec:exact}. The limit was further studied by
Caffarelli, Lederman and Wolanski
\cite{CaffarelliLedermanWolanski1997a} and, regarding
uniqueness, by Lederman, V\'azquez and Wolanski \cite{LedermanVazquezWolanski2001}
and Petrosyan \cite{Petrosyan2002}. In \cite{FernandezBonderWolanski2000} Wolanski
and the author showed that if $w^\ve/\ve\to w_0>-1$ then limits of $u^\ve$ solve
\begin{equation}\label{eq:P}\tag{$P$}
u_t=\Delta u\ \text{ in }\{u>0\},\qquad |\nabla u|=\sqrt{2M}\ \text{ on }\partial\{u>0\},\qquad
M=M(w_0)=\int_{-w_0}^1(s+w_0)f(s)\,ds,
\end{equation}
in a pointwise sense at regular points of the free boundary and, under additional
assumptions, in the viscosity sense, and in \cite{FernandezBonderWolanski2004} that, under
geometric assumptions, the limit is unique and coincides with the classical
solution when one exists. Thus the high activation energy limit is not determined by
the limit of $u^\ve_0$: perturbations of order $\ve$ of the initial data, which do
not change that limit, change the free boundary condition. In physical terms, $w_0$ is the
perturbation of the flame temperature measured in units of the width $\ve$ of the range of
temperatures where the reaction takes place (the scaling of near-equidiffusional flames
\cite{MatkowskySivashinsky1979}), and the dependence of $M$ on $w_0$ is the classical
sensitivity of the flame speed to the flame temperature
\cite{ZeldovichFrankKamenetskii1938}: formally, for the Arrhenius rate $f(s)=e^{-s}$, which
has no cut-off and is not covered by our assumptions, $M(w_0)=e^{w_0}M(0)$. The
restriction $w_0>-1$ means that the temperature of the burnt gas, $T_f+\ve w_0$, is above
the cut-off; it is a feature of the cut-off kinetics. Problem \eqref{eq:P} is a
one-phase problem because, to leading order, the burnt gas is at the uniform
temperature $T_f$; the set $\{u>0\}$ is the \emph{fresh} region.

The convergence results available for this limit
\cite{CaffarelliVazquez1995,CaffarelliLedermanWolanski1997a,FernandezBonderWolanski2000}
are obtained by compactness, and uniqueness of the limit
\cite{LedermanVazquezWolanski2001,Petrosyan2002,FernandezBonderWolanski2004} by comparison
with strict super- and subsolutions; none of them gives a rate. The same is true of
the finer analysis of the limit \cite{Weiss2003,PetrosyanYip2008,WeissZhang2009} and of
the recent regularity theory for singular perturbation problems of this type, in the
parabolic \cite{RicarteTeymurazyanUrbano2019,AudritoSanzPerela2024} and in the elliptic
case \cite{AudritoSerra2022,AraujoRicarteTeixeira2017}: the estimates are
uniform in $\ve$, but they do not measure the distance to the limit. For travelling waves, whose
high activation energy limit is studied in
\cite{BerestyckiNicolaenkoScheurer1985,Glangetas1992}, quantitative
information is available: the limit of the speed was
obtained in \cite{BerestyckiNicolaenkoScheurer1985}, and
a rigorous expansion up to order $\ve$ was recently proved for a related scalar
equation in \cite{JelbartKristiansenSzmolyan2025}.

Our first result is a rate of convergence for the evolution problem in any
dimension (Section~\ref{sec:rate}). We make the approximation theorems of
\cite{FernandezBonderWolanski2004} quantitative and prove that, if \eqref{eq:P} has a
classical solution $v$ that is regular up to its free boundary, and the data are
strictly starshaped, then
\[
-C\Big(\ve\log\frac1\ve+\eta_\ve\Big)\le u^\ve-v\le C(\ve+\eta_\ve),\qquad \eta_\ve=\|w^\ve/\ve-w_0\|_\infty
\]
(Theorem~\ref{q:thm:rate}). The setting is that of \cite{FernandezBonderWolanski2004}:
the fresh region $\{u_0>0\}$ is bounded and starshaped, a pocket of fresh gas surrounded
by burnt gas, which is consumed in finite time; flames with an unbounded fresh region
(planar or expanding fronts) are not covered. The order $\ve$ cannot be improved in the
sup norm when $w_0\not\equiv0$, for the simple reason that the temperature of the burnt
gas is $T_f+\ve w_0$; for travelling waves, which are outside this setting,
the position of the front also differs by a quantity of order $\ve$
(Remark~\ref{q:rem:rate}). The proof rests on
an explicit barrier, an inner profile composed with an extension of a strict
super- or subsolution, whose error does not depend on the strictness. Barriers of
this type are classical for Allen--Cahn type limits
\cite{Chen1992,AlfaroHilhorstMatano2008}. To our knowledge no rate has been stated for
the present singular limit beyond travelling waves. When $w^\ve/\ve$ is constant the
problem reduces to the scalar one of \cite{CaffarelliVazquez1995,Petrosyan2002}, and a
two-sided bound of order $\ve\log(1/\ve)$ can be extracted from the barriers of
\cite{Petrosyan2002}; what Theorem~\ref{q:thm:rate} adds is the upper bound of order
$\ve$, without logarithm, and the case in which $w^\ve/\ve$ is not constant. The
relation with that construction and with the one of \cite{FernandezBonderWolanski2004}
is discussed in Section~\ref{sec:rate}. For the rate we assume that $f$ is nonincreasing and of class
$C^2$ up to $1$, as $f(s)=(1-s)^+$; the second derivative is used only to differentiate
the inner profiles with respect to $w^\ve/\ve$, and is not needed when $w^\ve/\ve$ is
constant.

Our second result is the convergence of a fully discrete scheme to the free boundary
problem (Sections~\ref{sec:method} and \ref{sec:analysis}). The idea is simple: the
barriers of the rate theorem are explicit functions, so they can be evaluated on a
mesh, and if the scheme has a discrete comparison principle they can be used as discrete
barriers. The scheme that we use is linearly implicit in time: the reaction rate is
taken at the previous time level and the reactant at the new one. When $f$ is
nonincreasing this gives a comparison principle for every time step, and the error
analysis reduces to the construction of barriers. Since the reaction takes place in a
layer of width $O(\ve)$, the mesh has a fine band, of size $h_f$, around the front, and
is coarse, of size $h_c$, elsewhere; the band moves with the front.

For radial data (in dimension $1$, $2$ or $3$), a perturbation $w_0(x,t)$ that need not be
constant, and a band that moves with a prescribed velocity, we prove that the discrete
solution $U$ converges to the classical solution $v$ of \eqref{eq:P}, with
\[
|U-v|\le C\Big(\ve\log\frac1\ve+\frac{h_f}\ve+h_c+\Delta t+\mu\,\min\Big\{1,\frac{\Delta t}\ve\Big\}\Big),
\]
where $\mu$ is the largest difference between the velocity of the band and that of the
free boundary (Theorem~\ref{cv:thm:main}). Thus it is enough to resolve the layer,
$h_f\ll\ve$. The time step, instead, depends on the mesh: if the mesh does not move it
enters as $\Delta t/\ve$, and indeed the computed front lags, by about
$0.23\,(\Delta t/\ve)T$ for a travelling wave, while on a mesh that follows the front it
enters as $\Delta t$. For a travelling wave on a mesh that moves with it the time step
does not appear at all (Theorem~\ref{an:thm:ale}). In the code the velocity of the band
is computed, not prescribed, so that $\mu$ is only known a posteriori; this and the
other differences between the theorem and the method, as well as the terms of the bound
that are not sharp, are discussed in Section~\ref{an:sec:impl}.

Moving and adaptive meshes for thin flame fronts are classical (see
\cite{HuangRussell2011} and the references therein), but they work at a fixed
activation energy, and we are not aware of error estimates for them that are uniform in
$\ve$, or that measure the distance to the limit problem. The closest results that we
know are those for the Allen--Cahn and double obstacle approximations of mean curvature
flow \cite{NochettoVerdi1996,NochettoVerdi1997}, where fully discrete schemes are shown
to converge to the limit interface, with rates. Theorem~\ref{cv:thm:main} is a result of
this kind for the present problem.

Both results use exact solutions, which we collect in Section~\ref{sec:exact}: the
travelling waves of \eqref{eq:Peps}, for which we also obtain the first-order correction
of the flame speed (Proposition~\ref{prop:speed}), and the radial self-similar
solutions of \eqref{eq:P}. In Section~\ref{sec:experiments} we compare \eqref{eq:Peps}
with the limit problem, also for a perturbation $w_0$ that is not constant, and test the
scheme; the experiments on travelling waves and on the cost of the method are summarized in
Section~\ref{sec:exp:tw} and reported in full in a supplementary document available
with the code.

\paragraph{Code availability} The Python code for all the schemes of this paper, and
the scripts that produce every figure and table, are available at
\url{https://github.com/jfbonder/combustion}.

\section{Exact solutions}\label{sec:exact}

In this section $f$ satisfies
\begin{itemize}
\item[(F)] $f$ is Lipschitz (constant $L_f$) and nonincreasing, $f>0$ on $(-\infty,1)$ and
$f=0$ on $[1,\infty)$.
\end{itemize}
Monotonicity of $f$ is not needed in \cite{FernandezBonderWolanski2000}. It means that
the reaction rate increases with the temperature (recall that $u$ is minus the
temperature); it holds for $f(s)=(1-s)^+$ and for Arrhenius rates with a cut-off, such
as $f(s)=(e^{-s}-e^{-1})^+$ modified to be bounded for $s\to-\infty$, and it is what
makes the scheme of Section~\ref{sec:method} monotone. The rate $e^{-s}$ itself does not
vanish on $[1,\infty)$ and does not satisfy (F).

\subsection{Travelling waves}\label{sec:tw}
For \eqref{eq:P} in one dimension with $w_0>-1$ constant,
$u(x,t)=u_0(x-c_0t)$ with
\[
u_0(z)=A\big(1-e^{-c_0z}\big)^+,\qquad c_0=\frac{\sqrt{2M}}A ,
\]
is a solution moving into the fresh region $\{u>0\}$, where $A=\lim_{x\to\infty}u>0$ is
the difference between the reference temperature and the temperature of the fresh gas.
The formula for $c_0$ is that of Zeldovich and Frank-Kamenetskii
\cite{ZeldovichFrankKamenetskii1938}.

For \eqref{eq:Peps} with $w^\ve\equiv\ve w_0$ the equation is
$u_t=u_{xx}-g_\ve(u)$ with $g_\ve(s)=(s+\ve w_0)f_\ve(s)$, and a travelling wave
$u=\Phi(x-ct)$ satisfies
\begin{equation}\label{an:eq:tw}
\Phi''+c\Phi'=g_\ve(\Phi),\qquad \Phi(-\infty)=-\ve w_0,\quad\Phi(+\infty)=A .
\end{equation}
We normalize $\Phi(0)=\ve$: throughout the paper the \emph{front} of a solution of
\eqref{eq:Peps} is the level set $\{u=\ve\}$, the cold edge of the reaction zone. This is
a convention, and the distance from this level set to the free boundary of the limit,
which is of order $\ve$, depends on it (Remark~\ref{q:rem:rate}). Ahead of the reaction
zone there is no reaction,
so $\Phi(z)=A-(A-\ve)e^{-cz}$ for $z\ge0$. For $z\le0$ we use the inner variables
$\Phi(z)=\ve\varphi(z/\ve)$, $\zeta=z/\ve$, $\sigma=\ve c$:
\begin{equation}\label{eq:tw}
\varphi''+\sigma\varphi'=(\varphi+w_0)f(\varphi)\ \text{ on }(-\infty,0),\qquad
\varphi(-\infty)=-w_0,\quad\varphi(0)=1,\quad \varphi'(0)=c(A-\ve),
\end{equation}
the last condition because $\Phi\in C^1$. Let
\[
B(\varphi)=\int_{-w_0}^{\varphi}(s+w_0)f(s)\,ds\quad(\text{so }B(1)=M),\qquad
\alpha_0=\sqrt{f(-w_0)},\qquad I=\int_{-w_0}^1\sqrt{2B(\varphi)}\,d\varphi .
\]

\begin{lemma}[The travelling wave of \eqref{eq:Peps}]\label{an:lem:profile}
Assume (F), $w_0>-1$ and $A>0$. There are $\ve_0,K>0$ depending only on $f,w_0,A$,
and for every $Z\ge1$ a constant $\theta_Z>0$ depending only on $f,w_0,A,Z$, such
that for $0<\ve\le\ve_0$ problem \eqref{an:eq:tw} has a solution $(\Phi_\ve,c_\ve)$,
unique up to translations, and, with $\Phi=\Phi_\ve$, $c=c_\ve$:
\begin{enumerate}
\item $|c_\ve-c_0|\le K\ve$; in particular $c_0/2\le c_\ve\le 2c_0$;
\item $0<\Phi'\le K$ on $\R$, and $\Phi'\ge\theta_Z$ on $[-Z\ve,2\ve]$;
\item $\Phi\in C^{2,1}(\R)$, $|\Phi''(z)|\le K\ve^{-1}\Phi'(z)$ for all $z$ and
$|\Phi'''(z)|\le K\ve^{-2}\Phi'(z)$ for a.e.\ $z$; if in addition $f\in C^{1,1}(\R)$, then
$\Phi\in C^{3,1}(\R)$ and $|\Phi''''(z)|\le K\ve^{-3}\Phi'(z)$ a.e. Consequently
$\Phi'(z_1)\le e^{K|z_1-z_2|/\ve}\Phi'(z_2)$ for all $z_1,z_2$. On $(0,\infty)$,
$\Phi^{(k)}=(-c)^{k-1}\Phi'$ for $k\ge1$;
\item for $z\le0$, $0<\Phi(z)+\ve w_0\le K\ve\,e^{\alpha_0 z/(2\ve)}$;
\item $\sup_z|\Phi_\ve(z)-u_0(z)|\le K\ve$.
\end{enumerate}
Moreover, with $p=\varphi'$ regarded as a function of $\varphi\in(-w_0,1]$,
\begin{equation}\label{eq:energyphi}
\tfrac12p(\varphi)^2=B(\varphi)-\sigma J(\varphi),\qquad
J(\varphi)=\int_{-w_0}^{\varphi}p(\phi)\,d\phi ,\qquad 0\le J\le I .
\end{equation}
\end{lemma}

\begin{proof}
\emph{Existence.} With $\tilde u(\tilde x)=(A-\Phi(-\tilde x))/(A+\ve w_0)$,
\eqref{an:eq:tw} becomes $-\tilde u''+c\tilde u'=\tilde g(\tilde u)$,
$\tilde u(-\infty)=0$, $\tilde u(+\infty)=1$, with a Lipschitz
$\tilde g\ge0$ that vanishes on $[0,(A-\ve)/(A+\ve w_0)]$ and at $1$ and is positive in
between. This is the scalar problem with an ignition-type nonlinearity of
\cite[Thm.~3.1]{BerestyckiNicolaenkoScheurer1985}, which has a solution $(\tilde u,c)$ with
$\tilde u'>0$, unique up to translations. In particular $\Phi'>0$, and $p=\varphi'>0$ on
$(-\infty,0]$.

\emph{Energy and monotonicity.} Let $\mathcal H=\frac12p^2-B(\varphi)$. By \eqref{eq:tw},
$\frac{d\mathcal H}{d\zeta}=p\big(\varphi''-(\varphi+w_0)f(\varphi)\big)=-\sigma p^2\le0$, and $\mathcal H\to0$ as
$\zeta\to-\infty$. Since $p>0$, $\varphi$ is a valid independent variable,
$\int_{-\infty}^{\zeta}p^2\,d\zeta'=\int_{-w_0}^{\varphi(\zeta)}p\,d\phi$, and integrating we get the
first identity in \eqref{eq:energyphi}. It gives $p\le\sqrt{2B}$, hence $0\le J\le I$.

\emph{Item (1).} At $\varphi=1$, \eqref{eq:energyphi} and $p(0)=c(A-\ve)$ give
$2M-2\sigma I\le c^2(A-\ve)^2\le 2M$. With $\sigma=\ve c$ this yields
$c\le c_0A/(A-\ve)$ and $c(A-\ve)\ge\sqrt{2M}-2\ve cI/\sqrt{2M}$, so
$|c-c_0|\le K\ve$ for $\ve\le\ve_0$.

\emph{Lower bounds for $p$.} Since $f$ is nonincreasing, $f\le\alpha_0^2$ on
$[-w_0,\infty)$, so $B(\varphi)\le\frac12\alpha_0^2(\varphi+w_0)^2$,
$p\le\alpha_0(\varphi+w_0)$ and $J(\varphi)\le\frac12\alpha_0(\varphi+w_0)^2$; also
$B(\varphi)\ge\frac12f(\varphi)(\varphi+w_0)^2$. Fix $\varphi_*\in(-w_0,1)$ with
$f(\varphi_*)\ge\alpha_0^2/2$. By \eqref{eq:energyphi}, if $\sigma\le\alpha_0/4$,
\begin{equation}\label{eq:plow}
\begin{aligned}
p^2&\ge\big(\tfrac12\alpha_0^2-\sigma\alpha_0\big)(\varphi+w_0)^2\ge\tfrac14\alpha_0^2(\varphi+w_0)^2
&&\text{for }\varphi\le\varphi_*,\\
p^2&\ge2B(\varphi_*)-\sigma\alpha_0(1+w_0)^2\ge B(\varphi_*)&&\text{for }\varphi_*\le\varphi\le1,
\end{aligned}
\end{equation}
the last inequality for $\sigma$ small.

\emph{Item (4).} By \eqref{eq:plow}, $\frac{d}{d\zeta}\log(\varphi+w_0)\ge\frac{\alpha_0}2$
while $\varphi\le\varphi_*$, and the point $\zeta_*$ where $\varphi=\varphi_*$ satisfies
$|\zeta_*|\le(1-\varphi_*)/\sqrt{B(\varphi_*)}=:Z_*$. Hence
$\varphi+w_0\le(1+w_0)e^{\frac{\alpha_0}{2}(\zeta+Z_*)}$ for $\zeta\le0$, which is (4) with
$K=(1+w_0)e^{\alpha_0Z_*/2}$.

\emph{Item (2).} $p\le\sqrt{2M}$, and $\Phi'=(A-\ve)ce^{-cz}\le Ac$ for $z\ge0$. On
$[0,2\ve]$, $\Phi'\ge(A-\ve)ce^{-2c\ve}\ge Ac_0/4$. On $[-Z\ve,0]$, that is, for
$\zeta\in[-Z,0]$:
where $\varphi\ge\varphi_*$, $p\ge\sqrt{B(\varphi_*)}$ by \eqref{eq:plow}; where
$\varphi\le\varphi_*$, $p\ge\frac{\alpha_0}2(\varphi+w_0)\ge\frac{\alpha_0}2(\varphi_*+w_0)e^{-\alpha_0Z}$,
because $\frac{d}{d\zeta}\log(\varphi+w_0)=p/(\varphi+w_0)\le\alpha_0$.

\emph{Item (3).} Let $g(s)=(s+w_0)f(s)$, which is Lipschitz, with
$0\le g(\varphi)\le\alpha_0^2(\varphi+w_0)$. By \eqref{eq:plow}, $g(\varphi)\le Cp$ on
$(-w_0,1]$, so $|\varphi''|=|g(\varphi)-\sigma p|\le Cp$. Differentiating,
$\varphi'''=g'(\varphi)p-\sigma\varphi''$ a.e., so $|\varphi'''|\le Cp$; and if
$f\in C^{1,1}$, $\varphi''''=g''(\varphi)p^2+g'(\varphi)\varphi''-\sigma\varphi'''$ a.e., so
$|\varphi''''|\le Cp$. Since $\Phi^{(k)}(z)=\ve^{1-k}\varphi^{(k)}(z/\ve)$ and $\Phi'=p$, this
gives the bounds for $z<0$; for $z>0$ they follow from the explicit formula. If
$f\in C^{1,1}(\R)$ then $f'(1)=0$, $g_\ve\in C^{1,1}$ and $\Phi\in C^{3,1}(\R)$; in
general $\Phi\in C^{2,1}(\R)$. Finally $|(\log\Phi')'|\le K/\ve$ gives the Harnack
inequality.

\emph{Item (5).} For $z\ge0$,
$|\Phi(z)-u_0(z)|\le A|e^{-c_0z}-e^{-cz}|+\ve e^{-cz}\le A\,|c-c_0|\,\sup_{z\ge0}ze^{-\min(c,c_0)z}+\ve
\le\frac{2A}{ec_0}K\ve+\ve$. For $z\le0$, $u_0=0$ and $-\ve w_0\le\Phi\le\ve$.
\end{proof}

The identity $\frac{d\mathcal H}{d\zeta}=-\sigma p^2$ also gives a numerical method for $c_\ve$:
integrating \eqref{eq:tw} backwards from $\zeta=0$, if $\sigma$ is too large $\mathcal H$ remains
positive and $\varphi$ crosses $-w_0$, and if $\sigma$ is too small $\varphi'$ vanishes
before. Bisection on $\sigma$ gives $c_\ve=\sigma/\ve$ to machine precision,
independently of any PDE solver. We use it as a reference in
Section~\ref{sec:experiments}.

\begin{proposition}[First-order correction of the flame speed]\label{prop:speed}
Assume (F), $w_0>-1$ and $A>0$. There are $\ve_0,C>0$ such that for $0<\ve\le\ve_0$
\[
\Big|c_\ve-c_0\Big[1+\frac\ve A\Big(1-\frac{I}{\sqrt{2M}}\Big)\Big]\Big|\le C\ve^2 ,
\qquad
1-\frac{I}{\sqrt{2M}}=-w_0+\int_{-w_0}^1\Big(1-\sqrt{\frac{B(\varphi)}{M}}\Big)\,d\varphi .
\]
For $f(s)=(1-s)^+$, $I/\sqrt{2M}=\kappa_*(1+w_0)$ with
$\kappa_*=\frac{3\sqrt3-2}{5}\approx0.6392$.
\end{proposition}

\begin{proof}
We use the notation of Lemma~\ref{an:lem:profile}, with $c=c_\ve$ and $\sigma=\ve c\le2c_0\ve$.

\emph{Step 1: $J(1)$ is close to $I$.} Let $t=\varphi+w_0\in(0,1+w_0]$. Since
$f\ge m>0$ on $[-w_0,\frac{1-w_0}2]$ and $B$ is increasing, $B(\varphi)\ge b_0t^2$ with
$b_0>0$ depending only on $f,w_0$; and $J(\varphi)\le\frac12\alpha_0t^2$, as in the proof of
Lemma~\ref{an:lem:profile}. By \eqref{eq:energyphi}, $0\le2B-p^2=2\sigma J$, so
\[
0\le\sqrt{2B}-p=\frac{2\sigma J}{\sqrt{2B}+p}\le\frac{\sigma\alpha_0t^2}{\sqrt{2b_0}\,t}
=\frac{\alpha_0}{\sqrt{2b_0}}\,\sigma\,t ,
\qquad
0\le I-J(1)=\int_{-w_0}^1\big(\sqrt{2B}-p\big)\,d\varphi\le C\sigma ,
\]
that is, $J(1)=I-r_\ve$ with $0\le r_\ve\le C\ve$.

\emph{Step 2: the speed.} At $\varphi=1$, \eqref{eq:energyphi} reads
$\frac12c^2(A-\ve)^2=M-\ve cJ(1)$. With $\tau=\ve cJ(1)/M\in[0,C\ve]$ and
$\sqrt{1-\tau}=1-\frac\tau2+O(\tau^2)$,
\[
c(A-\ve)=\sqrt{2M}-\frac{\ve cJ(1)}{\sqrt{2M}}+O(\ve^2),
\qquad\text{i.e.}\qquad
c\Big(A-\ve+\frac{\ve J(1)}{\sqrt{2M}}\Big)=\sqrt{2M}+O(\ve^2),
\]
so that $c=c_0\big[1+\frac\ve A\big(1-\frac{J(1)}{\sqrt{2M}}\big)\big]+O(\ve^2)$. Replacing
$J(1)$ by $I$ costs $c_0\ve r_\ve/(A\sqrt{2M})=O(\ve^2)$. The second identity in the
statement holds because $I/\sqrt{2M}=\int_{-w_0}^1\sqrt{B/M}\,d\varphi$ and the interval
has length $1+w_0$.

\emph{Step 3: $f(s)=(1-s)^+$.} Let $a=1+w_0$ and $t=\varphi+w_0$. Then $f(\varphi)=a-t$ for
$t\le a$, $B=\frac a2t^2-\frac13t^3$, $M=B|_{t=a}=\frac{a^3}6$ and $\sqrt{2M}=a^{3/2}/\sqrt3$.
With $t=a\tau$ and then $v=1-\frac23\tau$,
\[
I=\int_0^at\sqrt{a-\tfrac23t}\,dt=a^{5/2}\int_0^1\tau\sqrt{1-\tfrac23\tau}\,d\tau
=\frac94a^{5/2}\int_{1/3}^1\big(v^{1/2}-v^{3/2}\big)\,dv
=\frac94a^{5/2}\Big[\frac4{15}-\frac{8}{5\cdot3^{5/2}}\Big],
\]
and $\frac{I}{\sqrt{2M}}=\frac{9\sqrt3}{4}a\Big[\frac4{15}-\frac{8}{5\cdot3^{5/2}}\Big]
=a\Big(\frac{3\sqrt3}5-\frac25\Big)=\kappa_*a$.
\end{proof}

Berestycki, Nicolaenko and Scheurer proved that $c_\ve\to c_0$ and that the
profiles converge, by compactness and without a rate
\cite[Thm.~5.2]{BerestyckiNicolaenkoScheurer1985} (for $w_0=0$). Two-term expansions of
the flame speed in the inverse activation energy are classical in combustion
theory \cite{BushFendell1970,BuckmasterLudford1982}, and a rigorous one was recently
obtained for a related scalar equation by geometric desingularization
\cite{JelbartKristiansenSzmolyan2025}; what Proposition~\ref{prop:speed}
adds, for the present model, is a rigorous remainder and the explicit dependence on the
$O(\ve)$ perturbation $w_0$. The two terms of the correction have different origins. The
term $-w_0$ is a matter of normalization: the temperature rise across the flame is
$A+\ve w_0$, not $A$, and the formula of Zeldovich and Frank-Kamenetskii with this value is
$\sqrt{2M}/(A+\ve w_0)=c_0(1-\ve w_0/A)+O(\ve^2)$. The integral is the correction due to the
finite width of the reaction zone; it is positive, and for $f(s)=(1-s)^+$ it equals
$(1-\kappa_*)(1+w_0)$. Hence, for this $f$, the first-order correction is positive for
$w_0<w_0^*=\frac1{\kappa_*}-1\approx0.564$ and negative for $w_0>w_0^*$; the change of sign
only reflects the choice of $A$ as the reference temperature rise.
For $f(s)=(1-s)^+$, $A=1$ and $w_0=-0.5$, $0$, $1$ we have compared the prediction with
the speed obtained by shooting, for $\ve$ from $0.1$ to $10^{-4}$: the differences
between $(c_\ve-c_0)/\ve$ and its limit value are $0.10\,\ve$, $0.12\,\ve$ and
$0.65\,\ve$ respectively, in agreement with the remainder $O(\ve^2)$ of the proposition. The values $\ve=0.1$ and
$0.05$ correspond to Zeldovich numbers of $10$ and $20$, which are realistic; there $c_0$
differs from $c_\ve$ by $2.5$ to $7\%$ and by $1.3$ to $3.4\%$ respectively, and the
first-order correction reduces these differences to at most $0.5\%$ and $0.12\%$.

\subsection{The radial self-similar solution}\label{sec:selfsimilar}
Caffarelli and V\'azquez \cite{CaffarelliVazquez1995} showed that, for $w_0$ constant,
\eqref{eq:P} has solutions
\[
V(x,t)=(T-t)^{1/2}\,h\big(|x|(T-t)^{-1/2}\big)
\]
with compact support shrinking to a point at $t=T$: a ball of fresh gas that is consumed
by the flame, whose front is concave towards the fresh gas and moves with speed
proportional to $(T-t)^{-1/2}$. (As a model of a flame, \eqref{eq:P} is meaningful while
the radius is large compared with $\ve$.) Following
\cite{GalaktionovHulshofVazquez1997,Petrosyan2002}, we call the instant at
which the fresh region disappears the \emph{extinction time}; note that it is the fresh
pocket that is extinguished, not the flame in the sense of combustion theory. These
solutions are used in
\cite[Prop.~4.3]{FernandezBonderWolanski2004} to show that limit solutions vanish in
finite time. Substituting, $V$ is caloric in its positivity set if and only if
\begin{equation}\label{eq:h}
h''+\Big(\frac{N-1}r-\frac r2\Big)h'+\frac12h=0.
\end{equation}
The change of
variables $z=r^2/4$ turns \eqref{eq:h} into Kummer's equation
$zh_{zz}+(\frac N2-z)h_z+\frac12h=0$, so that a solution with $h'(0)=0$ is
\[
h(r)=\beta\,{}_1F_1\big(-\tfrac12;\tfrac N2;\tfrac{r^2}4\big),
\qquad
{}_1F_1(a;b;z)=\sum_{k=0}^\infty\frac{(a)_k}{(b)_k}\frac{z^k}{k!},
\]
where ${}_1F_1$ is Kummer's confluent hypergeometric function,
$(a)_k=a(a+1)\cdots(a+k-1)$ is the Pochhammer symbol, and $\beta$ is a constant. Since $-1<-\frac12<0$, ${}_1F_1(-\frac12;\frac N2;\cdot)$ has exactly one
positive zero $z_N$ \cite[\S13.9]{DLMF}. Then $h>0$ on $[0,R)$ with $R=2\sqrt{z_N}$, and
$\beta>0$ is fixed by the free boundary condition $h'(R)=-\sqrt{2M}$. Numerically
$R=1.8483$, $2.5139$ and $3.0040$ for $N=1,2,3$. The free boundary is the sphere of
radius $R\sqrt{T-t}$.

\section{Rate of convergence of the singular limit}\label{sec:rate}

In this section we make the approximation theorems of
\cite[Theorems~4.1--4.2]{FernandezBonderWolanski2004} quantitative. There, a
strict classical supersolution of \eqref{eq:P} is shown to be the uniform limit
of supersolutions of \eqref{eq:Peps}, by solving \eqref{eq:Peps} in a layer
and checking the jump of the gradients by a blow-up argument by contradiction.
We replace both steps by an explicit barrier: an inner profile, written in the
variable $(u+w)/\ve$ (the rescaled distance to the burnt state), composed with
an extension of the supersolution whose gradient is monotone in the normal
direction. The approximation error is $O(\ve)$ with constants that do not
depend on the strictness of the super- or subsolution, and this yields a rate
for the convergence of $u^\ve$ to a classical solution
(Theorem~\ref{q:thm:rate}).

\emph{Relation with previous constructions.} For $w^\ve=0$ the uniqueness of the limit
for starshaped data is due to Petrosyan \cite{Petrosyan2002}. There
\cite[Proposition~2.5]{Petrosyan2002} the solutions of \eqref{eq:Peps} are compared with
an explicit supersolution, the composition of the strict supersolution $\tilde u$ with the
solution of the one-dimensional equation in which the reaction term is cut off below a
level $a\ve$ and multiplied by a constant, so that the profile is constant on the burnt
side and no extension of $\tilde u$ is needed. The statement is qualitative, but the
construction can be made quantitative: for a regular strict supersolution with
strictness $\delta$ one can check that it is a supersolution of \eqref{eq:Peps} when
$\ve\log(1/\delta)\le c\,\delta$, at distance $O(\ve\log(1/\delta))$ from $\tilde u$, because the
truncated profile has a length of order $\log(1/\delta)$; with the dilations used in the
proof of Theorem~\ref{q:thm:rate} this would give the upper bound $u^\ve\le v+C\ve\log(1/\ve)$ in the case
$w^\ve=0$. The analogous statement for strict subsolutions is asserted in
\cite[Remark~2.6]{Petrosyan2002}, and the same count gives
$u^\ve\ge v-C\ve\log(1/\ve)$. Moreover, if $w^\ve\equiv\ve w_0$ with $w_0$ constant, then
$Y^\ve=u^\ve+\ve w_0$ solves the scalar equation $Y_t=\Delta Y-\ve^{-1}\beta(Y/\ve)$ with
$\beta(s)=sf(s-w_0)$, $\int\beta=M(w_0)$, which is the problem of
\cite{CaffarelliVazquez1995,Petrosyan2002} with another constant in the free boundary
condition. Hence, for constant $w_0$, a two-sided bound of order $\ve\log(1/\ve)$ can be
extracted from those barriers, with less regularity of $\tilde u$ than we require (no
extension is needed), and what is new in Theorem~\ref{q:thm:rate} is the upper bound
without logarithm. When $W^\ve=w^\ve/\ve$ is not constant there is no such reduction, and, as
observed in \cite{FernandezBonderWolanski2004}, the construction does not extend without
further conditions on $f$.
In \cite{FernandezBonderWolanski2004}, following
\cite{LedermanVazquezWolanski2001}, the approximating function $u^\ve$ equals the strict
supersolution $\tilde u$ where $\tilde u\ge(1-\theta)\ve$ and solves \eqref{eq:Peps} in the
set $\{\tilde u<(1-\theta)\ve\}$, which contains the burnt region, so that $|u^\ve-\tilde u|\le C\ve$ is immediate; that
$u^\ve$ is a supersolution, that is, that the jump of its gradient across the boundary of
that set has the right sign, is proved for $\ve<\ve_0$ by a blow-up argument by
contradiction, which gives no information on how $\ve_0$ depends on the strictness
$\delta$. What is new here is that the explicit barriers of
Propositions~\ref{q:prop:super} and \ref{q:prop:sub}, which keep the whole profile and
are composed with an extension of $\tilde u$ to the burnt region, are super- and
subsolutions as soon
as $\ve\le c_*\delta$, with a threshold that is linear in $\delta$ and an error that does
not depend on $\delta$, also when $W^\ve$ is not constant; this is what
produces the rate, without logarithm in the upper bound. With respect to the barriers used
for Allen--Cahn type limits \cite{Chen1992,DeMottoniSchatzman1995,AlfaroHilhorstMatano2008},
where the profile is a function of the signed distance to the interface and the
strictness is provided by a shift of the interface that grows in time, there are three
differences: the inner profile is coupled to the gradient of the outer solution, through
the free boundary condition; it is not adjusted to the local value of that gradient but
to an energy level $G_\pm(W)$ at a distance of order $\delta$ from $2M(W)$; and the
strictness is provided by the dilations of starshaped data, as in
\cite{FernandezBonderWolanski2004}. The barriers are coarse: they give the order $\ve$,
but they do not capture the term of order $\ve$ of the solution, as the expansions of
\cite{DeMottoniSchatzman1995} do in their setting.

\subsection{Setting}
Let $Q_T=\R^N\times(0,T)$. For $W>-1$ and $t\in\R$ let
\[
b(t;W)=t\,f(t-W),\qquad \tilde B(\Theta;W)=\int_0^\Theta b(t;W)\,dt,\qquad
M(W)=\tilde B(1+W;W) .
\]
Then $M(W)=\int_{-W}^1(s+W)f(s)\,ds$, as in \eqref{eq:P}.
If $u=-w+\ve\Theta$, the reaction term of \eqref{eq:Peps} is
$(u+w)f_\ve(u)=\ve^{-1}b(\Theta;w/\ve)$; $\Theta=0$ is the burnt state, and
$b(\cdot;W)=0$ on $[1+W,\infty)$. (In the notation of Section~\ref{sec:tw},
$\tilde B(\Theta;w_0)=B(\Theta-w_0)$.)

\textbf{Assumptions.}
\begin{itemize}
\item[(F$_2$)] $f$ satisfies (F), and $f\in C^2((-\infty,1])$ (one-sided derivatives at $1$).
\item[(W)] $w^\ve$ is a bounded solution of the heat equation in $Q_T$,
$W^\ve:=w^\ve/\ve$ satisfies $\|\nabla W^\ve\|_{L^\infty(Q_T)}\le C_W$ and
$\|W^\ve-w_0\|_{L^\infty(Q_T)}\le\eta_\ve$, and $-1+\delta_1\le w_0\le \bar W$.
\end{itemize}
Assumption (F$_2$) holds for $f(s)=(1-s)^+$. Assumption (W) is a restriction on
the initial data of $w^\ve$; it holds with $\eta_\ve=0$ when $w^\ve(\cdot,0)=\ve w_0$ with
$w_0$ constant. We write $W=W^\ve(x,t)$ and
$\mathcal W=[-1+\delta_1/2,\bar W+1]$, and we assume $\eta_\ve\le\min\{\delta_1/2,1\}$, so that
$W\in\mathcal W$. Note that
$M'(W)=\int_{-W}^1f>0$, so $M$ is increasing and Lipschitz on $\mathcal W$. If
$W^\ve\equiv w_0$ is constant, then the term $\mathcal P$ in \eqref{q:eq:basic} below
vanishes, the derivatives of the profiles with respect to $W$ are not needed, and
assumption (F) is enough in this section.

\textbf{Comparison.} We say that $V$ is a \emph{strong supersolution} of
\eqref{eq:Peps} in an open set $D\subset Q_T$ if $V$ is bounded and continuous in
$\overline D$, locally $C^{1,1}$ in $x$ and Lipschitz in $t$ in $D$ (assumption (W) bounds
$\nabla W^\ve$ up to $t=0$, but $D^2W^\ve$ and $\partial_tW^\ve$ only for $t>0$), and
$\Delta V-V_t-(V+w^\ve)f_\ve(V)\le0$ a.e.\ in $D$; strong
subsolutions satisfy the reverse inequality. The initial data $u_0^\ve$ are continuous
and bounded, and solutions of \eqref{eq:Peps} are bounded classical solutions, continuous
up to $t=0$. For fixed $\ve$ the nonlinearity is
Lipschitz in $V$, with constant of order $\ve^{-2}$, so the maximum principle for strong
solutions in (not necessarily cylindrical) bounded domains
\cite[Chapters~II and VII]{Lieberman1996}, applied when $D$ is unbounded to the difference
multiplied by the weight $e^{-\beta t}/\cosh(|x|/R)$ with $R\to\infty$ (the functions being
bounded, this is the Phragm\'en--Lindel\"of argument), gives: if $V$ is a strong
supersolution in $D$,
$u$ is a bounded classical solution, and $V\ge u$ on the parabolic boundary of
$D$ (that is, on $\overline D\cap\{t=0\}$ and on $\partial D\cap Q_T$), then $V\ge u$ in
$D$; and similarly for subsolutions. When $V$ is the composition of a $C^{1,1}$ profile
with a $C^{1,1}$ function, the chain rule holds a.e.\ provided the set where the argument
takes the values at which the second derivative of the profile jumps has measure zero;
we check this in each case.

\textbf{The basic identity.} Let $\tilde v$ be $C^{1,1}$ in $x$ and Lipschitz in
$t$, and let $\Theta=\Theta(s;W)$ be $C^{1,1}$. For
\[
V=-w^\ve+\ve\,\Theta\Big(\frac{\tilde v}\ve;W\Big),\qquad s=\frac{\tilde v}{\ve},
\]
using that $w^\ve$ and $W$ are caloric, we have a.e.
\begin{equation}\label{q:eq:basic}
\begin{gathered}
\Delta V-V_t-(V+w^\ve)f_\ve(V)=\Theta_s(\Delta\tilde v-\tilde v_t)
+\frac1\ve\big(\Theta_{ss}|\nabla\tilde v|^2-b(\Theta;W)\big)+\mathcal P,\\
\mathcal P=2\Theta_{sW}\nabla\tilde v\cdot\nabla W+\ve\Theta_{WW}|\nabla W|^2 .
\end{gathered}
\end{equation}

\subsection{Inner profiles}
Fix $0<\gamma\le\gamma_0$ and let $G_\pm(W)=\sqrt{2M(W)\pm2\gamma}$. In this
subsection $c,C,\lambda_0$ denote positive constants that depend only on
$f,\delta_1,\bar W$ (\emph{not} on $\gamma$ or on the damping $\kappa$ below), and
$\gamma_0$ is small in terms of them. We let $2a_0=f(1-\delta_1/4)>0$, so that
$f\ge2a_0$ on $(-\infty,1-\delta_1/4]$.

\begin{lemma}[Inner profiles]\label{q:lem:inner}
Assume (F$_2$).

\textbf{(a) Supersolution profile.} For $W\in\mathcal W$ and $0\le\kappa\le c\gamma$
there are $s_m=s_m(W)\in[\frac12,1)$, with $1-s_m\le\delta_1/4$ and $|s_m'(W)|\le C$, and
$\Theta=\Theta^+(s;W)$ such that $\Theta(s)=s+W$ for $s\ge s_m$,
\[
G_-(W)^2\Theta_{ss}+\kappa\Theta_s=b(\Theta;W)\quad\text{for }s<s_m,\qquad\Theta(-\infty)=0,
\]
and
\begin{enumerate}
\item[(i)] $0<\Theta_s\le1$ and $0\le\Theta_{ss}\le C\Theta_s$; $0<\Theta\le Ce^{\lambda_0s}$ for
$s\le0$ and $c\,\Theta\le\Theta_s\le C\Theta$ for $s\le s_m$; $\Theta_s\ge c\,e^{-CS}$ on $[-S,s_m]$ for every $S\ge0$; and
$b(\Theta;W)\ge\gamma/4$ where $\frac{1+W}2\le\Theta\le\Theta(s_m)$;
\item[(ii)] $\Theta$ is $C^{1,1}$ in $(s,W)$, of class $C^2$ on $\{s\ne s_m(W)\}$, and
\[
|\Theta_W|\le C(1+|s|)\Theta_s,\qquad
|\Theta_{sW}|\le C(1+|s|)\Theta_{ss}+C\Theta_s,\qquad
|\Theta_{WW}|\le C(1+|s|)^2\Theta_s .
\]
\end{enumerate}
\textbf{(b) Subsolution profile.} For $W\in\mathcal W$, $0\le\kappa\le c\gamma$ and
$0<\bar q\le1$ there are $s_*=s_*(W)<s_c(W)<1$, Lipschitz in $W$ with constant
$C\log(1/\gamma)$, and
$\Theta=\Theta^-(s;W)$, defined for $s\ge s_*$, with $\Theta(s)=s+W$ for $s\ge2$,
\[
G_+(W)^2\Theta_{ss}-\kappa\Theta_s=b(\Theta;W)\quad\text{for }s_*<s<2,\qquad
\Theta(s_*)=-\bar q,\quad\Theta(s_c)=0,
\]
and
\begin{enumerate}
\item[(i)] $c\sqrt{\gamma+\Theta^2}\le\Theta_s\le C$ on $[s_*,2]$; $|\Theta_{ss}|\le C\Theta_s$;
$\Theta_{ss}\ge0$ where $\Theta\ge0$, and $|\Theta_{ss}|\le C(\bar q+\kappa\Theta_s)$ where $\Theta<0$;
$0\le\Theta\le(2+W)e^{c(s-2)}$ on $[s_c,2]$; and $s_*\ge-C\log(1/\gamma)$;
\item[(ii)] $\Theta$ is $C^{1,1}$ in $(s,W)$ on $\{s\ge s_*(W)\}$, of class $C^2$ on
$\{s\ne2,\ \Theta\ne1+W\}$, and
\[
|\Theta_W|\le C(1+|s|)\Theta_s,\qquad
|\Theta_{sW}|\le C(1+|s|)|\Theta_{ss}|+C\Theta_s,\qquad
|\Theta_{WW}|\le C(1+|s|)^2\Theta_s .
\]
\end{enumerate}
\end{lemma}

In (b) the profile is linear only for $s\ge2$, beyond the point $\Theta=1+W$ where the
reaction stops ($b(\Theta;W)=0$ for $\Theta\ge1+W$); this plays no role in this section, and
it is used in Section~\ref{sec:analysis}.

\begin{proof}
Both profiles are increasing in $s$, so we parametrize them by $\Theta$: with
$P(\Theta)=\Theta_s$ we have $\Theta_{ss}=PP_\Theta$, and the equations
$G^2\Theta_{ss}\pm\kappa\Theta_s=b$ become
$\frac{d}{d\Theta}\big(\frac{G^2}2P^2\big)=b(\Theta;W)\mp\kappa P$. We use that
\begin{equation}\label{q:eq:Bquad}
c\Theta^2\le\tilde B(\Theta;W)\le C\Theta^2\ \text{ for }-1\le\Theta\le1+W,\qquad
|\partial_W\tilde B|+|\partial^2_{W}\tilde B|\le C\Theta^2 ,
\end{equation}
which follow from $2a_0\le f(t-W)\le C$ for $-1\le t\le\frac{1+W}2$, from
$\partial_W\tilde B=-\int_0^\Theta tf'(t-W)\,dt$, and from (F$_2$). They also hold for
$1+W\le\Theta\le2+W$, where $\tilde B(\Theta;W)=M(W)$. Note also that
$G_\pm G_\pm'=M'$ is bounded, and that $G_\pm\ge c$.

\smallskip
(a) Let $G=G_-(W)$ and let $P$ be the branch $P>0$ of the unstable manifold of the saddle
$(0,0)$ in the $(\Theta,P)$-plane:
\begin{equation}\label{q:eq:Psup}
\tfrac{G^2}2P^2=\tilde B(\Theta;W)-\kappa J(\Theta),\qquad J(\Theta)=\int_0^\Theta P .
\end{equation}
By \eqref{q:eq:Bquad}, $P\le\sqrt{2\tilde B}/G\le C\Theta$, so $\kappa J\le C\kappa\Theta^2$ and, for
$\kappa$ small,
\begin{equation}\label{q:eq:Plin}
c\,\Theta\le P\le C\Theta\qquad\text{as long as $P$ is defined and }\Theta\le1+W .
\end{equation} Let $F(\Theta)=\tilde B(\Theta;W)-\kappa J(\Theta)-(M(W)-\gamma)$. Then $F(0)<0$ and
$F(1+W)=\gamma-\kappa J(1+W)\ge\gamma/2$ if $\kappa\le c\gamma$. Let $\Theta_m$ be the first
zero of $F$: then $P(\Theta_m)=1$ and $0<P<1$ on $(0,\Theta_m)$. Set $s_m=\Theta_m-W$.
Since
\[
\tfrac\gamma2\le\gamma-\kappa J(\Theta_m)=M-\tilde B(\Theta_m)=\int_{\Theta_m}^{1+W}tf(t-W)\,dt\le\gamma
\]
and $f>0$ on $(-\infty,1)$, we have $1-s_m\to0$ as $\gamma\to0$, uniformly in $W$; we take
$\gamma_0$ so small that $1-s_m\le\min\{\frac12,\frac{\delta_1}4\}$. Then
$\Theta_m\ge\frac{1+W}2$ and, since $f$ is nonincreasing,
$\frac\gamma2\le(1-s_m)(1+W)f(s_m)\le\Theta_mf(s_m)=b(\Theta_m;W)$, and
$b(\Theta;W)\ge\frac12b(\Theta_m;W)$ on $[\frac{1+W}2,\Theta_m]$. Hence, for $\kappa\le\gamma/8$,
\begin{equation}\label{q:eq:top}
b(\Theta;W)\ge\tfrac\gamma4\ \text{ on }[\tfrac{1+W}2,\Theta_m],\qquad
\partial_\Theta F(\Theta_m)=b(\Theta_m;W)-\kappa\ge\tfrac12b(\Theta_m;W)\ge\tfrac\gamma4 .
\end{equation}
Define $s(\Theta)=s_m-\int_\Theta^{\Theta_m}d\theta/P(\theta)$ for $0<\Theta\le\Theta_m$, let
$\Theta(\cdot)$ be its inverse on $(-\infty,s_m]$, and extend it by $\Theta=s+W$ for $s\ge s_m$.
Since $P(\Theta_m)=1$, $\Theta\in C^1$, and $\Theta_{ss}=(b-\kappa\Theta_s)/G^2$ is bounded.

(i) $\Theta_{ss}\ge0$: on $(0,\frac{1+W}2]$, $b\ge2a_0\Theta\ge cP\ge\kappa P$; on
$[\frac{1+W}2,\Theta_m]$, $b\ge\gamma/4\ge\kappa\ge\kappa P$ by \eqref{q:eq:top}. In particular $P$ is
nondecreasing in $\Theta$. Also $\Theta_{ss}\le b/G^2\le C\Theta\le C\Theta_s$ on the first interval
and $\Theta_{ss}\le C\le C\Theta_s$ on the second. By \eqref{q:eq:Plin}, $c\,\Theta\le\Theta_s\le C\Theta$. Since $s_m<1$, this gives
$\Theta(0)\ge\Theta_me^{-C}\ge c$, then $\Theta(-S)\ge c\,e^{-CS}$ and
$\Theta_s\ge c\,\Theta\ge c\,e^{-CS}$ on $[-S,s_m]$; and $\Theta(s)\le\Theta(0)e^{cs}$ for $s\le0$.

(ii) We first check that $P$ is twice continuously differentiable in $W$ for
$0<\Theta\le\Theta_m$, with $P_W=O(\Theta)$ and $P_{WW}=O(\Theta)$ as $\Theta\to0$. Write
$P=\Theta Q$. By \eqref{q:eq:Psup} and \eqref{q:eq:Plin}, on every interval $(0,\Theta_1]$
with $\Theta_1\le1+W$ the function $Q$ is the fixed point of
\begin{gather*}
\mathcal T(Q;W)(\Theta)=\frac{\sqrt2}{G(W)}\Big(\beta(\Theta;W)-\kappa\int_0^1\tau\,Q(\Theta\tau)\,d\tau\Big)^{1/2},\\
\beta(\Theta;W)=\frac{\tilde B(\Theta;W)}{\Theta^2}=\int_0^1\tau f(\Theta\tau-W)\,d\tau,
\end{gather*}
in the set $\{c\le Q\le C\}$ of $C([0,\Theta_1])$, where $\beta\ge c$ by \eqref{q:eq:Bquad}; for
$\kappa$ small $\mathcal T(\cdot\,;W)$ maps this set into itself and is a contraction with
constant $\frac12$. Fix $W_*$ and $\Theta_1<1+W_*$. By (F$_2$), $f''$ is uniformly
continuous on compact subsets of $(-\infty,1]$, so for $W$ near $W_*$ the map
$W\mapsto\beta(\cdot\,;W)$ is of class $C^2$ with values in $C([0,\Theta_1])$, and so is
$(Q,W)\mapsto\mathcal T(Q;W)$. By the uniform contraction principle
\cite[Section~2.2]{ChowHale1982}, $W\mapsto Q(\cdot\,;W)\in C([0,\Theta_1])$ is of class
$C^2$; with the equation for $P$ in $\Theta$, $P$ is $C^2$ jointly in $(\Theta,W)$ for
$0<\Theta<\Theta_1$. Since $\Theta_m<1+W$ and $\Theta_m$ is continuous in $W$ (it is a
transversal zero of $F$, by \eqref{q:eq:top}), this is the claim. We now bound the derivatives with
constants that do not depend on $\gamma$ and $\kappa$.
Differentiating \eqref{q:eq:Psup} with respect to $W$ at fixed $\Theta$,
\[
G^2PP_W=\partial_W\tilde B-M'P^2-\kappa\int_0^\Theta P_W .
\]
By \eqref{q:eq:Bquad}, $|\partial_W\tilde B-M'P^2|\le C(\Theta^2+P^2)$. If
$m(\Theta)=\sup_{(0,\Theta]}|P_W(\theta)|/\theta$, then by \eqref{q:eq:Plin}
$|P_W|\le C\Theta+\frac{C\kappa}{\Theta}\cdot\frac{m(\Theta)\Theta^2}2$ on $(0,\Theta_m]$, i.e.\
$m\le C+C\kappa m$, so $|P_W|\le C\Theta\le CP$ for $\kappa$ small ($m$ is finite by the
first step). Differentiating once more, and using $(G^2)'=2M'$,
\[
G^2PP_{WW}=\partial_W^2\tilde B-M''P^2-4M'PP_W-G^2P_W^2-\kappa\int_0^\Theta P_{WW} ,
\]
where $M''(W)=f(-W)$ is bounded and $|\partial_W^2\tilde B|\le C\Theta^2$ by
\eqref{q:eq:Bquad}. With the bound for $P_W$, the first four terms on the right are bounded
by $C\Theta^2$, and the same argument gives $|P_{WW}|\le C\Theta\le CP$.

Next, $\Theta_m$ is defined by $F=0$, and
$\partial_WF(\Theta_m)=\partial_W\tilde B(\Theta_m;W)-\partial_W\tilde B(1+W;W)-\kappa J_W
=\int_{\Theta_m}^{1+W}tf'(t-W)\,dt-\kappa J_W$, because $M'(W)=\partial_W\tilde B(1+W;W)$. Since
$f'\le0$, the integral is bounded in absolute value by $(1+W)f(s_m)\le2b(\Theta_m;W)$, and
$\kappa|J_W|\le C\kappa\le Cb(\Theta_m;W)$. With \eqref{q:eq:top},
$|\Theta_m'(W)|=|\partial_WF/\partial_\Theta F|\le C$, and so $|s_m'|\le C$.

At fixed $\Theta$, using $P(\Theta_m)=1$,
\[
s_W=s_m'-\Theta_m'+\int_\Theta^{\Theta_m}\frac{P_W}{P^2}\,d\theta
=-1+\int_\Theta^{\Theta_m}\frac{P_W}{P^2}\,d\theta,\qquad\text{so}\qquad
|s_W|\le1+C\int_\Theta^{\Theta_m}\frac{d\theta}P=1+C(s_m-s).
\]
Since $\Theta_W=-s_WP$ (at fixed $s$), this is the first bound. Next,
$\Theta_{sW}=\partial_W[P(\Theta(s;W);W)]=P_\Theta\Theta_W+P_W$ with
$P_\Theta=\Theta_{ss}/\Theta_s\in[0,C]$, which gives the second bound. Finally
$s_{WW}=\Theta_m'P_W(\Theta_m)+\int_\Theta^{\Theta_m}\big(P_{WW}/P^2-2P_W^2/P^3\big)d\theta$, so
$|s_{WW}|\le C(1+s_m-s)$, and
\[
\Theta_{WW}=\frac{P_W}P\Theta_W-s_{WW}P-s_W\big(P_\Theta\Theta_W+P_W\big)
\]
gives the third bound. At $s=s_m$ we have $\Theta_s=1$ and $\Theta_W=-s_WP=1$ from the
left, which are the values from the right; hence $\Theta\in C^1$, with bounded second
derivatives on each side.

\smallskip
(b) Let $G=G_+(W)$. Integrating backwards from $(\Theta,P)=(2+W,1)$ and using
$G^2/2-M=\gamma$ and $\tilde B(2+W;W)=M$,
\begin{equation}\label{q:eq:Psub}
\tfrac{G^2}2P^2=\gamma+\tilde B(\Theta;W)-\kappa\int_\Theta^{2+W}P .
\end{equation}
As long as $P\ge0$ this gives $\frac{G^2}2P^2\le\gamma+\tilde B$, so $P\le1$ for
$0\le\Theta\le2+W$ and $P\le C$ for $-1\le\Theta\le0$; then the last term is at most
$C\kappa\le\gamma/2$ for $\kappa\le c\gamma$, and by \eqref{q:eq:Bquad}
$c(\gamma+\Theta^2)\le P^2$ on $[-1,2+W]$. Hence the orbit reaches $\Theta=-\bar q$; with
$s(\Theta)=2-\int_\Theta^{2+W}d\theta/P$ we set $s_c=s(0)$, $s_*=s(-\bar q)$, and
\[
s_*\ge2-C\int_{-1}^{C}\frac{d\theta}{\sqrt{\gamma+\theta^2}}\ge-C\log\frac1\gamma ;
\]
moreover $s_c\le2-(2+W)<1$, because $P\le1$.
(i) $\Theta_{ss}=(b+\kappa\Theta_s)/G^2$ with $b\ge0$ for $\Theta\ge0$, $|b|\le C\bar q$ for
$-\bar q\le\Theta<0$, and $|b|\le C|\Theta|\le C\Theta_s$; and $\Theta_s\ge c\Theta$ gives the
exponential bound on $[s_c,2]$.

(ii) By \eqref{q:eq:Psub}, in which $\tilde B$ is of class $C^1$ in $W$ with a
derivative that is Lipschitz, and of class $C^2$ away from $\Theta=1+W$, the function $P$
is $C^1$ in $W$, with $P$, $P_\Theta$ and $P_W$ continuous across $\Theta=1+W$, and $C^2$ on each
side. Differentiating \eqref{q:eq:Psub},
$G^2PP_W=\partial_W\tilde B-M'P^2-\kappa\big(1+\int_\Theta^{2+W}P_W\big)$, because $\gamma$
does not depend on $W$. Now $|\partial_W\tilde B-M'P^2|\le CP^2$ and
$\kappa\le c\gamma\le CP^2$, so $|P_W|\le CP+\frac{C\kappa}{P}\int_\Theta^{2+W}|P_W|$ with
$\kappa/P\le C\sqrt\gamma$, and Gronwall's inequality (from $\Theta=2+W$ downwards)
gives $|P_W|\le CP$; similarly $|P_{WW}|\le CP$ a.e.\ ($\partial_W^2\tilde B$ may jump across
$\Theta=1+W$). Here $s_W=-1+\int_\Theta^{2+W}P_W/P^2$ and
$s_{WW}=P_W(2+W)+\int_\Theta^{2+W}(P_{WW}/P^2-2P_W^2/P^3)$, with no dependence on a point
$\Theta_m$, and the three bounds follow as in (a), using $|P_\Theta|=|\Theta_{ss}|/\Theta_s\le C$.
At $s=2$, $\Theta_s=1$ and $\Theta_W=1$ from both sides. Since $s_c(W)$ and $s_*(W)$ are the values
of $s$ at $\Theta=0$ and $\Theta=-\bar q$, their derivatives are values of $s_W$, which is bounded by
$C(1+|s|)\le C\log(1/\gamma)$.
\end{proof}

\subsection{Strict super- and subsolutions and their extensions}

For a function $\tilde u\ge0$ on $\overline{Q_T}$ we write $\Omega=\{\tilde u>0\}$,
$\Omega_t=\{x:\tilde u(x,t)>0\}$, $\Gamma_t=\partial\Omega_t$,
$\Gamma=\bigcup_t\Gamma_t\times\{t\}$, and $d(x,t)=\operatorname{dist}(x,\Gamma_t)$.

\begin{definition}\label{q:def:strict}
Let $\delta,\delta_0,r,K_0,C_1>0$. A continuous function $\tilde u\ge0$ on
$\overline{Q_T}$, with $\Omega_t$ bounded, is a \emph{regular strict supersolution with
margins $(\delta,\delta_0,r,K_0)$} if
\begin{enumerate}
\item[(R1)] $\tilde u|_{\overline\Omega}$ has an extension $E\tilde u$ to
$\overline\Omega\cup\{d<r\}$ which is $C^{1,1}$ in $x$ and Lipschitz in $t$ (on this set,
which is not a cylinder: $E\tilde u$ is Lipschitz in $(x,t)$ there), with
$|\nabla E\tilde u|+|D^2E\tilde u|+|\partial_tE\tilde u|\le K_0$, and $E\tilde u<0$ outside
$\overline\Omega$;
\item[(R2)] $|\nabla\tilde u|\ge2\delta_0$ on $\Gamma$ and $\tilde u\ge\delta_0\min\{d,r\}$ in $\Omega$;
\item[(R3)] $\tilde u\in C^{2,1}(\Omega)$ and $\Delta\tilde u-\tilde u_t\le0$ in $\Omega$;
\item[(R4)] $|\nabla\tilde u|\le\sqrt{2M(w_0)}-\delta$ on $\Gamma$.
\end{enumerate}
It is a \emph{regular strict subsolution} if (R1), (R2) hold,
$\Delta\tilde u-\tilde u_t\ge0$ in $\Omega$, and
$\sqrt{2M(w_0)}+\delta\le|\nabla\tilde u|\le\sqrt{2M(w_0)}+C_1\delta$ on $\Gamma$.
\end{definition}

By (R1)--(R2) and the implicit function theorem, $\Gamma_t=\{E\tilde u(\cdot,t)=0\}$ is a
$C^{1,1}$ hypersurface, which moves with normal velocity bounded by $K_0/(2\delta_0)$; in
particular $d$ is Lipschitz in $t$. Conversely, (R1) holds when $\tilde u$ is $C^{2,1}$ up to
$\Gamma$ and $\Gamma$ is locally the graph of a function $x_N=g(x',t)$ with bounded
$D^2_{x'}g$ and $\partial_tg$: one flattens $\Gamma$ and uses the reflection
$Eu(y',y_N)=3u(y',-y_N)-2u(y',-2y_N)$, as in \cite[Lemma~6.37]{GilbargTrudinger2001}.

\begin{lemma}[Monotone extensions]\label{q:lem:ext}
Let $\tilde u$ satisfy (R1), (R2), and let $\Lambda\ge1$. There are $r_1\in(0,r]$ and
$K_1$, depending only on $\delta_0,r,K_0,\Lambda$, and functions $\tilde v^+$ and
$\tilde v^-$ on $\overline{Q_T}$, $C^{1,1}$ in $x$ and Lipschitz in $t$, such that
\begin{enumerate}
\item[(i)] $\tilde v^\pm=\tilde u$ in $\overline\Omega$, and
$-2K_0d\le\tilde v^\pm\le-\frac{\delta_0}2\min\{d,r_1\}$ outside $\Omega$ (the lower bound for
$d\le r_1$);
\item[(ii)] $|\nabla\tilde v^\pm|\le2K_0$ and
$|D^2\tilde v^\pm|+|\partial_t\tilde v^\pm|\le K_1$ a.e.; in particular
$|\Delta\tilde v^\pm-\tilde v^\pm_t|\le K_1$;
\item[(iii)] if $x\notin\Omega_t$, $d=d(x,t)\le r_1$ and $y$ is a closest point of
$\Gamma_t$ to $x$, then
\[
|\nabla\tilde v^+(x,t)|^2\le|\nabla\tilde u(y,t)|^2-\Lambda d,\qquad
|\nabla\tilde u(y,t)|^2+\Lambda d\le|\nabla\tilde v^-(x,t)|^2\le|\nabla\tilde u(y,t)|^2+K_1d .
\]
\end{enumerate}
The constants $2K_0$ and $\delta_0/2$ in (i) and (ii) do not depend on $\Lambda$ (only $r_1$ and
$K_1$ do).
\end{lemma}

\begin{proof}
Let $x\notin\Omega_t$ with $d\le r$ and let $y$ be a closest point of $\Gamma_t$. Then
$x-y$ is normal to $\Gamma_t$, pointing outwards, and $\nabla\tilde u(y)$ is normal, pointing
inwards, so $E\tilde u(x)=-|\nabla\tilde u(y)|d+R$ with $|R|\le\frac12K_0d^2$ and
$2\delta_0\le|\nabla\tilde u(y)|\le K_0$. Hence
\begin{equation}\label{q:eq:Eu}
-\tfrac32K_0d\le E\tilde u\le-\delta_0d\qquad\text{for }x\notin\Omega_t,\ d\le r_0:=\min\{r,\delta_0/K_0,1\}.
\end{equation}
Let $\Lambda'=(\Lambda+2K_0^2+\delta_0K_0)/(8\delta_0^3)$ and
$m=\min\{\delta_0r_0/4,\,1/(8\Lambda')\}$. Take $\psi\in C^\infty(\R)$ nondecreasing with
$\psi(\tau)=\tau$ for $\tau\ge-m/2$, $\psi(\tau)=-m$ for $\tau\le-2m$ and $0\le\psi'\le1$, and set
$\hat E=\psi(E\tilde u)$ on $\overline\Omega\cup\{d<r_0\}$ and $\hat E=-m$ elsewhere. By
\eqref{q:eq:Eu}, $\hat E=-m$ where $r_0/2\le d<r_0$ outside $\Omega$, so $\hat E$ is $C^{1,1}$
in $x$ and Lipschitz in $t$ on $\overline{Q_T}$; $\hat E=\tilde u$ in $\overline\Omega$; $-m\le\hat E<0$
outside; and $\hat E=E\tilde u$ where $d\le r_1:=m/(3K_0)$, because there
$E\tilde u\ge-m/2$ by \eqref{q:eq:Eu}. Let
\[
\tilde v^\pm=\hat E\pm\Lambda'\big(\hat E^-\big)^2,\qquad \hat E^-=\min\{\hat E,0\}.
\]
Then $\nabla\tilde v^\pm=(1\mp2\Lambda'|\hat E^-|)\nabla\hat E$ with
$2\Lambda'|\hat E^-|\le\frac14$, which gives (ii) and
$\frac34|\hat E|\le|\tilde v^\pm|\le\frac54|\hat E|$ outside $\Omega$; with \eqref{q:eq:Eu} and
the monotonicity of $\psi$ this gives (i) (note that $\frac54\cdot\frac32<2$ and
$\frac34>\frac12$). For (iii), $\hat E=E\tilde u$, and
$\big||\nabla E\tilde u(x)|^2-|\nabla\tilde u(y)|^2\big|\le2K_0^2d$. Therefore, using
\eqref{q:eq:Eu} and $|\nabla\tilde u(y)|\ge2\delta_0$,
\[
|\nabla\tilde v^+|^2\le(1-2\Lambda'\delta_0d)\big(|\nabla\tilde u(y)|^2+2K_0^2d\big)
\le|\nabla\tilde u(y)|^2-(8\Lambda'\delta_0^3-2K_0^2)d,
\]
\[
|\nabla\tilde v^-|^2\ge(1+2\Lambda'\delta_0d)\big(|\nabla\tilde u(y)|^2-2K_0^2d\big)
\ge|\nabla\tilde u(y)|^2+(8\Lambda'\delta_0^3-2K_0^2-4\Lambda'\delta_0K_0^2r_1)d,
\]
and $4\Lambda'\delta_0K_0^2r_1\le\delta_0K_0$. The upper bound for $|\nabla\tilde v^-|^2$ is
immediate.
\end{proof}

\subsection{Quantitative approximation}

\begin{proposition}[Quantitative version of {\cite[Theorem~4.1]{FernandezBonderWolanski2004}}]\label{q:prop:super}
Assume (F$_2$), (W), and let $\tilde q\ge1$. There are $\Lambda,K,c_1,c_*,\delta_*,C$
depending only on
$f,\delta_1,\bar W,C_W,\delta_0,r,K_0$ and $\tilde q$ (\emph{not} on $\delta$) such that, if
$\tilde u$ is a regular strict supersolution with margins $(\delta,\delta_0,r,K_0)$ and
\[
\delta\le\delta_*,\qquad\ve+\eta_\ve\le c_*\delta,
\]
then, with $\tilde v=\tilde v^+$ from Lemma~\ref{q:lem:ext} (with this $\Lambda$) and
$\Theta=\Theta^+$ from Lemma~\ref{q:lem:inner}(a) with $\gamma=c_1\delta$ and $\kappa=K\ve$,
\[
V^\ve=-w^\ve+\ve\,\Theta\Big(\frac{\tilde v}\ve;W^\ve\Big)+\ve^2\tilde q
\]
is a strong supersolution of \eqref{eq:Peps} in $Q_T$, and
$\|V^\ve-\tilde u\|_{L^\infty(Q_T)}\le C\ve$.
\end{proposition}

\begin{proof}
In this proof $C,c$ denote constants that depend on $f,\delta_1,\bar W,C_W,\delta_0,r,K_0$
but not on $\Lambda$, $K$, $\tilde q$, $\delta$, $\ve$. We write $s=\tilde v/\ve$.

\emph{Margin.} We claim that, where $s\le s_m$ and $d\le r_1$,
\begin{equation}\label{q:eq:marginsup}
G_-(W)^2-|\nabla\tilde v|^2\ge c\,\delta+\Lambda_1\ve|s|,\qquad \Lambda_1=\frac{\Lambda-C}{2K_0}.
\end{equation}
Let $y$ be a closest point of $\Gamma_t$. Since $|W(x,t)-w_0(y,t)|\le\eta_\ve+C_Wd$ and $M$ is
Lipschitz, $G_-(W)^2\ge2M(w_0(y))-2\gamma-C(\eta_\ve+d)$, while by (R4)
$|\nabla\tilde u(y)|^2\le2M(w_0(y))-c\,\delta$ for $\delta\le\delta_*$. In $\Omega$, $0<s\le s_m<1$ means
$\tilde u\le\ve$, so $d\le\ve/\delta_0$ by (R2), and
$|\nabla\tilde v|^2\le|\nabla\tilde u(y)|^2+C\ve$. Outside $\Omega$, Lemma~\ref{q:lem:ext} gives
$|\nabla\tilde v|^2\le|\nabla\tilde u(y)|^2-\Lambda d$ and $d\ge\ve|s|/(2K_0)$. With
$\gamma=c_1\delta$, $c_1$ small, and $\ve+\eta_\ve\le c_*\delta$, this gives \eqref{q:eq:marginsup}.

\emph{The equation where $s<s_m$.} Write $V=-w+\ve(\Theta+\ve\tilde q)$. Since $f$ is
Lipschitz and $\Theta>0$,
\[
(V+w)f_\ve(V)=\tfrac1\ve(\Theta+\ve\tilde q)f(\Theta+\ve\tilde q-W)
\ge\tfrac1\ve b(\Theta;W)-L_f\tilde q\,\Theta+\tilde q\,f(\Theta+\ve\tilde q-W).
\]
By \eqref{q:eq:basic} and the equation of $\Theta$,
\[
\mathcal R:=\Delta V-V_t-(V+w)f_\ve(V)\le\Theta_s(\Delta\tilde v-\tilde v_t)-K\Theta_s
-\frac1\ve\Theta_{ss}\big(G_-^2-|\nabla\tilde v|^2\big)+|\mathcal P|+L_f\tilde q\,\Theta-\tilde q\,f(\Theta+\ve\tilde q-W),
\]
and, by Lemma~\ref{q:lem:inner}(a)(ii), $|\nabla\tilde v|\le2K_0$ and $|\nabla W|\le C_W$,
\begin{equation}\label{q:eq:Pbound}
|\mathcal P|\le C(1+|s|)\Theta_{ss}+C\Theta_s+C\ve(1+|s|)^2\Theta_s .
\end{equation}
Recall that $\Delta\tilde v-\tilde v_t\le K_1$ (and $\le0$ in $\Omega$). Let
$\Theta_0=\min\{\delta_1/8,a_0/L_f\}$. We distinguish three regions.

\emph{(A) $\Theta\le\Theta_0$ and $s\ge-\ve^{-1/2}$.} By Lemma~\ref{q:lem:ext}(i),
$d\le r_1$ for $\ve$ small, so \eqref{q:eq:marginsup} holds and
$\frac1\ve\Theta_{ss}(G_-^2-|\nabla\tilde v|^2)\ge\Theta_{ss}(c\,\delta/\ve+\Lambda_1|s|)\ge C(1+|s|)\Theta_{ss}$
if $\Lambda$ is large and $c_*$ small. Also $\ve(1+|s|)^2\le4$. Since
$\Theta+\ve\tilde q-W\le1-\delta_1/4$ for $\ve\tilde q\le\delta_1/8$,
$f(\Theta+\ve\tilde q-W)\ge2a_0$, and $L_f\tilde q\,\Theta\le a_0\tilde q$. Hence
$\mathcal R\le(K_1+C-K)\Theta_s-a_0\tilde q\le0$ if $K\ge K_1+C$.

\emph{(B) $\Theta>\Theta_0$ (and $s<s_m$).} By Lemma~\ref{q:lem:inner}(a)(i), $s\ge-S_0$ with
$S_0$ depending only on $\Theta_0$, and $\Theta_s\ge p_0:=c\,e^{-CS_0}$. As in (A) the margin
absorbs the term $C(1+|s|)\Theta_{ss}$, and $L_f\tilde q\,\Theta\le C\tilde q\le(C\tilde q/p_0)\Theta_s$. Hence
$\mathcal R\le(K_1+C+C\tilde q/p_0-K)\Theta_s\le0$ if $K$ is large.

\emph{(C) $s<-\ve^{-1/2}$.} Here $\Theta$, $\Theta_s$ and $\Theta_{ss}$ are bounded by
$Ce^{-\lambda_0\ve^{-1/2}}$ and $|s|\le C/\ve$, so every term except
$-\tilde q\,f(\Theta+\ve\tilde q-W)\le-2a_0\tilde q$ is bounded by
$C(K_1+K)\ve^{-3}e^{-\lambda_0\ve^{-1/2}}\le a_0$ for $\ve$ small, and $\mathcal R\le0$.

\emph{The equation where $s\ge s_m$.} Here $V=\tilde u+\ve^2\tilde q$. Either $s\ge1$ and
$f_\ve(V)=0$, or $V+w=\ve(s+W+\ve\tilde q)\ge\ve(s_m-1+\delta_1/2)>0$; in both cases
$\mathcal R\le\Delta\tilde u-\tilde u_t\le0$.

\emph{Regularity and choice of the constants.} By Lemma~\ref{q:lem:inner}(a)(ii), $V$ is
$C^{1,1}$ in $x$ and Lipschitz in $t$, and the set $\{\tilde v=\ve s_m(W)\}$ has measure zero,
because $|\nabla(\tilde v-\ve s_m(W))|\ge2\delta_0-C\ve>0$ near it. So $V$ is a strong
supersolution. The constants are chosen in this order: $\Lambda$ (which fixes $r_1$, $K_1$
and $\Lambda_1$); then $K$; then $c_1$; then $c_*$, so that $\kappa=K\ve\le Kc_*\delta\le c\gamma$
(that is, $c_*\le c\,c_1/K$) and the
conditions above hold; finally $\delta_*$ and $\ve$ small. In this proof the monotonicity of
$f$ is used only through Lemma~\ref{q:lem:inner}(a) and the constant $a_0$.

\emph{Error.} With $\tilde u=\ve s^+$, $V-\tilde u=\ve(\Theta-W-s^+)+\ve^2\tilde q$, and
$\Theta-W-s^+$ vanishes for $s\ge s_m$, while $0<\Theta\le1+W$ for $s<s_m$.
\end{proof}

\begin{proposition}[Quantitative version of {\cite[Theorem~4.2]{FernandezBonderWolanski2004}}]\label{q:prop:sub}
Assume (F$_2$), (W), and let $\tilde q\ge1$. There are $c_*,\delta_*,C$ as in
Proposition~\ref{q:prop:super}, depending also on $C_1$, such that, if $\tilde u$ is a
regular strict subsolution with margins $(\delta,\delta_0,r,K_0)$, $\delta\le\delta_*$ and
$\ve+\eta_\ve\le c_*\delta$, then, with $\tilde v=\tilde v^-$ from Lemma~\ref{q:lem:ext} and
$\Theta=\Theta^-$ from Lemma~\ref{q:lem:inner}(b) (for constants $\Lambda$, $K$, $c_1$ as in
Proposition~\ref{q:prop:super}, $\gamma=c_1\delta$,
$\kappa=K\ve$, and $\bar q=\ve\tilde q$), the function
\[
V_1=-w^\ve+\ve\,\Theta\Big(\frac{\tilde v}\ve;W^\ve\Big)
\]
is a strong subsolution of \eqref{eq:Peps} in the bounded open set
$D=\{(x,t)\in Q_T:\tilde v>\ve s_*(W^\ve)\}$, and $V_1=-w^\ve-\ve^2\tilde q$ on $\partial D\cap Q_T$.
Consequently, if $u$ is a bounded solution of \eqref{eq:Peps} with $u+w^\ve\ge0$, and
\[
\underline V^\ve=\begin{cases}V_1&\text{in }D,\\ -w^\ve-\ve^2\tilde q&\text{in }Q_T\setminus D,\end{cases}
\]
satisfies $\underline V^\ve(\cdot,0)\le u(\cdot,0)$, then $\underline V^\ve\le u$ in $Q_T$. Moreover
$\|\underline V^\ve-\tilde u\|_{L^\infty(Q_T)}\le C\ve$.
\end{proposition}

\begin{proof}
We use the conventions of the previous proof. On $\partial D\cap Q_T$, $s=s_*(W)$ and
$\Theta=-\bar q$, so $V_1=-w-\ve^2\tilde q$.

\emph{Margin.} In $D$, where $s\le2$ we have
$d\le C\ve(1+|s_*|)\le C\ve\log(1/\gamma)\le r_1$, and, as in the previous proof, now with
Lemma~\ref{q:lem:ext}(iii) for $\tilde v^-$ and the two-sided bound for $|\nabla\tilde u|$ on
$\Gamma$,
\begin{equation}\label{q:eq:marginsub}
c\,\delta+\Lambda_1\ve|s|\le|\nabla\tilde v|^2-G_+(W)^2\le C_\Lambda(\delta+\ve|s|),
\end{equation}
where $C_\Lambda$ depends also on $\Lambda$ and $C_1$.

\emph{The equation where $s_*<s<2$.} By \eqref{q:eq:basic} and the equation of $\Theta$,
\[
\mathcal R=\Theta_s(\Delta\tilde v-\tilde v_t)+K\Theta_s+\frac1\ve\Theta_{ss}\big(|\nabla\tilde v|^2-G_+^2\big)+\mathcal P ,
\]
with $\mathcal P$ bounded as in \eqref{q:eq:Pbound} (with $|\Theta_{ss}|$), $\Delta\tilde v-\tilde v_t\ge-K_1$,
and $\ve(1+|s|)^2\le C\ve\log^2(1/\gamma)\le1$. Where $\Theta\ge0$, $\Theta_{ss}\ge0$ and
\eqref{q:eq:marginsub} gives
$\frac1\ve\Theta_{ss}(|\nabla\tilde v|^2-G_+^2)\ge C(1+|s|)\Theta_{ss}$ for $\Lambda$ large and $c_*$
small, so $\mathcal R\ge(K-K_1-C)\Theta_s\ge0$. Where $-\ve\tilde q\le\Theta<0$, by
Lemma~\ref{q:lem:inner}(b)(i), $|\Theta_{ss}|\le C(\ve\tilde q+K\ve\Theta_s)$ and
$\Theta_s\ge c\sqrt\gamma$; with the upper bound in \eqref{q:eq:marginsub} and
$|s|\le C\log(1/\gamma)$,
\begin{align*}
\mathcal R&\ge(K-K_1-C)\Theta_s-C_{\Lambda,K}\,(\tilde q+\Theta_s)\big(\delta+\ve\log(1/\gamma)\big)
\ge\Theta_s-C_{\Lambda,K}\,\tilde q\,\big(\sqrt{\delta_*}+\sqrt{c_*}\big)\sqrt\delta\\
&\ge\big(c\sqrt{c_1}-C_{\Lambda,K}\,\tilde q\,(\sqrt{\delta_*}+\sqrt{c_*})\big)\sqrt\delta\ge0,
\end{align*}
if $K\ge K_1+C+2$ and then $\delta_*$, $c_*$ are small; we used $\gamma\ge\ve$ (that is,
$c_*\le c_1$) and
$\ve\log(1/\ve)\le\sqrt{\ve}\le\sqrt{c_*\delta}$. Note that $C_{\Lambda,K}$ does not depend on
$c_1$, and that $\delta_*$ and $c_*$, which are chosen last, depend on $\tilde q$. The set $D$ is
bounded because $\tilde v^-$ is a negative constant away from $\overline\Omega$
(proof of Lemma~\ref{q:lem:ext}).

\emph{The equation where $s\ge2$.} $V_1=\tilde u\ge2\ve$, there is no reaction, and
$\mathcal R=\Delta\tilde u-\tilde u_t\ge0$. By Lemma~\ref{q:lem:inner}(b)(ii), $V_1$ is $C^{1,1}$ in
$x$ and Lipschitz in $t$ in $D$, and the sets where $s=2$ or $\Theta=1+W$ have measure zero,
because $\nabla\tilde v\ne0$ there and $\Theta$ is strictly increasing in $s$ with $|\Theta_W|$
bounded.

\emph{Comparison.} The function $-w-\ve^2\tilde q\le-w\le u$. On the parabolic boundary of
$D$ we have $V_1\le u$: at $t=0$ by assumption, and on $\partial D\cap Q_T$ because
$V_1=-w-\ve^2\tilde q$ there. Hence $V_1\le u$ in $D$.

\emph{Error.} In $D$, $V_1-\tilde u=\ve(\Theta-W-s^+)$ with $-\ve\tilde q\le\Theta\le2+W$ for
$s\le2$; outside $D$, $\tilde u=0$ and $\underline V^\ve=-w-\ve^2\tilde q$.
\end{proof}

\subsection{Rate of convergence to a classical solution}

We recall the setting of \cite{FernandezBonderWolanski2004}, to which we add a
Lipschitz condition on $w_0$: $u_0\ge0$ is Lipschitz with
compact support, and $w_0$ satisfies the monotonicity condition
\begin{equation}\label{q:eq:mono}
\lambda\mapsto w_0(\lambda x,\lambda^2t)\ \text{ is nondecreasing for }\lambda>0,
\qquad |w_0(\lambda x,\lambda^2t)-w_0(x,t)|\le C_w|\lambda-1|
\end{equation}
for $\frac12\le\lambda\le2$ and $(x,t)$ in bounded sets (for instance, $w_0$ constant).
When $\eta_\ve\to0$, $w_0$ is a bounded caloric function, and so is
$(x,t)\mapsto w_0(\lambda x,\lambda^2t)$; hence the monotonicity in \eqref{q:eq:mono} holds for
all $t$ if it holds at $t=0$, by the maximum principle, as in
\cite{FernandezBonderWolanski2004}.
Since $M$ is increasing, $M(w_0(\lambda x,\lambda^2t))\le M(w_0(x,t))$ for $\lambda\le1$, with the
reverse inequality for $\lambda\ge1$.

We say that $v$ is a \emph{regular classical solution} of \eqref{eq:P} on $[0,T^*]$ if
$v\ge0$ is continuous on $\R^N\times[0,T^*]$, $v\in C^{2,1}(\{v>0\})$, $v_t=\Delta v$ in
$\{v>0\}$, $|\nabla v|=\sqrt{2M(w_0)}$ on $\partial\{v>0\}$, and $v$ satisfies (R1), (R2) on
$[0,T^*]$ for some $\delta_0,r,K_0>0$. Regularity up to $t=0$ is a compatibility
condition on the data: $v(\cdot,0)=u_0$ is $C^{1,1}$ up to
$\partial\Omega_0$ and satisfies the free boundary condition $|\nabla u_0|=\sqrt{2M(w_0)}$
there. Recall that $\Omega_t=\{v(\cdot,t)>0\}$, the fresh region, is bounded
(Definition~\ref{q:def:strict}).

\begin{theorem}[Rate of convergence]\label{q:thm:rate}
Assume (F$_2$), (W), \eqref{q:eq:mono}, and
\begin{itemize}
\item[(H)] $\Omega_0=\{u_0>0\}$ is strictly starshaped with respect to $0$
($x\cdot\nu\ge\varrho_0>0$ on $\partial\Omega_0$), $u_0(\lambda x)\ge u_0(x)$ for
$0<\lambda<1$, and $x\cdot\nabla u_0\le-\varrho_0$ in
$\Omega_0\cap\{\operatorname{dist}(x,\partial\Omega_0)<\varrho_0\}$.
\end{itemize}
Let $v$ be a regular classical solution of \eqref{eq:P} on $[0,T^*]$ with $v(\cdot,0)=u_0$,
and let $u^\ve$ be the solution of \eqref{eq:Peps} with initial data such that
\begin{equation}\label{q:eq:data}
Y_0^\ve:=u_0^\ve+w^\ve(\cdot,0)\ge0,\qquad |u_0^\ve-u_0|\le C_0\ve\ \text{ in }\Omega_0,\qquad
Y_0^\ve\le C_0\ve^2\ \text{ outside }\Omega_0 .
\end{equation}
Assume that (W) holds in $Q_{T^*}$. Then for every $T<T^*$ there are $C,\ve_1>0$ such that,
if $\ve+\eta_\ve\le\ve_1$,
\[
-C\Big(\ve\log\frac1\ve+\eta_\ve\Big)\le u^\ve-v\le C(\ve+\eta_\ve)\qquad\text{in }\R^N\times[0,T].
\]
\end{theorem}

The data of Section~\ref{sec:method}, $u_0^\ve=u_0-\ve w_0(1-u_0/\ve)^+$ and
$w^\ve(\cdot,0)=\ve w_0$ with $w_0$ constant, satisfy \eqref{q:eq:data} and (W) with
$\eta_\ve=0$ and $Y_0^\ve=0$ outside $\Omega_0$.

\begin{proof}
Let $\ve'=\ve+\eta_\ve$. The constants $C,c$ below do not depend on $\ve$ or on the
parameters $a,a_1$.

\emph{Rescalings.} For $\lambda$ close to $1$ let $\lambda'=1+(\lambda-1)/2$ and
$v_\lambda(x,t)=\frac1{\lambda'}v(\lambda x,\lambda^2t)$, which is defined on $[0,T]$ if
$\lambda^2T\le T^*$. It is caloric in its positivity set, and on its free boundary
$|\nabla v_\lambda|(x,t)=\frac\lambda{\lambda'}\sqrt{2M(w_0(\lambda x,\lambda^2t))}$. By \eqref{q:eq:mono},
as in \cite[Section~2]{FernandezBonderWolanski2004}, $v_\lambda$ is, for $\lambda<1$, a regular
strict supersolution, and for $\lambda>1$ a regular strict subsolution, with margins
$(\delta,\delta_0,r,K_0)$, where $\delta_0,r,K_0,C_1$ do not depend on $\lambda$ and
$c|\lambda-1|\le\delta\le C|\lambda-1|$. Since $v$ is bounded and Lipschitz in $x$ and $t$
(because $v_t$ is bounded in $\{v>0\}$ and $v$ is continuous),
\begin{equation}\label{q:eq:vlam}
|v_\lambda-v|\le C|\lambda-1|\qquad\text{in }\R^N\times[0,T].
\end{equation}
By (H), there is $c>0$ such that, for $\lambda<1$ and for $\lambda>1$ respectively,
\begin{equation}\label{q:eq:init}
v_\lambda(\cdot,0)\ge u_0+c(1-\lambda)\ \text{ in }\overline{\Omega_0},\qquad
u_0\ge v_\lambda(\cdot,0)+c(\lambda-1)\ \text{ in }\lambda^{-1}\Omega_0=\{v_\lambda(\cdot,0)>0\}.
\end{equation}
Indeed, near $\partial\Omega_0$, $\frac{d}{d\tau}u_0(\tau x)=\frac1\tau(\tau x)\cdot\nabla u_0(\tau x)\le-\varrho_0/\tau$,
so $|u_0(\lambda x)-u_0(x)|\ge c|\lambda-1|$ with the appropriate sign, and
$v_\lambda(x,0)=u_0(\lambda x)/\lambda'$; away from $\partial\Omega_0$, $u_0\ge c$, $u_0$ is monotone
along rays, and $|1/\lambda'-1|\ge c|\lambda-1|$.

\emph{Upper bound.} Let $\lambda=1-a\ve'$, so that $\delta\ge c\,a\ve'$. We apply
Proposition~\ref{q:prop:super} with $\tilde u=v_\lambda$ and $\tilde q=\max\{1,C_0\}$, taking $a$ so
large that $\ve'\le c_*\delta$: $V^\ve$ is a strong supersolution with
$|V^\ve-v_\lambda|\le C\ve$. At $t=0$: in $\Omega_0$, by \eqref{q:eq:init},
$V^\ve\ge v_\lambda-C\ve\ge u_0+(ca-C)\ve'\ge u_0^\ve$ if $a$ is large; outside $\Omega_0$,
$V^\ve\ge-w^\ve+\ve^2\tilde q\ge u_0^\ve$, because $\Theta>0$. By comparison and
\eqref{q:eq:vlam}, $u^\ve\le V^\ve\le v_\lambda+C\ve\le v+C\ve'$.

\emph{Lower bound.} Let $\lambda=1+a\ve'+a_1\ve\log(1/\ve)$, so that
$c(\lambda-1)\le\delta\le C(\lambda-1)\le\delta_*$ for $\ve$ small, and apply
Proposition~\ref{q:prop:sub} with $\tilde u=v_\lambda$, $\tilde q=1$, and $a$ large. Note that
$u^\ve+w^\ve\ge0$ by the maximum principle. We check that
$\underline V^\ve(\cdot,0)\le u_0^\ve$. Outside $D$ this is $Y_0^\ve\ge0$. In $D$,
$\underline V^\ve=V_1=-w+\ve\Theta$; where $s\le s_c$ we have $\Theta\le0$ and again
$V_1\le-w\le u_0^\ve$. It remains to consider the points with $s>s_c$. Since
$\Theta\le2+W$ for $s\le2$, we have $V_1\le\max\{v_\lambda,2\ve\}$.
\begin{itemize}
\item If $x\in\lambda^{-1}\Omega_0$: by \eqref{q:eq:init},
$u_0^\ve\ge u_0-C_0\ve\ge v_\lambda+c\,a\ve'-C_0\ve\ge v_\lambda+2\ve\ge V_1$ for $a$ large.
\item If $x\notin\lambda^{-1}\Omega_0$: let $\varrho$ be the distance from $x$ to
$\lambda^{-1}\partial\Omega_0$, and let $d_0$ be the distance from $x$ to $\partial\Omega_0$ if
$x\in\Omega_0$, and $d_0=0$ otherwise. By strict starshapedness,
$d_0+\varrho\ge c(\lambda-1)$. If $d_0\ge c(\lambda-1)/2$, then, since
$u_0\ge c\operatorname{dist}(x,\partial\Omega_0)$ near $\partial\Omega_0$ by (H), $u_0\ge c\,d_0\ge c\,a\ve'$, so
$u_0^\ve\ge u_0-C_0\ve\ge2\ve\ge V_1$ for $a$ large. Otherwise
$\varrho\ge c(\lambda-1)/2\ge ca_1\ve\log(1/\ve)$, and by Lemma~\ref{q:lem:ext}(i)
$s=\tilde v(x,0)/\ve\le-ca_1\log(1/\ve)$, which is smaller than $s_c$ if $a_1$ is large, because
$s_c\ge-C\log(1/\gamma)$ and $\gamma\ge c\,\ve$ (Lemma~\ref{q:lem:inner}(b)(i)): such
points do not have $s>s_c$.
\end{itemize}
Therefore $u^\ve\ge\underline V^\ve\ge v_\lambda-C\ve\ge v-C(\lambda-1)-C\ve$, by
\eqref{q:eq:vlam}.
\end{proof}

\begin{remark}\label{q:rem:rate}
\begin{enumerate}
\item \emph{Why the rate is $\ve$ and not $\sqrt\ve$.} The strictness $\delta$ enters only
through the conditions $\delta\le\delta_*$ and $\ve+\eta_\ve\le c_*\delta$; the constants of
Propositions~\ref{q:prop:super} and \ref{q:prop:sub} do not depend on $\delta$. This is what
allows $|1-\lambda|\sim\ve$ in the rescaling argument. It rests on two choices: the
profiles are parametrized by $W$ only, with $G_\pm(W)^2-2M(W)=\pm2\gamma$ independent of
$W$, so that their derivatives with respect to $W$ are controlled uniformly in $\gamma$;
and the extension $\tilde v^\pm$ has a gradient that is monotone in the normal direction
in the burnt region, which provides the margin \eqref{q:eq:marginsup},
\eqref{q:eq:marginsub} that absorbs the terms coming from the variation of $W$ across
the layer.
\item \emph{The logarithm.} It comes only from the initial ordering: when
$Y_0^\ve$ vanishes in the burnt region, the subsolution must reach the burnt state
$u=-w$ within the initial gap, and its profile takes a length $O(\log(1/\gamma))$ to do
so. If the data satisfy $Y_0^\ve\ge c\,\ve\, e^{-\beta\operatorname{dist}(x,\Omega_0)/\ve}$ with
$\beta$ small, the bound $\Theta\le(2+W)e^{c(s-2)}$ of Lemma~\ref{q:lem:inner}(b) gives the
initial ordering with $a_1=0$, and the lower bound becomes $u^\ve-v\ge-C(\ve+\eta_\ve)$.
Such data, however, do not satisfy the third condition in \eqref{q:eq:data}, which we
used for the upper bound; a two-sided $O(\ve)$ estimate for data with the tail of a
travelling wave would require sharp decay rates of the profiles $\Theta^\pm$, which we do
not pursue. We do not know whether the logarithm is needed when $Y_0^\ve=0$ in the burnt
region; the experiments of Section~\ref{sec:experiments}, which use such data, do not show
it. A natural way to remove it would be to wait for a time of order $\ve^2\log(1/\ve)$, in
which $Y^\ve$ develops an exponential tail in the burnt region, and then to start the
comparison from that time. This requires two estimates on that initial layer that we
have not been able to prove: a lower bound for the tail of $Y^\ve$ with the same exponent as
the barrier, whose decay rate differs from the natural one, $\sqrt{f(-w_0)}/\ve$
per unit length, only by a factor $1+o(1)$, so that a crude bound does not suffice; and a
lower bound for $u^\ve$ in the fresh
region that loses only $O(\ve)$.
\item \emph{The order $\ve$.} By the upper barrier, $0\le u^\ve+w^\ve\le C\ve^2$
at the points of $\{v=0\}$ at distance at least $C(\ve\log(1/\ve)+\eta_\ve)$ from the free
boundary.
There $u^\ve-v=-\ve W^\ve+O(\ve^2)$, so the order $\ve$ cannot be improved in the sup norm
when $w_0\not\equiv0$. This difference is the offset between the temperature of the burnt
gas, $T_f+\ve w_0$, and the reference $T_f$; when $w_0\equiv0$ it gives no information, and
a more significant quantity is the position of the front. For it we only have the example
of the travelling waves of Section~\ref{sec:tw}, which are not covered by the theorem,
since their fresh region is not bounded. By Proposition~\ref{prop:speed},
$c_\ve-c_0=\beta\,\ve+O(\ve^2)$ with $\beta=\frac{c_0}A\big(1-\frac I{\sqrt{2M}}\big)$, which is not
zero in general (for $f(s)=(1-s)^+$ it vanishes only for $w_0=w_0^*$). Hence the distance
travelled in a time $t$ by any level set of $u^\ve$ differs from that of the free
boundary by $\beta\,\ve\,t+O(\ve^2t)$, whatever level is used to define the front. In
addition, at each time the front $\{u^\ve=\ve\}$ is ahead of the free boundary of the
outer solution: for $z\ge0$, $\Phi_\ve(z)=A(1-e^{-c_\ve(z-z_\ve)})$ with
$z_\ve=\log(1-\ve/A)/c_\ve=-\ve/\sqrt{2M}+O(\ve^2)$; this shift depends on the level that
defines the front, and it is observed for general data in Section~\ref{sec:limit}.
\item \emph{Existence of regular classical solutions.} The theorem is conditional on
the existence of a regular classical solution, which includes the compatibility of the
data at $t=0$. Local existence of solutions of class $C^{1+\alpha/2,2+\alpha}$
up to $t=0$ is proved in \cite[Theorem~2.8]{BaconneauLunardi2004} for free boundary
conditions of the form $u=g_0(x,t)$, $\partial_nu=g_1(x,t)$ with $g_1\ne\partial_ng_0$, when
$\partial\Omega_0$ and $u_0$ are of class $C^{3+\alpha}$ and satisfy the compatibility
conditions; with $g_0=0$ and $g_1=-\sqrt{2M(w_0)}$ this covers \eqref{eq:P}, for constant
$w_0$ and also for a sufficiently smooth $w_0(x,t)$. Hence, if in addition (H),
\eqref{q:eq:mono} and \eqref{q:eq:data} hold, the theorem applies on a short time
interval (on which the solution is positive in $\Omega_t$, by continuity). The radial
self-similar solutions of Section~\ref{sec:selfsimilar} are
explicit regular classical solutions that satisfy all the assumptions of the theorem
on $[0,T^*]$ for every $T^*<T$, that is, up to any time before extinction. In all cases
the fresh region is a bounded pocket.
\item \emph{The constants.} The constant $C$ depends on $f,\delta_1,\bar W,C_W,C_w,C_0$ and
$\varrho_0$, on the margins $\delta_0,r,K_0$ of $v$ on $[0,T^*]$, and on the Lipschitz constant
of $v$ in $(x,t)$ and the diameter of its support, through \eqref{q:eq:vlam}; $\ve_1$
depends in addition on $T^*-T$, because the dilations must satisfy $\lambda^2T\le T^*$.
Assumption (H) is used only at $t=0$. No Gronwall argument is involved, so $C$ does not
grow exponentially with $T$.
\end{enumerate}
\end{remark}

\section{The numerical method}\label{sec:method}

We solve \eqref{eq:Peps} in planar ($N=1$) or radial ($N=2,3$) symmetry on $[0,L]$,
with zero flux at both ends, for $u$ and $Y=u+w$. The method has three ingredients:
a mesh with a fine band around the flame front
(Section~\ref{sec:meshes}), a linearly implicit time discretization that satisfies a
discrete comparison principle for every time step (Section~\ref{sec:theory}), and
a motion of the mesh that keeps the band on the front
(Section~\ref{sec:ale}). We assume (F) throughout.

\subsection{Space discretization and band meshes}\label{sec:meshes}

On a mesh $0=x_0<\dots<x_J=L$ with $h_i=x_{i+1}-x_i$ we use vertex-centred finite
volumes: with faces $x_{i\pm1/2}$ at the midpoints and
$|K_i|=(x_{i+1/2}^N-x_{i-1/2}^N)/N$,
\[
(\Delta_hu)_i=\frac1{|K_i|}\Big(x_{i+1/2}^{N-1}\frac{u_{i+1}-u_i}{h_i}
-x_{i-1/2}^{N-1}\frac{u_i-u_{i-1}}{h_{i-1}}\Big),
\]
with the obvious modification (zero flux) at $i=0$ and $i=J$. The discretization is
conservative and needs no special treatment at $r=0$, where the weight $r^{N-1}$
vanishes. First derivatives, which appear when the mesh moves, are discretized by
the three-point formula that is exact on quadratics,
\[
(D u)_i=\frac{h_{i-1}}{h_i(h_{i-1}+h_i)}(u_{i+1}-u_i)+\frac{h_i}{h_{i-1}(h_{i-1}+h_i)}(u_i-u_{i-1}),
\qquad 0<i<J .
\]
The matrix $-\Delta_h$ has positive diagonal, nonpositive off-diagonal entries and zero
row sums; the rows of $D$ also add up to zero.

The \emph{front} $x_f$ is the point where $u=\ve$; for $u\ge\ve$ there is no reaction.
A \emph{band mesh} has three zones:
\begin{itemize}
\item a \emph{band} around the front with uniform spacing $h_f$, which covers the
reaction layer ($5\ve$ on the fresh side and $20\ve$ on the burnt side, where $u$ decays
exponentially to $-\ve w_0$) plus a margin of $15\ve$ on each side;
\item a \emph{transition} on each side of the band, where the cells grow geometrically
by a factor $1.1$;
\item a \emph{coarse} zone with spacing $h_c$.
\end{itemize}
Our default values are $h_f=\ve/16$ and $h_c=0.01$. The number of nodes is
\[
\frac{55\,\ve}{h_f}+O\Big(\log\frac{h_c}{h_f}\Big)+O\Big(\frac L{h_c}\Big),
\]
which, for a fixed ratio $h_f/\ve$, grows only logarithmically as $\ve\to0$, and does not
depend on how far the front travels. The band of the analysis
(Theorem~\ref{cv:thm:main}) is wider than this one; see Section~\ref{an:sec:impl}.

\subsection{A monotone linearly implicit scheme}\label{sec:theory}

Let $x^n$ be the mesh at time $t_n=n\Delta t$ and let the nodes move with velocities
$\dot x^n_i$ during the step, with $\dot x^n_0=\dot x^n_J=0$, so that the new mesh is
$x^{n+1}_i=x^n_i+\Delta t\,\dot x^n_i$. Along the trajectory of a node,
$\frac{d}{dt}u(x_i(t),t)=u_t+\dot x_iu_x$, so \eqref{eq:Peps} becomes
$\frac{d}{dt}u=\Delta u+\dot x\,u_x-(u+w)f_\ve(u)$. Let
\[
\mathcal D^{n+1}=\Delta_h+\operatorname{diag}(\dot x^n)\,D\qquad\text{on the mesh }x^{n+1}.
\]
Given $(U^n,Y^n)$, the scheme (S) computes $(U^{n+1},Y^{n+1})$ from
\begin{align}
\frac{Y^{n+1}-Y^n}{\Delta t}&=\mathcal D^{n+1}Y^{n+1}-f_\ve(U^n)\,Y^{n+1},\label{eq:SY}\\
\frac{U^{n+1}-U^n}{\Delta t}&=\mathcal D^{n+1}U^{n+1}-f_\ve(U^n)\,Y^{n+1}.\label{eq:SU}
\end{align}
The reaction is linearly implicit: its rate is frozen at the previous time level
while the reactant is implicit. This is the standard linearization of sink terms that
preserves the positivity of concentrations
\cite{Patankar1980}; what we add is the observation that, when $f$ is
monotone, it gives a comparison principle for the coupled system
(Lemma~\ref{lem:comparison}). The scheme is written along the trajectories of the nodes,
in non-conservative form: $\Delta_h$ is conservative, but the convection term is not, so
the total amount of reactant is not conserved exactly when the mesh moves, and no
geometric conservation law is needed. Each step requires two tridiagonal solves. Subtracting
the equations, $W_h^n=Y^n-U^n$ solves the scheme without reaction, exactly as
$w^\ve=Y^\ve-u^\ve$ solves the heat equation; in particular, if $W_h^0\equiv\ve w_0$ is
constant, then $W_h^n\equiv\ve w_0$ for all $n$. With $\dot x\equiv0$ the scheme is
Eulerian. Occasionally the mesh is replaced by a new one after a step
(\emph{remeshing}); then $U^{n+1}$ and $Y^{n+1}$ are transferred by piecewise
\emph{linear} interpolation $\Pi$, whose weights are nonnegative and add up to one.

We assume that the mesh velocities are such that
\begin{equation}\label{eq:peclet}
\text{the off-diagonal entries of }-\mathcal D^{n+1}\text{ are nonpositive.}
\end{equation}
For $N=1$ these entries are
$-\frac{2+\dot x_ih_{i-1}}{h_i(h_{i-1}+h_i)}$ and $-\frac{2-\dot x_ih_i}{h_{i-1}(h_{i-1}+h_i)}$, so
\eqref{eq:peclet} is $-2\le\dot x_ih_{i-1}$ and $\dot x_ih_i\le2$, which is implied by the
cell P\'eclet condition $|\dot x_i|\max\{h_{i-1},h_i\}\le2$. In
radial symmetry the number $2$ in these two entries is replaced by $2\theta_i^+$ and
$2\theta_i^-$ respectively, where
$\theta_i^\pm=x_{i\pm1/2}^{N-1}(h_{i-1}+h_i)/(2|K_i|)=1+O(h/x_i)$, and \eqref{eq:peclet} reads
\[
-2\theta_i^+\le\dot x_i\,h_{i-1},\qquad \dot x_i\,h_i\le2\theta_i^- .
\]
The code checks \eqref{eq:peclet} whenever a matrix is assembled. A trial velocity that
violates it is discarded by the iteration of Section~\ref{sec:ale}, and if no admissible
velocity is found the run stops with an error; the time step is not reduced.
Since the rows of $\Delta_h$ and of $D$ add up to zero,
so do the rows of $\mathcal D^{n+1}$. All the statements below that hold ``for every time
step and every mesh'' are subject to \eqref{eq:peclet}. With $h_c=0.01$ it allows mesh
velocities up to about $200$; it becomes restrictive on coarser meshes when the front is
fast, as at the end of the extinction of a fresh pocket, where the speed of the front is
unbounded. An upwind discretization of the convection term outside the band would be
monotone without restriction, at the price of first order accuracy there; we have not
used it.

\begin{lemma}[Positivity and bounds]\label{lem:bounds}
Assume (F), \eqref{eq:peclet}, $\Delta t>0$ and $Y^0\ge0$. Then for all $n$:
\begin{enumerate}
\item (S) has a unique solution, and $Y^n\ge0$;
\item $\min W_h^0\le W_h^n\le\max W_h^0$;
\item $-\max W_h^0\le U^n\le\max U^0$.
\end{enumerate}
In particular the bounds are uniform in $\ve$, in the mesh and in $\Delta t$.
\end{lemma}

\begin{proof}
For $F\ge0$ the matrix $\mathsf A_F=\Delta t^{-1}I-\mathcal D^{n+1}+\operatorname{diag}F$
has nonpositive off-diagonal entries and row sums $\Delta t^{-1}+F_i>0$, so it is
strictly diagonally dominant with positive diagonal, hence a nonsingular
M-matrix: $\mathsf A_F^{-1}\ge0$. By \eqref{eq:SY} with $F=f_\ve(U^n)$,
$Y^{n+1}=\mathsf A_F^{-1}Y^n/\Delta t\ge0$. Since $\mathsf A_0\mathbf 1=\Delta t^{-1}\mathbf 1$, the solution of
$\mathsf A_0W_h^{n+1}=W_h^n/\Delta t$ satisfies the maximum principle. From \eqref{eq:SU} and
$f_\ve(U^n)Y^{n+1}\ge0$, $\mathsf A_0U^{n+1}\le U^n/\Delta t$, hence $\max U^{n+1}\le\max U^n$;
and $U^{n+1}=Y^{n+1}-W_h^{n+1}\ge-\max W_h^0$. The interpolation $\Pi$ preserves all
these bounds.
\end{proof}

We say that a sequence $V=(V^n)$ on the same meshes is a \emph{discrete
supersolution} if $V^{n+1}+W_h^{n+1}\ge0$ and
\[
\frac{V^{n+1}-V^n}{\Delta t}-\mathcal D^{n+1}V^{n+1}+f_\ve(V^n)\,(V^{n+1}+W_h^{n+1})\ \ge\ 0
\]
for every $n$ (and $V$ after a remeshing is at least $\Pi V$), and a \emph{discrete
subsolution} if the reverse inequalities hold; no sign condition is needed for
subsolutions.

\begin{lemma}[Discrete comparison principle]\label{lem:comparison}
Assume (F), \eqref{eq:peclet} and $Y^0\ge0$, and let $U$ be the solution of (S). If $V$ is
a discrete supersolution with $V^0\ge U^0$, then $V^n\ge U^n$ for all $n$. If $V$ is a
discrete subsolution with $V^0\le U^0$, then $V^n\le U^n$ for all $n$. This holds for
every $\Delta t>0$ and every sequence of meshes.
\end{lemma}

\begin{proof}
Assume $U^n\le V^n$ and let $Z=U^{n+1}-V^{n+1}$, $F_U=f_\ve(U^n)$, $F_V=f_\ve(V^n)$.
Since $f_\ve$ is nonincreasing, $F_U\ge F_V$. Subtracting the inequality for $V$ from
\eqref{eq:SU} (recall $Y=U+W_h$),
\[
\Big(\frac1{\Delta t}-\mathcal D^{n+1}\Big)Z+F_U(U^{n+1}+W_h^{n+1})-F_V(V^{n+1}+W_h^{n+1})
\le\frac{U^n-V^n}{\Delta t}\le0 .
\]
The reaction terms equal $F_UZ+(F_U-F_V)(V^{n+1}+W_h^{n+1})\ge F_UZ$, so
$\mathsf A_{F_U}Z\le0$, and $Z\le0$ because $\mathsf A_{F_U}^{-1}\ge0$. For subsolutions, with
$Z'=V^{n+1}-U^{n+1}$ the reaction terms are
$F_VZ'+(F_V-F_U)(U^{n+1}+W_h^{n+1})\ge F_VZ'$, because now $F_V\ge F_U$ and
$U^{n+1}+W_h^{n+1}=Y^{n+1}\ge0$ by Lemma~\ref{lem:bounds}. Finally $\Pi$ is order preserving.
\end{proof}

Lemma~\ref{lem:comparison} is the discrete form of the comparison principle for
\eqref{eq:Peps}, on which the barrier arguments of
\cite{LedermanVazquezWolanski2001,Petrosyan2002,FernandezBonderWolanski2004} and of
Section~\ref{sec:rate} rest. It does not require $\Delta t\lesssim\ve^2$; a fully
implicit treatment of $(u+w)f_\ve(u)$ would, because that function is not
monotone in $u$. It reduces error estimates to the construction of discrete
barriers, which is the subject of Section~\ref{sec:analysis}.

\subsection{Following the flame: the ALE version}\label{sec:ale}

On a mesh that is fixed between remeshings the front crosses $O(\ve^{-1})$ fine cells
per unit time, and the lag of the reaction rate in (S) produces an error of order
$\Delta t/\ve$ (Sections~\ref{sec:analysis} and \ref{sec:experiments}). We avoid it by letting the
band travel with the front. In an \emph{arbitrary Lagrangian--Eulerian} (ALE)
formulation the computational nodes move with a velocity that is chosen freely: it
is neither zero (Eulerian description, fixed nodes) nor the velocity of the
material (Lagrangian description), and the equation is written along the node
trajectories, so that the motion of the nodes appears as a convection term
\cite{HuangRussell2011}; this is the term $\operatorname{diag}(\dot x)D$ in (S). Here the node
velocity is chosen so that the band follows the flame.

\emph{Mesh velocity.} During the step from $t_n$ to $t_{n+1}$ node $i$ moves with
velocity $\dot x^n_i=v_n\chi(x^n_i)$, where $v_n$ is a scalar and $\chi=1$ within a distance
$35\ve+\ell_0$ of the reference point of the band (the \emph{rigid zone}, which contains the
band), and $\chi$ decays linearly to $0$ over a further distance $\ell_0$; nodes far from the front, and the ends
of $[0,L]$, do not move. (Near the ends $\chi$ is multiplied by a linear cut-off, so
that no node leaves $[0,L]$.) The band therefore translates rigidly. By default
$\ell_0=\max\{2,3/c\}$, where $c$ is an estimate of the front speed: the rigid zone
must contain several decay lengths $1/c$ of the profile ahead of the front
(Section~\ref{an:sec:impl}).

\emph{Implicit mesh velocity.} The band has a reference point $x_c$, initially the
front. At each step $v_n$ is chosen so that the computed front ends the step at
the reference point: $x_f^{n+1}(v_n)=x_c+\Delta t\,v_n$, and then $x_c$ is advanced by
$\Delta t\,v_n$. This scalar equation is solved by the secant method, starting from
$v_{n-1}$, until the front is within $10^{-3}h_f$ of the reference point. Each evaluation
is a step of (S), that is, two tridiagonal solves; the iteration takes at least three
evaluations and at most ten, after which the last iterate is accepted and the event is
counted. For travelling waves it takes three evaluations per step, and for the
accelerating fronts of Section~\ref{sec:exp:thm} between $3$ and $4.2$ on average,
depending on the time step. The tolerance allows an error of $10^{-3}h_f/\Delta t$ in the
velocity of a single step, which does not accumulate, because the equation is for the
position of the front. The same idea underlies the freezing method
for travelling waves \cite{BeynThummler2004}. An explicit choice of $v_n$ (the front
speed measured in the previous steps) lags behind accelerating fronts by
(acceleration)$\times\Delta t$, which is too much for small $\ve$; the implicit choice
removes this lag.

\emph{Remeshing.} The cells of the decay zones are stretched or compressed by the
motion. The mesh is rebuilt when it has travelled a distance $\ell_0/4$, or if the front
drifts from the reference point by more than $15\ve/4$, or if a cell has been compressed
below $h_f/2$; then $x_c$ is reset to the
front. (A trial velocity for which two nodes would cross is discarded.) The new band is placed on the lattice of the old band, so the transfer by $\Pi$
is exact in the band, and interpolation only acts in the transition and coarse zones,
where the solution is smooth on the scale of the mesh.

Since $\Pi$ is order preserving, and \eqref{eq:peclet} is enforced,
Lemmas~\ref{lem:bounds} and \ref{lem:comparison} hold for the ALE scheme, for every time
step. For a travelling wave the solution is steady in the frame of the band, and
the error of the ALE scheme does not contain $\Delta t$ (Theorem~\ref{an:thm:ale}). For
general data the error contains $\Delta t$, but not $\Delta t/\ve$ (Theorem~\ref{cv:thm:main}
and Section~\ref{sec:experiments}).

\subsection{Eulerian version, reference solver and initial data}\label{sec:eulerian}

For comparison we also use two solvers on meshes that do not move.

\emph{Eulerian (S).} The scheme (S) with $\dot x\equiv0$ and a fixed time step, on a
band mesh that is rebuilt around the current front after the time in which the front,
at its measured speed, travels half the margin of the band ($7.5\,\ve$), and in any case
at least $20$ times in a run, with transfer by $\Pi$. It is covered by
Lemmas~\ref{lem:bounds} and \ref{lem:comparison}.

\emph{BDF.} A method-of-lines solver: the finite volume discretization on the same
band meshes, rebuilt in the same way, integrated with the variable order, variable
step BDF method of SciPy with the analytic sparse Jacobian
($\mathrm{rtol}=10^{-5}$, $\mathrm{atol}=10^{-4}\ve$), and transfer by monotone cubic
(PCHIP) interpolation. It is not monotone and is not covered by the analysis; we use
it as an accurate reference. Its time step must resolve the passage of the front
through each fine cell, so the number of steps grows faster than $\ve^{-1}$ (like
$\ve^{-1.4}$ in the experiments of Section~\ref{sec:exp:tw}).

\emph{Initial data.} Given the initial datum $u_0\ge0$ of \eqref{eq:P} and $w_0$ we take
\begin{equation}\label{eq:initdata}
u_0^\ve=u_0-\ve w_0\Big(1-\frac{u_0}\ve\Big)^+,\qquad Y_0^\ve=u_0^\ve+\ve w_0,
\end{equation}
so that $(Y_0^\ve-u_0^\ve)/\ve=w_0$ exactly, $Y_0^\ve\ge0$, $Y_0^\ve=0$ in the burnt region
$\{u_0=0\}$, and $u_0^\ve\to u_0$ uniformly with uniformly bounded Lipschitz constants, as
required in \cite{FernandezBonderWolanski2004}; these data satisfy \eqref{q:eq:data}.

\section{Convergence of the scheme to the free boundary problem}\label{sec:analysis}

The barriers of Section~\ref{sec:rate} are explicit functions, so they can be
evaluated on the mesh. In this section we show that, with a dilation that depends on
time and a damping of the inner profile that is proportional to the truncation errors
in the layer, they are discrete super- and subsolutions of the scheme (S). By the
discrete comparison principle this gives the convergence of the fully discrete solution
to the solution of the free boundary problem \eqref{eq:P}, for radial data and a
perturbation $w_0$ that need not be constant, with
an explicit error bound (Theorem~\ref{cv:thm:main}). For
travelling waves on meshes that move with the wave we also prove an estimate in which
the time step does not appear at all (Theorem~\ref{an:thm:ale}).

\subsection{Setting and main result}\label{cv:sec:main}

Let $N\in\{1,2,3\}$ and assume (F$_2$). Let $w_0=w_0(x,t)$ be a radial solution of the heat
equation in $\{|x|<L\}\times(0,T^*]$ with zero flux on $\{|x|=L\}$ (for $N=1$, an even
function of $x$), and let $0<L_w\le L$, such that
\begin{itemize}
\item[(W1)] $-1+\delta_1\le w_0\le\bar W$, and $D_x^3w_0$, $\partial_t\nabla w_0$, $\partial_t^2w_0$ are
bounded;
\item[(W2)] $x\cdot\nabla w_0+2t\,\partial_tw_0\ge0$ in $\{|x|\le L_w\}\times[0,T^*]$.
\end{itemize}
Condition (W2) says that $\lambda\mapsto w_0(\lambda x,\lambda^2t)$ is nondecreasing while $(\lambda x,\lambda^2t)$ stays
in that set; it is the monotonicity in \eqref{q:eq:mono}, in the region where it will be
used. A constant $w_0$ satisfies (W1), (W2) with $L_w=L$. Let $w^\ve=\ve w_0$, and let $v$ be a
radial regular classical solution on $[0,T^*]$, in the sense of Section~\ref{sec:rate}, of
\eqref{eq:P} with this $w_0$, that is, with $|\nabla v|=\sqrt{2M(w_0(x,t))}$ on the free boundary,
with $u_0=v(\cdot,0)$, such that
\begin{itemize}
\item[(V1)] $\{v(\cdot,t)>0\}=\{r<R(t)\}$ with $0<R_{\min}\le R(t)\le L_w-\varrho_L$, $\varrho_L>0$; $u_0$ is
nonincreasing in $r$; and $R(t)-2t\dot R(t)\ge c_R>0$ on $[0,T^*]$;
\item[(V2)] $D_x^3v$, $\partial_t\nabla v$ and $\partial_t^2v$ are bounded in $\{v>0\}$, and the
extension $Ev$ of (R1) is $C^{1,1}$ jointly in $(x,t)$: $\partial_t\nabla Ev$ and $\partial_t^2Ev$
are bounded too.
\end{itemize}
For $N=1$ this is a slab $\{|x|<R(t)\}$ of fresh gas, computed on $[0,L]$. By (V1), (V2) and
(R2), $R\in C^{1,1}$ and (H) holds. The last condition in (V1) holds, in particular, when
the pocket does not expand, $\dot R\le0$; it is used through the following consequence,
which is a strict form of the monotonicity of $v$ under parabolic dilations.

\begin{lemma}\label{cv:lem:dil}
Under (V1), (V2) there is $c_v>0$ such that
$x\cdot\nabla v+2t\,v_t-\frac12v\le-c_v$ in $\overline{\{v>0\}}$.
\end{lemma}

\begin{proof}
The function $\psi=x\cdot\nabla v+2t\,v_t-\frac12v$ is caloric in $\{v>0\}$, because $v$ is, and by
(V2) it is continuous up to the boundary. On the free boundary $v=0$,
$x\cdot\nabla v=-R|v_r|$ and, differentiating $v(R(t),t)=0$, $v_t=|v_r|\dot R$; hence
$\psi=-|v_r|(R-2t\dot R)\le-2\delta_0c_R$ there, by (R2). At $t=0$,
$\psi=x\cdot\nabla u_0-\frac12u_0\le-\frac12u_0$ in $\overline{\Omega_0}$, which is negative in $\Omega_0$, and
$\psi=-R|\partial_ru_0|<0$ on $\partial\Omega_0$; so $\psi(\cdot,0)\le-c$ in $\overline{\Omega_0}$. The claim follows
from the maximum principle in the non-cylindrical domain $\{v>0\}$.
\end{proof}

The scheme (S) is run on $[0,L]$ with $U^0_j=u_0^\ve(x^0_j)$ and
$Y^0_j=U^0_j+\ve w_0(x^0_j,0)$, where $u_0^\ve$ is given by \eqref{eq:initdata} with
$w_0(\cdot,0)$; these data satisfy \eqref{q:eq:data}. We may and do assume that
$r\le\varrho_L/4$ in (R1), (R2). We make the following
assumptions on the meshes, in which $h_f\le h_c\le1$, $\ell_b>0$ and $v_{\max}$ are given.
\begin{itemize}
\item[(M1)] There is no remeshing, \eqref{eq:peclet} holds, all the spacings are at most
$h_c$, and $|\dot x^n_j|\le v_{\max}$.
\item[(M2)] (\emph{Band.}) For every $n$, every node of $x^{n+1}$ with
$|x-R(t_{n+1})|\le\ell_b$ has its two neighbours at distance $h_f$, and during the step
the three of them moved with the same velocity $v_n$. We call these the \emph{band
nodes} of the step.
\end{itemize}
Nothing else is assumed outside the band. Let
\[
\mu=\max_n|v_n-\dot R(t_n)|,\qquad
\mathcal E=\ve\log\frac1\ve+\frac{h_f}\ve+h_c+\Delta t+\mu\,\min\Big\{1,\frac{\Delta t}\ve\Big\}.
\]
In (M2), $\ell_b<\min\{R_{\min},\varrho_L\}$. A uniform mesh that does not move satisfies
(M1), (M2) with $h_c=h_f$, $v_n=0$ and $\ell_b=\min\{R_{\min},\varrho_L\}/2$, and then
$\mu=\max|\dot R|$.

\begin{theorem}[Convergence of the scheme]\label{cv:thm:main}
Assume (F$_2$), (W1), (W2), (V1), (V2), (M1), (M2), and let $T<T^*$. There are positive
constants $\ve_0,c_0,C_b,C$, depending only on $f,w_0,N,v,T,T^*,L$ and $v_{\max}$, such that, if
$0<\ve\le\ve_0$ and
\begin{equation}\label{cv:eq:hyp}
\mathcal E\le c_0,\qquad \Delta t\,(\mu+\mathcal E)\le c_0\,\ve,\qquad \ell_b\ge C_b\,\mathcal E,
\end{equation}
then the solution of (S) satisfies
\[
\big|U^n_j-v(x^n_j,t_n)\big|\le C\,\mathcal E\qquad\text{for all }j\text{ and all }t_n\le T.
\]
\end{theorem}

\begin{remark}\label{cv:rem:main}
\begin{enumerate}
\item \emph{Meshes that do not move.} Here $\mu$ is of order one and the bound is
$C(\ve\log\frac1\ve+\frac{h_f+\Delta t}{\ve}+h_c)$: the time step is divided by $\ve$. This is
the lag of the reaction rate $f_\ve(U^n)$ on a front that crosses the cells, and it is
observed (Sections~\ref{sec:exp:tw} and \ref{sec:exp:thm}).
\item \emph{Meshes that follow the front.} Since $\mu\min\{1,\Delta t/\ve\}\le\mu$, the bound is
at most $C(\ve\log\frac1\ve+\frac{h_f}{\ve}+h_c+\Delta t+\mu)$: when $\mu$ is small the time step
enters as $\Delta t$, not as $\Delta t/\ve$.
\item \emph{The time step.} The second condition in \eqref{cv:eq:hyp} says that in one
step the mesh and the fronts of the barriers, whose velocities differ from $\dot R$ by a
multiple of $\mathcal E$, move with respect to each other by a small fraction of the width
$\ve$ of the layer. Since $\mathcal E\ge h_f/\ve+h_c+\Delta t$, it implies
$\Delta t\le c_0\min\{\ve^2/h_f,\ve/h_c\}$ and $\Delta t^2\le c_0\ve$: for a fixed ratio $h_f/\ve$ the theorem
covers time steps up to a multiple of $\ve$, and larger ones only on finer meshes.
\item \emph{Mesh sizes.} Negative powers of $\ve$ multiply only the mesh size of the band,
through the ratio $h_f/\ve$: for $\Delta t\lesssim\ve$ it is enough to resolve the layer, and
$h_f\ll\ve^2$ is not needed; for a fixed ratio $h_f/\ve$ the bound is uniform in $\ve$. On
the other hand, the band has to contain the fronts of the two barriers, so its width
must exceed a multiple of the bound itself (for a uniform mesh, this is a further
smallness condition on $\mathcal E$); and we use the pointwise truncation error of $\Delta_h$,
which is of first order on arbitrary meshes.
\item \emph{The mesh velocity.} In the theorem the velocities $v_n$ are given. In the
method of Section~\ref{sec:ale} they are computed, and $\mu$ is an output: for the
method itself the estimate is a posteriori. See Section~\ref{an:sec:impl}.
\item \emph{The hypotheses on $w_0$ and $v$.} As usual in error estimates, the solution
$v$ of the limit problem is assumed to be regular; (V2) includes compatibility conditions
on the data at $t=0$. The radial self-similar solutions of
Section~\ref{sec:selfsimilar}, for which $w_0$ is constant, satisfy (V1), (V2) on $[0,T^*]$ for
every $T^*$ smaller than their extinction time $T_e$: their radius is $\varrho\sqrt{T_e-t}$, so
that $R-2t\dot R=\varrho\,T_e/\sqrt{T_e-t}$. They also give an example with a non-constant
perturbation: if $V$ is the self-similar solution for a constant $\bar w$, continued to
$r>R(t)$ by the same formula, and $b<0$ is small, then $w_0=\bar w+b\,V$ is caloric, equals
$\bar w$ on the free boundary, so that $V$ solves \eqref{eq:P} with this $w_0$, and satisfies
(W1), (W2), because $x\cdot\nabla V+2t\,V_t\le0$; only the zero-flux condition at $|x|=L$ fails.
Functions $w_0$ that satisfy (W1), (W2) and that condition are also easy
to give: for $N=3$, $w_0=\bar w-b\,e^{-k^2t}\sin(kr)/(kr)$ with $b>0$ and $\tan(kL)=kL$ has zero
flux at $r=L$ and satisfies (W2) with $L_w=\pi/k\approx0.7\,L$. In a bounded domain (W2) is not
inherited from $t=0$, so it is a condition for all times. Finally, when $w_0$ is constant
the last condition in (V1) can be replaced by the assumption that $v$ is nonincreasing
in $r$: in Step~1 below one takes $\tau=\int_0^t\lambda^2$ instead of $\tau=\lambda^2t$, and then
$\tilde u_t-\Delta\tilde u=\frac{\dot\lambda}{\lambda'}\big[\frac{y\cdot\nabla v}\lambda-\frac v{2\lambda'}\big]$, whose bracket is negative; this is
not compatible with (W2) when $w_0$ depends on $t$.
\end{enumerate}
\end{remark}

\subsection{Proof of Theorem~\ref{cv:thm:main}}

We write $W_h^n=Y^n-U^n$ for the discrete perturbation,
$\Lop[V]^n_j=\Delta t^{-1}(V^{n+1}_j-V^n_j)-(\mathcal D^{n+1}V^{n+1})_j
+f_\ve(V^n_j)(V^{n+1}_j+W^{n+1}_{h,j})$, and $\Lop_0[V]$ for the same expression without the last
term. For a function $V(x,t)$ we set $V^n_j=V(x^n_j,t_n)$. In a step $n\to n+1$ and for a
node $j$ we write $x'=x^{n+1}_j$ and $x_j(t)=x'-(t_{n+1}-t)\dot x^n_j$, and we call
\emph{evaluation points} of the node the points $(x,t_{n+1})$ with
$x^{n+1}_{j-1}\le x\le x^{n+1}_{j+1}$, and the points $(x_j(t),t)$ with $t_n\le t\le t_{n+1}$.
Quantities evaluated at $(x',t_{n+1})$ carry a prime, and those evaluated at $(x^n_j,t_n)$
a subscript $j$. Unless otherwise stated, $C,c$ denote constants that depend only on the
quantities listed in the theorem. The discrete $W_h$ is not $\ve w_0$, but it is close to it.

\begin{lemma}\label{cv:lem:heat}
Let $w_0$ and the initial data of the scheme be as in Section~\ref{cv:sec:main}, and
assume (W1) and (M1). Then
\[
\varpi:=\max_{t_n\le T}\max_j\Big|\frac{W^n_{h,j}}\ve-w_0(x^n_j,t_n)\Big|\le C_wT\,(\Delta t+h_c),
\]
with $C_w$ depending only on the bounds in (W1), on $N$, $L$ and $v_{\max}$.
\end{lemma}

\begin{proof}
$W_h$ solves the scheme without reaction, $\Lop_0[W_h]=0$, and $W_h^0=\ve w_0(\cdot,0)$ on
the mesh. By (W1), $|\Lop_0[w_0]^n_j|\le C_w(\Delta t+h_c)$: the difference quotient in time along
the path of a node is $(\partial_t+\dot x_j\partial_x)w_0$ at $(x',t_{n+1})$ up to
$\frac{\Delta t}2\sup|(\partial_t+\dot x_j\partial_x)^2w_0|$, $D$ is exact on quadratic functions, and $\Delta_h$ is
consistent of first order on arbitrary meshes for radial functions (even functions, if
$N=1$) whose third derivatives in Cartesian coordinates are bounded; this includes the
rows next to $r=0$, the row $j=0$, which is $2N(u_1-u_0)/h_0^2$ and is exact on quadratic
functions of $r$, and the row $j=J$, because $\partial_rw_0=0$ at $r=L$. Hence $Z=W_h/\ve-w_0$
satisfies $|\Lop_0[Z]|\le C_w(\Delta t+h_c)$ and $Z^0=0$, and, since $\mathsf A_0^{-1}\ge0$ and
$\mathsf A_0\mathbf 1=\Delta t^{-1}\mathbf 1$,
$\max_j|Z^{n+1}_j|\le\max_j|Z^n_j|+C_w\Delta t(\Delta t+h_c)$.
\end{proof}

We set $\varpi'=\varpi+C\Delta t$, which bounds $|W_{h,j}/\ve-w_0'|$ at every node, because $w_0$ is
Lipschitz along the paths of the nodes. By the lemma and \eqref{cv:eq:hyp},
$\varpi'\le C\mathcal E\le Cc_0$.

\emph{Step 1. Dilations that depend on time.} For $\beta,a>0$ let
$\lambda(t)=1-a\,\mathcal E\,e^{\beta t}$ or $\lambda(t)=1+a\,\mathcal E\,e^{\beta t}$, and
\[
\tilde u(x,t)=\frac1{\lambda'(t)}\,v\big(\lambda(t)x,\lambda(t)^2t\big),\qquad
\lambda'=\frac{1+\lambda}2,
\]
which is defined on $[0,T]$, with $\tilde R\le L_w-\varrho_L/2$ below, if $c_0$ is small (depending
on $a$ and $\beta$), because then $\lambda^2T\le T^*$. With $y=\lambda x$ and $\tau=\lambda^2t$,
\[
\tilde u_t-\Delta\tilde u=\frac{\dot\lambda}{\lambda'}\Big[\frac{y\cdot\nabla v+2\tau v_t}{\lambda}-\frac v{2\lambda'}\Big](y,\tau).
\]
By Lemma~\ref{cv:lem:dil} the bracket is at most $-c_v/2$ in $\{v>0\}$ when $|1-\lambda|$ is small.
Hence, for the first
choice of $\lambda$, $\tilde u_t-\Delta\tilde u\ge\frac{c_v}4|\dot\lambda|$ in $\{\tilde u>0\}$, and for the second
one $\tilde u_t-\Delta\tilde u\le-\frac{c_v}4|\dot\lambda|$, where $|\dot\lambda|\ge\beta a\mathcal E$. The free boundary
of $\tilde u$ is $r=\tilde R(t)=R(\lambda^2t)/\lambda$, and there
$|\nabla\tilde u|=\frac\lambda{\lambda'}\sqrt{2M(w_0(\lambda x,\lambda^2t))}$. By (W2) and the monotonicity of $M$, for
$\lambda<1$ this is at most $\frac\lambda{\lambda'}\sqrt{2M(w_0(x,t))}\le\sqrt{2M(w_0(x,t))}-c\,(1-\lambda)$, and for
$\lambda>1$ it is at least $\sqrt{2M(w_0(x,t))}+c\,(\lambda-1)$ and, since $w_0$ is Lipschitz in $(x,t)$,
at most $\sqrt{2M(w_0(x,t))}+C(\lambda-1)$. So
$\tilde u$ is a regular strict supersolution (resp.\ subsolution) with margins
$(\delta,\delta_0,r,K_0)$, $\delta=c\,a\mathcal E$, uniformly in $t$ (for the subsolution, with
$C_1=Ce^{\beta T}$); this is the reason for the choice of the time $\lambda^2t$, for which the
dilation acts on $w_0$ as in \eqref{q:eq:mono}. By (V2),
$E\tilde u(x,t)=Ev(\lambda x,\lambda^2t)/\lambda'$ is $C^{1,1}$ in $(x,t)$, and so are
$\hat E=\psi(E\tilde u)$ and $\tilde v^\pm=\hat E\pm\Lambda'(\hat E^-)^2$ in the proof of
Lemma~\ref{q:lem:ext}: in addition to the bounds of that lemma,
$|\partial_t\nabla\tilde v^\pm|+|\partial_t^2\tilde v^\pm|\le K_1$. (The derivatives of $\lambda$ are bounded by
$C_\beta a\mathcal E\le1$ for $c_0$ small, so these constants do not depend on $a$.) Moreover
$|\tilde u-v|\le C_\beta\,a\mathcal E$, where $C_\beta$ depends also on $\beta$, and
\[
\dot{\tilde R}=\lambda\dot R(\tau)-\frac{\dot\lambda}{\lambda^2}\big(R(\tau)-2\tau\dot R(\tau)\big),\qquad
|\lambda\dot R(\tau)-\dot R(t)|\le C_R|1-\lambda(t)| ,
\]
because $|\tau-t|\le3t|1-\lambda(t)|$. By (V1), $R-2\tau\dot R\ge c_R$; we fix $\beta$ with
$\beta c_R/4\ge C_R+2$. Then
\begin{equation}\label{cv:eq:speed}
\pm\big(\dot{\tilde R}-\dot R\big)(t)\ge2|1-\lambda(t)|\ge2a\mathcal E,\qquad|\dot{\tilde R}-\dot R|\le C_\beta\,a\mathcal E,
\end{equation}
with the sign $+$ for the supersolution: its front moves outwards with respect to the
front of $v$, and the front of the subsolution inwards.

\emph{Step 2. The barriers and the parameters.} Let $\tilde q=\max\{1,\sup|w_0|\}$,
$s=\tilde v^\pm/\ve$, and
\[
V=-W_h+\ve\,\Theta^+(s;w_0)+\ve^2\tilde q,\qquad
\underline V=-W_h+Z,\quad Z=\max\big\{\ve\,\Theta^-(s;w_0),\,-\ve^2\tilde q\big\},
\]
where $W_h$ is the discrete function and the profiles are evaluated with the exact
$w_0$ at the node, $w_0=w_0(x^n_j,t_n)$. These are the barriers of
Propositions~\ref{q:prop:super} and \ref{q:prop:sub} with $W^\ve=w_0$, $\gamma=c_1\delta$,
$\bar q=\ve\tilde q$ and $\kappa=K\ve$, except that the offset is $W_h$ instead of $\ve w_0$. The
profiles carry the exact $w_0$, so that their derivatives are controlled by
Lemma~\ref{q:lem:inner}; the offset is discrete, so that $\Lop_0[W_h]=0$ exactly and
\begin{equation}\label{cv:eq:split}
\Lop[V]^n_j=\Lop_0[\mathcal V]^n_j+f_\ve(V_j)\,(\mathcal V'+\ve^2\tilde q),\qquad
\mathcal V=\ve\,\Theta^+(s;w_0),\qquad
f_\ve(V_j)=\frac1{\ve^2}f\Big(\Theta_j+\ve\tilde q-\frac{W_{h,j}}\ve\Big),
\end{equation}
and similarly $\Lop[\underline V]^n_j=\Lop_0[Z]^n_j+f_\ve(\underline V_j)\,Z'$. Here $\Theta^-$ is continued to
$s<s_*$ by the solution of
the same equation down to $\Theta^-=-1$, and by $-1$ further on. On the continuation the
bounds $c\sqrt\gamma\le\Theta^-_s\le C$ and $|\Theta^-_{ss}|\le C\Theta^-_s$ and those of
Lemma~\ref{q:lem:inner}(b)(ii) still hold, by the proof of that lemma, which covers
$-1\le\Theta\le2+W$ (the bound $|\Theta^-_{ss}|\le C(\bar q+\kappa\Theta^-_s)$ does not, and it
is not used there), and its length in $s$ is at least $c_e>0$; the margin
\eqref{q:eq:marginsub} also holds there, since its proof only uses that
$d\le C\ve\log(1/\gamma)\le r_1$ and Lemma~\ref{q:lem:ext}(iii). The new
point is the size of the damping: we take
\begin{equation}\label{cv:eq:K}
K=C_K\Big(\log\frac1\ve+\frac{h_f}{\ve^2}+\frac{\Delta t+h_c}\ve+\frac{\mu_*\Delta t}{\ve^2}\Big),\qquad
\mu_*=\begin{cases}\mu&\text{if }\Delta t<\ve,\\ 0&\text{if }\Delta t\ge\ve,\end{cases}
\end{equation}
so that $K\ve\le5C_K\mathcal E$.
Lemma~\ref{q:lem:inner} requires $\kappa\le c\gamma$, that is,
$5C_K\mathcal E\le c\,c_1\delta$, which holds if $a$ is large; the conditions $\delta\le\delta_*$ and
$\ve\le c_*\delta$ of the margins hold too, because $\delta=c\,a\mathcal E\ge c\,a\,\ve\log\frac1\ve$, for
$a$ large and then $c_0$ and $\ve_0$ small; this is where the truncation errors of
the layer enter the error bound. From Section~\ref{sec:rate} we only take $\Lambda$, $c_1$
and the thresholds $\delta_*$, $c_*$ under which the margins \eqref{q:eq:marginsup} and
\eqref{q:eq:marginsub} hold (with $\eta_\ve=0$, and $C_1$ fixed by $\beta$); these do not involve
the damping.
The further restrictions on $c_*$ and $\delta_*$ in Propositions~\ref{q:prop:super} and
\ref{q:prop:sub}, which depend on $K$, are not used: here $\kappa\le c\gamma$ is obtained by
taking $a$ large, and the term that they absorb in Proposition~\ref{q:prop:sub} is
estimated in a different way in Step~6. The constants are chosen in the order $\beta$;
$\Lambda$, $c_1$, $\delta_*$, $c_*$; $C_S$, $C_K$; $a$;
$c_0$, $\ve_0$, $C_b$. Recall that $\Theta^+$ is convex in $s$, so $\Theta^+\ge s+w_0$, with equality for
$s\ge s_m$; that $\Theta^-$ is convex on $[s_c,2]$ with $s_c<1$, so $\Theta^-\ge s+w_0$ for $s\ge1$, with
equality for $s\ge2$; and that $c\,\Theta^+\le\Theta^+_s\le C\Theta^+$ for $s\le s_m$, by \eqref{q:eq:Plin}.

\emph{Step 3. Elementary facts.} (i) If $g\in C^{1,1}([t_n,t_{n+1}])$ then
$g'(t_{n+1})-\frac{\Delta t}2\sup(g'')^+\le\Delta t^{-1}(g(t_{n+1})-g(t_n))\le
g'(t_{n+1})+\frac{\Delta t}2\sup(g'')^-$. (ii) On a uniform stencil of size $h$, if
$V$ is $C^{1,1}$,
\[
\delta^2V(x)=\int_{-h}^h\frac{h-|y|}{h^2}V_{xx}(x+y)\,dy,\qquad|DV(x)-V_x(x)|\le\frac h2\sup|V_{xx}| .
\]
On an arbitrary mesh $|\Delta_hV_j|\le C(\sup|V_{xx}|+\sup|V_x|)$ and $|DV_j|\le\sup|V_x|$ if
$x_j\ge R_{\min}/2$ and $j<J$, the suprema being taken on the stencil; and for $N\ge2$, on a uniform
stencil with $x_j\ge R_{\min}/2$,
$\Delta_h=(1+\epsilon_j)\delta^2+(\frac{N-1}{x_j}+\epsilon_j')D$ with $|\epsilon_j|+|\epsilon_j'|\le Ch^2$. (iii) If all the
evaluation points of a node lie in $\overline{\{\tilde u>0\}}$, then by (V2), as in the proof of
Lemma~\ref{cv:lem:heat},
\begin{equation}\label{cv:eq:outer}
\Lop_0[\tilde u]^n_j=(\tilde u_t-\Delta\tilde u)(x',t_{n+1})+\vartheta,\qquad|\vartheta|\le C_\beta(\Delta t+h_c).
\end{equation}
By Step~1, if $a$ is large, $\Lop_0[\tilde u]\ge\frac{c_v}8\beta a\mathcal E$ for the supersolution and
$\Lop_0[\tilde u]\le-\frac{c_v}8\beta a\mathcal E$ for the subsolution; and
$\ve|\Lop_0[w_0]|\le C\ve(\Delta t+h_c)$ is smaller than half of these quantities. (iv) Since
$|\Theta^\pm_{ss}|\le C\Theta^\pm_s$, we have
$\Theta^\pm_s(s_1)\le e^{C|s_1-s_2|}\Theta^\pm_s(s_2)$ at a fixed $W$; and $(1+|s|)\Theta^+\le C$ for $s\le1$,
$(1+|s|)|\Theta^-|\le C$ for $s_*\le s\le2$. (v) Let $\mathcal V=\ve\,\Theta^\pm(\tilde v^\pm/\ve;w_0)$. A.e.,
\begin{equation}\label{cv:eq:Vder}
\mathcal V_x=\Theta_s\tilde v_x+\ve\Theta_W\partial_xw_0,\qquad
\mathcal V_{xx}=\frac1\ve\Theta_{ss}\tilde v_x^2+\Theta_s\tilde v_{xx}+\mathcal P_x,
\end{equation}
with $\mathcal P_x=2\Theta_{sW}\tilde v_x\partial_xw_0+\ve\Theta_{WW}(\partial_xw_0)^2+\ve\Theta_W\partial_x^2w_0$, and, by
Lemma~\ref{q:lem:inner}(ii), where $|s|\le C\log\frac1\ve$,
\begin{equation}\label{cv:eq:Px}
|\mathcal P_x|\le C_3(1+|s|)|\Theta_{ss}|+C\Theta_s,\qquad|\ve\Theta_W|\le C\ve\log\tfrac1\ve\,\Theta_s .
\end{equation}
In the same way, along the path of a node, $g(t)=\mathcal V(x_j(t),t)$ satisfies
$g'=\Theta_s\dot{\tilde v}+\ve\Theta_W\dot w_0$ and
$g''=\ve^{-1}\Theta_{ss}\dot{\tilde v}^2+\Theta_s\ddot{\tilde v}+\mathcal P_t$ with
$|\mathcal P_t|\le C\big((1+|s|)|\Theta_{ss}|+\Theta_s\big)(1+|\dot{\tilde v}|)$, where the dots are derivatives along
the path and $\dot w_0$, $\ddot w_0$ are bounded by (W1) and (M1). On the linear parts of the
profiles ($s\ge s_m$ for $\Theta^+$, $s\ge2$ for $\Theta^-$), $\mathcal V=\tilde v+\ve w_0$.

\emph{Step 4. The band nodes.} By (M2), a band node has a
uniform stencil of size $h_f$, whose nodes moved with velocity $v_n$, so that
$\mathcal D^{n+1}=\Delta_h+v_nD$ there. Along its path let
$\dot{\tilde v}=\tilde v_t+v_n\tilde v_x$ and $\ddot{\tilde v}$ be the derivatives of
$t\mapsto\tilde v(x_j(t),t)$; then $|\ddot{\tilde v}|\le C$. Since
$\tilde v(\tilde R(t),t)=0$ and $\tilde v_r<0$ near the front,
$\dot{\tilde v}=(\dot{\tilde R}(t)-v_n)|\tilde v_r|$ at the front, and
$\dot{\tilde R}(t)-v_n=(\dot{\tilde R}-\dot R)(t)+(\dot R(t)-\dot R(t_n))+(\dot R(t_n)-v_n)$. By
\eqref{cv:eq:speed}, for $a$ large, at distance $d$ from the front
\begin{equation}\label{cv:eq:rel}
|\dot{\tilde v}|\le\mu':=C_{\beta,a}(\mu+\mathcal E)+Cd,\qquad
(\dot{\tilde v})^\mp\le C(\mu_*+d),
\end{equation}
with the sign $-$ for the supersolution and $+$ for the subsolution: if $\Delta t\ge\ve$ then
$\mathcal E\ge\mu+\Delta t$, and the front of the supersolution moves outwards with respect to the nodes,
and that of the subsolution inwards. This is the favourable direction for the lagged
reaction rate, as we shall see. Let
\[
\zeta=\frac{2K_0h_f+\Delta t\sup\mu'}{\ve}+C(h_f+\Delta t)\log\frac1\ve
\]
(the supremum over $d\le C\ve\log(1/\ve)$). The first term bounds the variation of $s$ on the
evaluation points of a band node at distance $d\le C\ve\log(1/\ve)$ from the front, and the
second one the variation there of $s_m(w_0)$, $s_c(w_0)$, $s_*(w_0)$, whose derivatives with
respect to $W$ are bounded by $C\log(1/\gamma)$, and that of $\Theta^\pm(s;w_0)$ through $w_0$, relative
to $\Theta^\pm_s$. By
\eqref{cv:eq:hyp}, which implies $\Delta t\log(1/\ve)\le c_0$, $\zeta$ is as small as we wish if
$c_0$ is small (depending on $a$); in particular $\zeta\le\min\{\frac14,c_e\}/2$ and
$e^{C\zeta}\le2$ for the constants $C$ of (iv) and of Lemma~\ref{q:lem:inner}(ii), and then, by
(iv) and Lemma~\ref{q:lem:inner}(ii), $\Theta_s$ varies by a factor at most $2$ on these
evaluation points; moreover
$\Delta t\,\mu'^2\le\ve$ there (note that $\mu\le v_{\max}+\sup|\dot R|$ by (M1), so $\mu'$ is
bounded). For $N\ge2$ the terms with $\epsilon_j,\epsilon_j'$ and $\frac{N-1}{x_j}D$ in (ii) are bounded by
$C\sup\Theta_s$, and they only change the constants below; we write the proof for $N=1$.
The decomposition in (ii) is used only at band nodes at distance at most
$C\ve\log(1/\ve)$ from the front or on the burnt side of it, where $x_j\ge R_{\min}/2$; the
other band nodes on the fresh side, which may be close to $x=0$, are treated with (iii).

\emph{Step 5. The supersolution.} Here $\Theta=\Theta^+$, $\tilde v=\tilde v^+$, $\mathcal V=\ve\Theta^+(s;w_0)$. Let
$S_1=\lambda_0^{-1}\log(C_S/\ve)$. Where $s<s_m$, the equation of $\Theta$ and \eqref{q:eq:marginsup} give
$\frac1\ve\Theta_{ss}\tilde v_x^2\le\frac1\ve b(\Theta;w_0)-K\Theta_s-\Theta_{ss}\big(\frac{c\,\delta}\ve+\Lambda_1|s|\big)$, and
$\Theta_{ss}\ge0$; for $\Lambda$ large and $a$ large the last term absorbs the first term of
\eqref{cv:eq:Px}, as in region (A) of the proof of Proposition~\ref{q:prop:super}. Hence, a.e.\
where $s\ge-3S_1$,
\begin{equation}\label{cv:eq:Vxx}
\mathcal V_{xx}\le\Big[\frac1\ve\,b(\Theta;w_0)-K\Theta_s\Big]\mathbf 1_{\{s<s_m\}}+K_2\Theta_s ,
\end{equation}
with $K_2=K_1+C$, and the bracket, which equals $\ve^{-1}G_-^2\Theta_{ss}$, is nonnegative. Let
$\omega$ be the weight, in the
formula for $\delta^2$ in (ii), of the part of the stencil where $s<s_m(w_0)$. As $\tilde v$ is
monotone near the front, $\omega\ge\frac12$ if $s'<s_m$ and $\omega\le\frac12$ if $s'\ge s_m$.

\smallskip
(a) \emph{Band nodes with $-2S_1\le s'<s_m$.} By (i), (v) and (iv),
\[
\frac{\mathcal V'-\mathcal V_j}{\Delta t}\ge\Theta_s'\dot{\tilde v}'-C\Big(\frac{\Delta t\,\mu'^2}\ve+\big(\ve+\Delta t(1+\mu')\big)\log\frac1\ve\Big)\Theta_s'
\ge\Theta_s'\dot{\tilde v}'-C\,\Theta_s' .
\]
By (ii), $|v_n(\mathcal V_x'-D\mathcal V')|\le C\frac{h_f}\ve\Theta_s'$; and by \eqref{cv:eq:Vxx}, since
$x\mapsto b(\Theta;w_0)$ is Lipschitz with constant $C\ve^{-1}\sup\Theta_s$ on the stencil (we use
$|\partial_Wb|\le C\Theta\le C\Theta_s$), $b\ge0$ and $\omega\ge\frac12$,
\[
\delta^2\mathcal V'\le\frac1\ve\,b(\Theta';w_0')+C\frac{h_f}{\ve^2}\Theta_s'-\frac K4\Theta_s'+2K_2\Theta_s' .
\]
For the reaction term in \eqref{cv:eq:split}: $f$ is nonincreasing and Lipschitz and
$W_{h,j}/\ve\ge w_0'-\varpi'$, so
$f(\Theta_j+\ve\tilde q-W_{h,j}/\ve)\ge f(\Theta'-w_0')-L_f(\Theta_j-\Theta')^+-L_f(\ve\tilde q+\varpi')$, and
\[
f_\ve(V_j)(\mathcal V'+\ve^2\tilde q)
\ge\frac1\ve\,b(\Theta';w_0')-\frac{L_f}\ve(\Theta_j-\Theta')^+\Theta'-L_f\tilde q\,\Theta'+\tilde q\,f\Big(\Theta_j+\ve\tilde q-\frac{W_{h,j}}\ve\Big)
-L_f\frac{\varpi'}\ve\,\Theta' .
\]
The second term on the right is the lag. Its main part vanishes when $s_j\le s'$, that
is, when the
node moves towards the fresh side with respect to the barrier. In general, by
\eqref{cv:eq:rel}, $d\le C\ve(1+|s|)$, (iv), and the bound for $\Theta_W$ in
Lemma~\ref{q:lem:inner}(ii) for the variation of $w_0$ along the path,
\[
(\Theta_j-\Theta')^+\Theta'\le\Big(2(s_j-s')^++C(1+|s'|)\Delta t\Big)\Theta_s'\,\Theta'
\le C\Big(\frac{\mu_*\Delta t}\ve+\Delta t\Big)\Theta_s' .
\]
The last term is at most $C(\varpi'/\ve)\Theta_s'$, because $\Theta\le C\Theta_s$. Adding up, and using
$|\tilde v_t|\le K_1$ and Lemma~\ref{cv:lem:heat},
\[
\Lop[V]^n_j\ge\Big[\frac K4-C_2\Big(1+\frac{h_f}{\ve^2}+\frac{\Delta t+h_c}\ve+\frac{\mu_*\Delta t}{\ve^2}\Big)\Big]\Theta_s'
-L_f\tilde q\,\Theta'+\tilde q\,f\Big(\Theta_j+\ve\tilde q-\frac{W_{h,j}}\ve\Big).
\]
The last two terms are treated as in regions (A) and (B) of the proof of
Proposition~\ref{q:prop:super}: where $\Theta'\le\Theta_0/2$ their sum is at least $a_0\tilde q$
(because then $\Theta_j\le2\Theta'\le\Theta_0$, by $\Theta_s\le C\Theta$ and the choice of $\zeta$, so that
$\Theta_j+\ve\tilde q-W_{h,j}/\ve\le1-\delta_1/4$ for $c_0$ small and the value of $f$ is at least $2a_0$), and
where $\Theta'>\Theta_0/2$ we have $L_f\tilde q\,\Theta'\le C\Theta_s'$. Hence $\Lop[V]^n_j\ge0$ if $C_K$ is large.

\smallskip
(b) \emph{Band nodes with $s'\ge s_m$ and $\omega>0$.} Then
$s'\le s_m+2K_0h_f/\ve+Ch_f$, all the
evaluation points have $|s-s_m|\le2\zeta$, and there $\frac12\le\Theta_s\le1$. As in (a), now with
$\omega\le\frac12$ and dropping the term with $K$,
\begin{align*}
\Lop[V]^n_j&\ge-C\Big(1+\frac{\Delta t\,\mu'^2}\ve\Big)-\frac1{2\ve}\Big(b(\Theta';w_0')+C\frac{h_f}\ve\Big)
+\frac1\ve\Big(b(\Theta';w_0')-C\frac{\mu_*\Delta t}\ve-C\Delta t-C\ve-C\varpi'\Big)\\
&\ge\frac1{2\ve}\Big(b(\Theta';w_0')-C\,\frac{K\ve}{C_K}\Big).
\end{align*}
By Lemma~\ref{q:lem:inner}(a), $b(\Theta';w_0')\ge b(\Theta(s_m);w_0')-Ch_f/\ve\ge\frac\gamma4-Ch_f/\ve$, and
$\gamma\ge K\ve/c$, so $\Lop[V]^n_j\ge0$ if $C_K$ is large: at these nodes the slack is the
reaction term itself, which is of order $\gamma/\ve$ because $s_m<1$.

\smallskip
(c) \emph{Nodes with $\omega=0$ and $s'\ge s_m$, and nodes outside the band on the fresh side.}
If $C_b$ is large, the evaluation points of a node with $x'<R(t_{n+1})-\ell_b$ lie in
$\{\tilde u\ge2\ve\}$, because $|\tilde R-R|\le C_\beta a\mathcal E$ and the nodes move by at most
$v_{\max}\Delta t$. For a band node of this case, if its path $(x_j(t),t)$ met $\{s<s_m\}$
at some time of the step, then by \eqref{cv:eq:rel} the node would stay within a distance
$C\ve$ of the front during the whole step, and $|s-s'|\le\zeta<s_m$ on the path; hence in
all cases the evaluation points lie in $\overline{\{\tilde u>0\}}$, $s_j>0$, and $s_j<s_m$ is
possible only for a band node. In both cases $V=\tilde u+(\ve w_0-W_h)+\ve^2\tilde q$ on the
stencil, and
$V_j=\tilde u_j+(\ve w_{0,j}-W_{h,j})+\ve^2\tilde q+\ve g_j$ with $g=\Theta-s-w_0\ge0$, $g=0$ for $s\ge s_m$ and
$g\le C(s_m-s)^2$ for $s<s_m$. Since $\Lop_0[W_h]=0$,
$\Lop_0[V]^n_j=\Lop_0[\tilde u]^n_j+\ve\Lop_0[w_0]^n_j-\ve g_j/\Delta t$. If $s_j\ge s_m$, then $g_j=0$ and
$\Lop[V]^n_j\ge\Lop_0[\tilde u]^n_j+\ve\Lop_0[w_0]^n_j\ge0$ by
(iii), because the reaction term is nonnegative. If $s_j<s_m$ (a band node), then
$s_m(w_{0,j})-s_j\le(s'-s_j)+C\Delta t\le\Delta t(\mu'/\ve+C)$, because $s'\ge s_m(w_0')$ and
$|s_m'(W)|\le C$, so that $\ve g_j/\Delta t\le C(\Delta t\,\mu'^2/\ve+\ve\Delta t)\le C$, while $\Theta_j<\Theta(s_m;w_{0,j})$ and
$\Theta'\ge\Theta(s_m;w_0')$ give, with $f$ nonincreasing and Lipschitz and $|s_m'(W)|\le C$,
$f_\ve(V_j)(V'+W_h')\ge\ve^{-1}\big(b(\Theta(s_m);w_0')-C(\ve+\Delta t+\varpi')\big)\ge\gamma/(8\ve)$ for $a$ large;
hence
$\Lop[V]^n_j\ge\Lop_0[\tilde u]^n_j+\ve\Lop_0[w_0]^n_j-C+\gamma/(8\ve)\ge0$.

\smallskip
(d) \emph{The tail.} The remaining nodes are the band nodes with $s'<-2S_1$ and the
nodes with $x'>R(t_{n+1})+\ell_b$. All their evaluation points have $s\le-S_1+1$. For the
former: where $d\le r_1$, by \eqref{cv:eq:rel}, $d\le C\ve|s|$ and Gronwall's inequality
along the path, $s$ changes by at most
$\zeta+C\Delta t|s'|\le\zeta+|s'|/2$ on the evaluation points, and where $d>r_1$ we have
$s\le-\delta_0r_1/(2\ve)$. For the latter by Lemma~\ref{q:lem:ext}(i), if $C_b$ is large. There
$\Theta$, $\Theta_s$ and $\Theta_{ss}$ are at most $Ce^{\lambda_0s}\le Ce^{-\lambda_0S_1}$, and the derivatives of $\Theta$ with
respect to $W$ are at most $C(1+|s|)^2e^{\lambda_0s}$ by Lemma~\ref{q:lem:inner}(ii), so that
$|\mathcal P_x|+|\mathcal P_t|\le CS_1^2e^{-\lambda_0S_1}$ in \eqref{cv:eq:Vder}, without the factor $\ve^{-1}$ of
the terms with $\Theta_{ss}$; so by the mean
value theorem, (ii), (v) and
$|\dot x_j|\le v_{\max}$ the terms of $\Lop_0[\mathcal V]^n_j$ are bounded by
$C\ve^{-1}e^{-\lambda_0S_1}+CS_1^2e^{-\lambda_0S_1}\le a_0$ if $C_S$ is large, while, by \eqref{cv:eq:split},
$f_\ve(V_j)(\mathcal V'+\ve^2\tilde q)\ge f(\Theta_j+\ve\tilde q-W_{h,j}/\ve)\,\tilde q\ge2a_0$, because
$-W_{h,j}/\ve\le1-\delta_1+\varpi$. Within a distance
$\varrho_L/4$ of $x=L$ the function $\tilde v$ is constant, by the proof of Lemma~\ref{q:lem:ext}
and $r\le\varrho_L/4$, and $\partial_rw_0=0$ at $r=L$, so the zero-flux row $j=J$ is bounded in the
same way.

\smallskip
Since $V+W_h=\ve(\Theta^++\ve\tilde q)>0$, $V$ is a discrete supersolution.

\emph{Step 6. The subsolution.} Here $\Theta=\Theta^-$, $\tilde v=\tilde v^-$, $\mathcal V_1=\ve\Theta^-(s;w_0)$ and
$Z=\max\{\mathcal V_1,-\ve^2\tilde q\}$. At a node with $s'<s_*(w_0')$,
$Z'=-\ve^2\tilde q$ is the minimum of $Z$, so $Z'\le Z_j$ and
$(\mathcal D^{n+1}Z^{n+1})_j\ge0$ by \eqref{eq:peclet}, and
$\Lop[\underline V]^n_j\le-f_\ve(\underline V_j)\,\ve^2\tilde q\le0$. At a node with $s'\ge s_*(w_0')$ we have
$Z'=\mathcal V_1'$, and $Z\ge\mathcal V_1$ elsewhere, so that, again by \eqref{eq:peclet},
\begin{equation}\label{cv:eq:sub1}
\Lop[\underline V]^n_j\le\frac{\mathcal V_1'-\mathcal V_{1,j}}{\Delta t}-(\mathcal D^{n+1}\mathcal V_1^{n+1})_j
+f_\ve(\underline V_j)\,\mathcal V_1'.
\end{equation}

(a) \emph{Nodes with $s'\ge2$} (in the band, or outside it on the fresh side). All their
evaluation points have $s\ge\frac32$: near the front by the choice of $\zeta$, and at a distance
$d\ge C\ve\log(1/\ve)$ from it because there $s'\ge\delta_0\min\{d,r\}/\ve$ and, as in
Step~5(d), by Gronwall's inequality $s$ varies by at most $\zeta+C\Delta t\,s'\le\zeta+s'/2$ on
the evaluation points; for the nodes outside the band, whose stencils are not uniform,
the evaluation points lie in $\{\tilde u\ge2\ve\}$ if $C_b$ is large, as in Step~5(c). There $\mathcal V_1=\tilde u+\ve w_0+\ve g$ with $g=\Theta-s-w_0\ge0$, and $g$
vanishes at $(x',t_{n+1})$. Hence $\Lop_0[\ve g]^n_j\le0$ by \eqref{eq:peclet}, so
$\Lop_0[\mathcal V_1]^n_j\le\Lop_0[\tilde u]^n_j+\ve\Lop_0[w_0]^n_j\le0$ by (iii), and the reaction term
vanishes because $(\underline V_j+W_{h,j})/\ve\ge\Theta_j\ge s_j+w_{0,j}$, so that
$\underline V_j/\ve\ge s_j-\varpi\ge1$. Note that this uses that the profile is
linear only for $s\ge2$: a node that leaves the set $\{s<1\}$ during the step still carries
a positive lagged reaction rate, and it is then inside the damped part of the profile,
which is case (b).

(b) \emph{Band nodes with $s_*(w_0')\le s'<2$.} All the evaluation points have
$s\ge s_*(w_0)-c_e$, by the choice of $\zeta$. By the equation of $\Theta$,
$\frac1\ve\Theta_{ss}\tilde v_x^2=\frac{\tilde v_x^2}{G_+^2}\big(\frac1\ve b(\Theta;w_0)+K\Theta_s\big)$, and by
\eqref{q:eq:marginsub}, $1\le\tilde v_x^2/G_+^2\le1+C(\delta+\ve|s|)$. Where $\Theta\ge0$ we have $\Theta_{ss}\ge0$,
and the lower bound in \eqref{q:eq:marginsub} absorbs the first term of \eqref{cv:eq:Px}, as
in Step~5. Where $\Theta<0$ we use $|\Theta_{ss}|\le C\Theta_s$ in \eqref{cv:eq:Px}, so that
$|\mathcal P_x|\le C\log\frac1\ve\,\Theta_s$, which is at most $\frac K8\Theta_s$ because $K\ge C_K\log\frac1\ve$. Hence, a.e.,
\[
\mathcal V_{1,xx}\ge\Big[\frac1\ve\,b(\Theta;w_0)+\frac K2\Theta_s\Big]\mathbf 1_{\{s<2\}}-K_2\Theta_s
-C\Big(\delta+\ve\log\frac1\ve\Big)\frac{b(\Theta;w_0)^-}\ve ,
\]
and $b(\Theta;w_0)^-\le C(\ve\tilde q+\frac{h_f}\ve\Theta_s')$ on the stencil. If the stencil meets $\{s\ge2\}$ then
$b(\Theta;w_0)=0$ on it. Arguing as in Step~5(a), with the reverse inequalities (the
estimates of the difference quotient in time and of the convection term only use
$|\Theta_{ss}|\le C\Theta_s$ and (v)), and using \eqref{cv:eq:sub1},
\[
\Lop[\underline V]^n_j\le-\Big[\frac K8-C_2\Big(1+\frac{h_f}{\ve^2}+\frac{\Delta t+h_c}\ve+\frac{\mu_*\Delta t}{\ve^2}\Big)\Big]\Theta_s'
+C\Big(\delta+\ve\log\frac1\ve+\varpi'\Big)\mathbf 1_{\{\Theta<0\text{ on the stencil}\}} .
\]
(The constants depend on $C_1$, hence on $\beta$ and $T$.) For the reaction term we used
the following. Let $\underline\Theta_j=\max\{\Theta_j,-\ve\tilde q\}$, so that
$f_\ve(\underline V_j)=\ve^{-2}f(\underline\Theta_j-W_{h,j}/\ve)$, and recall that
$|W_{h,j}/\ve-w_0'|\le\varpi'$. Where $\Theta'\ge0$,
$f(\underline\Theta_j-W_{h,j}/\ve)\Theta'\le b(\Theta';w_0')+L_f\big((\Theta'-\Theta_j)^++\varpi'\big)\Theta'$; the lag is
now controlled by
$(\dot{\tilde v})^+$ in \eqref{cv:eq:rel} and by the variation of $w_0$ along the path, as in
Step~5(a), and $\varpi'\Theta'\le C\varpi'\Theta_s'$ because $\Theta^-_s\ge c\,\Theta^-$ where $\Theta^-\ge0$. Where
$\Theta'<0$, $|f(\underline\Theta_j-W_{h,j}/\ve)-f(\Theta'-w_0')|\le L_f(|\underline\Theta_j-\Theta'|+\varpi')$ and
$|\Theta'|\le\ve\tilde q$, so that
$f(\underline\Theta_j-W_{h,j}/\ve)\Theta'\le b(\Theta';w_0')+L_f\big(|\Theta_j-\Theta'|+\varpi'\big)\ve\tilde q$, and
$|\Theta_j-\Theta'|\le C\zeta\,\Theta_s'$ goes into the term with $C_2$. Where the stencil meets
$\{\Theta<0\}$ we have $\frac K{16}\Theta_s'\ge c\sqrt\gamma\ge C(\delta+\ve\log\frac1\ve+\varpi')$ for $c_0$ small
(depending on $a$), because
$\delta+\ve\log\frac1\ve+\varpi'\le C(a+1)\mathcal E$ and $\gamma=c_1c\,a\mathcal E$. Hence $\Lop[\underline V]^n_j\le0$ if $C_K$ is large.

If $C_b$ is large, every node with $s'\ge s_*(w_0')$ is of one of these two kinds, because
$s_*\ge-C\log(1/\ve)$. Since no sign condition is required for subsolutions, $\underline V$ is a
discrete subsolution.

\emph{Step 7. Conclusion.} At $t=0$, $W_h^0=\ve w_0(\cdot,0)$ on the mesh, and
$|1-\lambda(0)|=a\mathcal E\ge a\,\ve\log\frac1\ve$. The verification of the initial ordering in the proof of
Theorem~\ref{q:thm:rate}, which only uses (H), \eqref{q:eq:data}, Lemma~\ref{q:lem:ext}(i) and
the bounds $s_c\ge-C\log(1/\gamma)$ and $\Theta^-\le0$ for $s\le s_c$ of Lemma~\ref{q:lem:inner}(b),
gives $\underline V(\cdot,0)\le u_0^\ve\le V(\cdot,0)$ in $[0,L]$ if $a$ is large. By
Lemma~\ref{lem:comparison}, $\underline V^n_j\le U^n_j\le V^n_j$ for $t_n\le T$, and
$|V-\tilde u|+|\underline V-\tilde u|\le C\ve+|\ve w_0-W_h|\le C\ve$, $|\tilde u-v|\le C_\beta a\mathcal E$. \qed

\subsection{Travelling waves: no time step in the error}\label{an:sec:ale}

For a travelling wave, on a mesh that translates with the wave, the solution is steady
in the frame of the mesh, and the term $\Delta t$ of Theorem~\ref{cv:thm:main} disappears. Here
we compare with the travelling wave $\Phi=\Phi_\ve$ of \eqref{eq:Peps} itself
(Lemma~\ref{an:lem:profile}), with speed $c=c_\ve$, so that the estimate has no term
$\ve\log(1/\ve)$ either; the price is a worse dependence on the mesh size. We take $N=1$,
$w^\ve\equiv\ve w_0$, assume only (F), and consider (S) for bounded sequences on the
infinite grid $x^n_j=jh+X_n$, with $X_{n+1}=X_n+v_n\Delta t$, $|v_n|\le v_{\max}$,
$v_{\max}h\le2$ and $h\le\ve\le\ve_0$. On the infinite grid the results of
Section~\ref{sec:theory} hold in $\ell^\infty$: for a bounded $F\ge0$ the matrix
$\mathsf A_F=(\Delta t^{-1}+F)I-\delta^2-v_nD$ of the proof of Lemma~\ref{lem:bounds} is
$\mathsf A_F=\mathfrak d\,(I-\mathfrak d^{-1}\mathfrak n)$, where
$\mathfrak d=\operatorname{diag}(\Delta t^{-1}+F_j+2h^{-2})$ and $\mathfrak n\ge0$ (because
$v_{\max}h\le2$) has row sums $2h^{-2}$; since
$\|\mathfrak d^{-1}\mathfrak n\|_{\ell^\infty\to\ell^\infty}\le2h^{-2}/(\Delta t^{-1}+2h^{-2})<1$, $\mathsf A_F$ is
invertible in $\ell^\infty$ and $\mathsf A_F^{-1}=\sum_{k\ge0}(\mathfrak d^{-1}\mathfrak n)^k\mathfrak d^{-1}\ge0$.
Hence the scheme has a unique bounded solution, and
Lemmas~\ref{lem:bounds} and \ref{lem:comparison}, whose proofs only use $\mathsf A_F^{-1}\ge0$ and
$\mathsf A_0\mathbf 1=\Delta t^{-1}\mathbf 1$, hold for bounded sequences. With
$\delta^2V_j=(V_{j+1}-2V_j+V_{j-1})/h^2$ and $DV_j=(V_{j+1}-V_{j-1})/(2h)$ the scheme reads
$\Lop_v[U]=0$, where
\begin{equation}\label{an:eq:Lv}
\Lop_v[V]^n_j=\frac{V^{n+1}_j-V^n_j}{\Delta t}-\delta^2V^{n+1}_j-v_n\,DV^{n+1}_j
+f_\ve(V^n_j)\big(V^{n+1}_j+\ve w_0\big).
\end{equation}
Let $\rho_h=h/\ve^2$, or $\rho_h=h^2/\ve^3$ if $f\in C^{1,1}(\R)$ (which forces $f'(1)=0$ and
excludes $f(s)=(1-s)^+$), and let
$\mathcal M_n=\Delta t\sum_{k<n}|v_k-c_\ve|$ be the accumulated mesh-velocity error.

\begin{theorem}[Travelling waves on moving meshes]\label{an:thm:ale}
Assume (F) and let $U$ be the solution of (S) on the grid $x^n_j=jh+X_n$ with
$U^0_j=\Phi_\ve(x^0_j)$. There are $C_A,\nu>0$, depending only on $f,w_0,A$ and $v_{\max}$, such
that, with $\sigma_n=|v_n-c_\ve|+C_A\rho_h$ and
$\xi_n=\Delta t\sum_{k<n}\sigma_k=\mathcal M_n+C_A\rho_ht_n$, if
\begin{equation}\label{an:eq:dtale}
\sigma_n\,\Delta t\le\nu\,\ve\qquad\text{for all }n<\bar n,
\end{equation}
then for all $n\le\bar n$ and $j$
\[
\Phi_\ve\big(x^n_j-c_\ve t_n-\xi_n\big)\ \le\ U^n_j\ \le\ \Phi_\ve\big(x^n_j-c_\ve t_n+\xi_n\big).
\]
In particular $|U^n_j-\Phi_\ve(x^n_j-c_\ve t_n)|\le K(\mathcal M_n+C_A\rho_ht_n)$.
\end{theorem}

As in Theorem~\ref{cv:thm:main}, the velocities $v_n$ are given; if they are computed from
the solution, the statement is a posteriori, since $\mathcal M_n$ depends on the realised
velocities. There is no circularity: the comparison principle is applied to the realised
sequence $(v_n)$. The barriers are translates of the wave, as in the stability theory of
travelling fronts of Fife and McLeod \cite{FifeMcLeod1977}; what is specific to the scheme
is the sign of the lag of the reaction rate.

\begin{proof}
Let $V^n_j=\Phi(z^n_j)$ with $z^n_j=x^n_j-ct_n+\xi_n$. Then
$z^{n+1}_j-z^n_j=d_n\Delta t$ with $d_n=v_n-c+\sigma_n\in[0,2\sigma_n]$: the barrier moves
towards the burnt side with respect to the nodes. Write $z'=z^{n+1}_j$ and
$z=z^n_j\in[z'-2\nu\ve,z']$. Now $\Phi(z)\le\Phi(z')$, so
$f_\ve(\Phi(z))\ge f_\ve(\Phi(z'))$ and the reaction term is \emph{at least}
$g_\ve(\Phi(z'))$: the lag of the reaction rate has the favourable sign. By Taylor's
formula and Lemma~\ref{an:lem:profile}(3) (all the points involved are within $\ve$ of
$z'$, so $\Phi'$ there is at most $e^K\Phi'(z')$),
\begin{gather*}
\frac{\Phi(z')-\Phi(z)}{\Delta t}\ge d_n\Phi'(z')-2\sigma_n^2\Delta t\,\tfrac K\ve e^K\Phi'(z'),\\
|D\Phi(z')-\Phi'(z')|+|\delta^2\Phi(z')-\Phi''(z')|\le2Ke^K\rho_h\Phi'(z') ;
\end{gather*}
for the second differences we used
$\delta^2\Phi(z')-\Phi''(z')=h^{-2}\int_{-h}^h(h-|y|)\big(\Phi''(z'+y)-\Phi''(z')\big)dy$, which is
bounded by $\frac h3\sup|\Phi'''|$ and, if $f\in C^{1,1}$, by $\frac{h^2}{12}\sup|\Phi''''|$. Since
$d_n\Phi'-\Phi''-v_n\Phi'+g_\ve(\Phi)=\sigma_n\Phi'$ by \eqref{an:eq:tw}, we get, using
\eqref{an:eq:dtale},
\[
\Lop_v[V]^n_j\ge\big(\sigma_n(1-2Ke^K\nu)-C\rho_h\big)\Phi'(z')\ge0
\]
if $2Ke^K\nu\le\frac12$ and $C_A\ge2C$. So $V$ is a supersolution, and $U\le V$. The
subsolution is $\Phi(x^n_j-ct_n-\xi_n)$: it moves towards the fresh side with respect to
the nodes, so $\Phi(z)\ge\Phi(z')$, the reaction term is at most $g_\ve(\Phi(z'))$, and the
same computation gives $\Lop_v[V]\le-(\sigma_n(1-2Ke^K\nu)-C\rho_h)\Phi'(z')\le0$.
\end{proof}

If the mesh moves exactly with the wave, $v_n\equiv c_\ve$, the error is at most
$KC_A\rho_ht_n$ and \eqref{an:eq:dtale} is $\Delta t\le\nu\ve/(C_A\rho_h)$; when $\rho_h\lesssim\ve$
this is $\Delta t\le\nu'$, independent of $\ve$. The barriers are translates of $\Phi_\ve$, whose
slack $\sigma_n\Phi'$ is spread over the whole profile; this is why the mesh size enters
through $\rho_h$, with one more power of $\ve^{-1}$ than in Theorem~\ref{cv:thm:main}, where
the slack is concentrated in the layer. For $v_n\equiv0$ the theorem gives nothing,
because $\mathcal M_n=c_\ve t_n$: a barrier that keeps the favourable direction with respect to
fixed nodes cannot follow the wave, and one is in the situation of
Remark~\ref{cv:rem:main}(1).

The absence of the time step has a simple reason. When $v_n\equiv v$ is constant, the scheme
on the translating grid is autonomous, and its steady states, the bounded solutions of
$-\delta^2U-v\,DU+f_\ve(U)(U+\ve w_0)=0$, do not depend on $\Delta t$: a linearly implicit
scheme leaves steady states invariant. A discrete travelling wave that moves with the
grid is such a steady state, and Theorem~\ref{an:thm:ale} shows that, uniformly in the time
step, the solution of the scheme stays close to the continuous wave, up to the mesh size
and the accumulated error of the mesh velocity. Accordingly, it
is a statement about travelling waves on rigidly translating grids, and not about the
method of Section~\ref{sec:ale} in general: there the velocity is computed, only a
neighbourhood of the band moves rigidly, the mesh is rebuilt, and for other data the
solution is not steady in any frame; in all these cases the time step enters the error
(Sections~\ref{an:sec:impl} and \ref{sec:experiments}). In the experiments the speed of the
discrete wave depends on the position $\vartheta h$ of the front inside its cell, and its
difference with $c_\ve$ is the lattice pinning error of Section~\ref{sec:exp:tw}.

\subsection{The scheme as implemented, and what is not sharp}\label{an:sec:impl}

The code of Section~\ref{sec:method} differs from the setting of
Theorem~\ref{cv:thm:main} in three respects.

\emph{Mesh velocity.} The velocity $v_n$ is computed, by the condition that the front
stays on the reference point of the band. The barriers only use the realised sequence
$(v_n)$, so the theorem applies with the realised $\mu$, but this is an a posteriori
statement: the barriers confine the front to a distance $C\mathcal E$ of $R(t)$, which
bounds the mesh displacement over many steps but not the velocity of a single step,
and an a priori bound for the implicit velocity is open. In the experiments $\mu$ is
small (Section~\ref{sec:exp:tw}).

\emph{Band and remeshing.} The experiments use $h_f=\ve/16$, a band of fixed width
$55\ve$ and $h_c=0.01$, so the band is narrower than $C_b\mathcal E$ when $h_c\gg\ve$; and
the mesh is rebuilt from time to time. The transfer $\Pi$ is order preserving and is the
identity on the band, but outside the band the transferred barrier is not the barrier
on the new nodes; the difference, an interpolation error of order $h_c^2$, is not
covered by the proof. Also, only a neighbourhood of the band moves rigidly, which is
allowed by (M1)--(M2); but if the rigid zone is short, the nodes ahead of it see a
solution that changes in time at a rate of order one, and an error $O(\Delta t)$ appears
even for travelling waves. Hence the rigid zone must contain several decay lengths $1/c$
of the profile. Finally, Theorem~\ref{an:thm:ale} is stated on the whole line, while the
code works on $[0,L]$ with zero flux at both ends; the translates of $\Phi_\ve$ do not
satisfy this condition, with a defect that is exponentially small in $L/\ve$ at the burnt
end and of the size of $\Phi_\ve'$ at the fresh end.

\emph{What is sharp.} The barriers absorb the truncation error pointwise. In the layer
this error, relative to the derivative of the profile, is of order
$\rho=\Delta t/\ve^2$ for the lagged reaction on a fixed mesh and $h_f/\ve^2$ for the second
differences, and a damping $K\sim\rho$ of the inner profile costs $K\ve$ in the free
boundary condition: this is the origin of the terms $\Delta t/\ve$ and $h_f/\ve$. A formal
expansion, which we do not justify, gives the same factor $\ve$. Let
$\mathcal A=\partial_z^2+c\,\partial_z-g_\ve'(\Phi)$ be the linearization at the travelling
wave, and $\tau$ the truncation error of the scheme on it; $\mathcal A\Phi'=0$,
$\psi^*=e^{cz}\Phi'$ spans the kernel of the adjoint, and the drift of the front is
\[
\lambda=-\frac{\int\tau\,\psi^*\,dz}{\int\Phi'\psi^*\,dz},
\]
an average of $\tau$ in which the layer has weight $\ve$. (This is only heuristic: $\psi^*$
does not decay as $z\to+\infty$ and the essential spectrum of $\mathcal A$ reaches the
origin, so that a proof would require exponentially weighted spaces. The absence of a
spectral gap is also behind the slow relaxation, on the time scale $4/c^2$, observed in
Section~\ref{sec:exp:tw}.) For the lagged reaction we have
evaluated it on the exact wave ($f(s)=(1-s)^+$, $A=1$, $w_0=0$):
\begin{equation}\label{an:eq:solv}
\lambda=-0.235\,\frac{\Delta t}\ve\,(1+o(1))\qquad(\ve=0.02,\ 0.01,\ 0.005:\ -0.2351,\ -0.2345,\ -0.2342),
\end{equation}
which is what the experiments show for travelling waves (Section~\ref{sec:exp:tw}); in
the setting of Theorem~\ref{cv:thm:main} the error on a fixed mesh is also proportional
to $\Delta t/\ve$ (Section~\ref{sec:exp:thm}). The terms $h_f/\ve$ and $h_c$ are not sharp: the second
differences of a smooth profile are of second order, their average over the position
of the front in its cell is smaller still, and the experiments show errors of order
$(h_f/\ve)^2$ and $h_c^2$.

\section{Numerical experiments}\label{sec:experiments}

In all the experiments $f(s)=(1-s)^+$, so that $M(w_0)=(1+w_0)^3/6$. Unless stated
otherwise the meshes are those of Section~\ref{sec:meshes} with $h_f=\ve/16$ and
$h_c=0.01$, and the initial data are \eqref{eq:initdata}. We use the three solvers
of Section~\ref{sec:method}, and say which one in each case: \emph{ALE} (the scheme
(S) on the moving mesh, with implicit mesh velocity and a fixed time step),
\emph{Eulerian (S)}, and the reference solver \emph{BDF}. The exact speed $c_\ve$ of the
travelling wave of \eqref{eq:Peps} is computed by shooting (Section~\ref{sec:tw}). With
$h_f=\ve/16$ the layer is resolved, as Theorem~\ref{cv:thm:main} requires; the experiments
differ from its setting in the width of the band, the remeshing, the computed mesh
velocity (Section~\ref{an:sec:impl}), and in that many of them use $\Delta t\gg\ve$.

\emph{The range of parameters.} Lengths are measured in units for which the flame speed
is of order one, so the thermal thickness of the flame is $1/c_0$ ($1.7$ for $w_0=0$ and
$A=1$). Most experiments use $\ve$ between $10^{-2}$ and $10^{-3}$, that is, Zeldovich numbers
between $100$ and $1000$, far above those of real flames, which are between $5$ and $20$:
our subject is the regime $\ve\to0$, in which the estimates are stated and in which methods
that resolve the layer on a fixed mesh become expensive. For $\ve\ge0.05$ the band, of
width $55\ve$, covers the whole flame, and the method reduces to a fine uniform mesh that
travels with it. The value $h_c=0.01$ is small with respect to the thermal thickness; it
was chosen so that the error of the coarse zone is negligible in the tables.

\subsection{The meshes: refinement ahead of the flame, coarsening behind it}

\begin{figure}[ht]
\centering
\includegraphics[width=\textwidth]{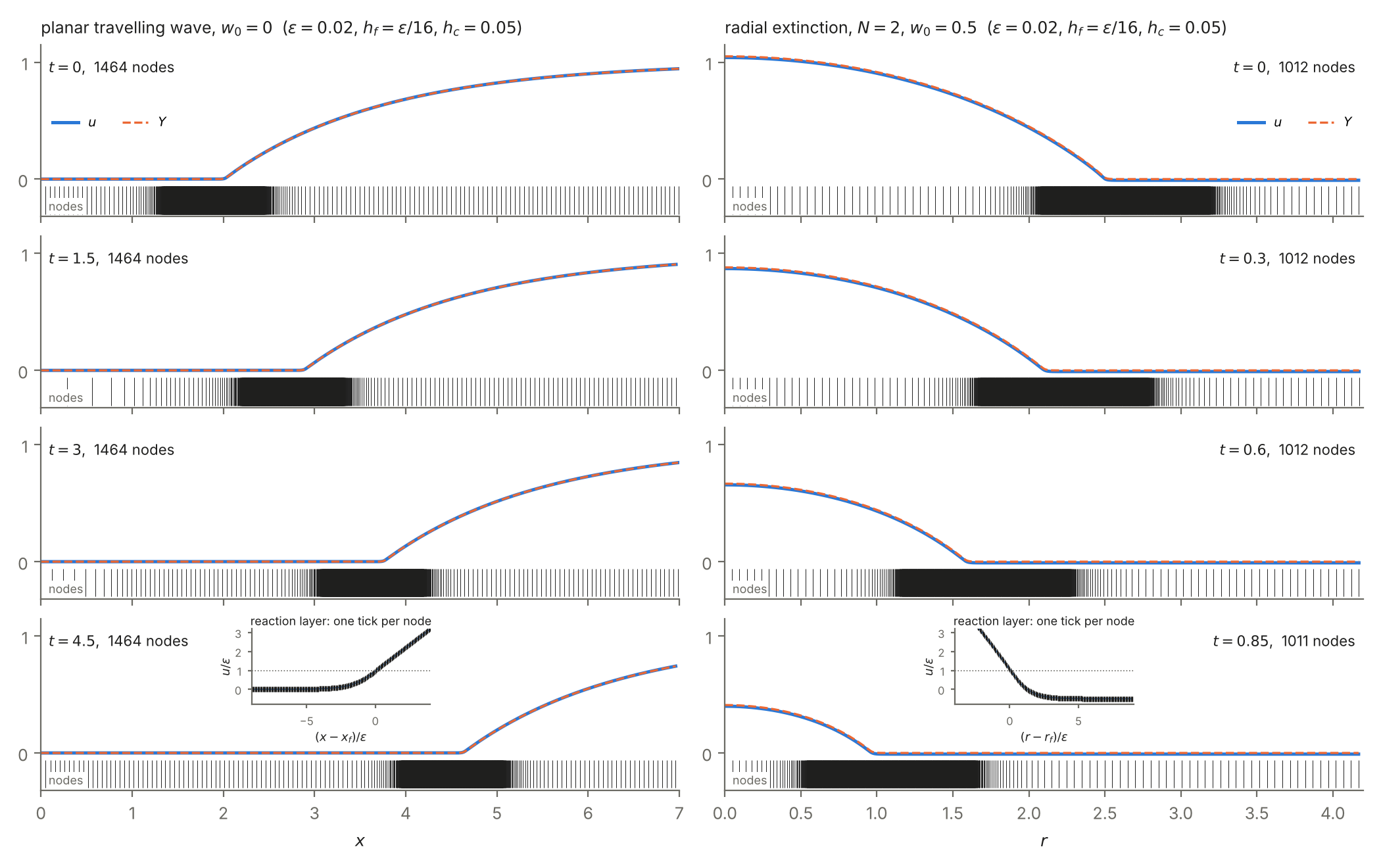}
\caption{ALE: the solution ($u$ solid, $Y$ dashed) and the mesh (one tick per node) at four
times, with $\ve=0.02$, $h_f=\ve/16$ and $h_c=0.05$ (to make the coarse nodes visible).
Left: planar travelling wave. Right: radial extinction of the self-similar solution,
$N=2$. Insets: the reaction layer, in units of $\ve$.}\label{fig:A5}
\end{figure}

Figure~\ref{fig:A5} shows the solution and the mesh at four times, for a planar
travelling wave and for the extinction of a ball in radial symmetry. The fine band
stays on the reaction layer; the cells are refined where the flame arrives and
coarsened where it has passed, and the number of nodes changes by at most one. The inset
shows the reaction layer, of width a few $\ve$, resolved by the band.

In the ALE version the band, the transitions and the coarse nodes near the front
translate with the flame, and the band nodes stay on their lattice at remeshing: for the
planar wave with $\ve=0.01$ and $0\le t\le4$ the mesh is rebuilt once. On meshes that do not
move it is rebuilt $31$ times, and at each rebuild the band jumps ahead of the front.

\emph{Number of nodes.} For a fixed ratio $h_f/\ve$ the number of nodes is constant in
time and grows only logarithmically in $1/\ve$: for the planar travelling wave ($w_0=0$,
$L=27$, $0\le t\le4$) the mesh has $3573$ nodes for $\ve=0.01$ and $3671$ for $\ve=0.001$, of
which $882$ are in the band; the difference is in the transitions. The cost of the method
is compared with that of other methods in Section~\ref{sec:exp:tw}.

\subsection{Comparison with the computed limit problem}\label{sec:limit}
To compare \eqref{eq:Peps} with \eqref{eq:P} for data that are not exact solutions we
solve \eqref{eq:P} directly, by front fixing. In radial symmetry with
$\{u>0\}=\{r<s(t)\}$ we set $r=s\xi$, so that $v(\xi,t)=u(s\xi,t)$ solves
\[
v_t=\frac{1}{s^2}\Big(v_{\xi\xi}+\frac{N-1}{\xi}v_\xi\Big)+\xi\frac{\dot s}{s}v_\xi,\quad 0<\xi<1,
\qquad v_\xi(0)=0,\quad v(1)=0,\quad v_\xi(1)=-s\sqrt{2M}.
\]
Both free boundary conditions are imposed, the second one through a ghost
node, and differentiating $v(1,t)=0$ gives $\dot s=\Delta u|_{r=s}/\sqrt{2M}$.
In the planar case we use $x=s+(L-s)g(\xi)$ with $g$ clustered at the front.
Second-order differences in $\xi$ and BDF in time give a second-order method: against
the travelling wave and the self-similar solution the errors in the free boundary are
at most $6\cdot10^{-6}$ with $800$ points.

\begin{figure}[ht]
\centering
\includegraphics[width=0.9\textwidth]{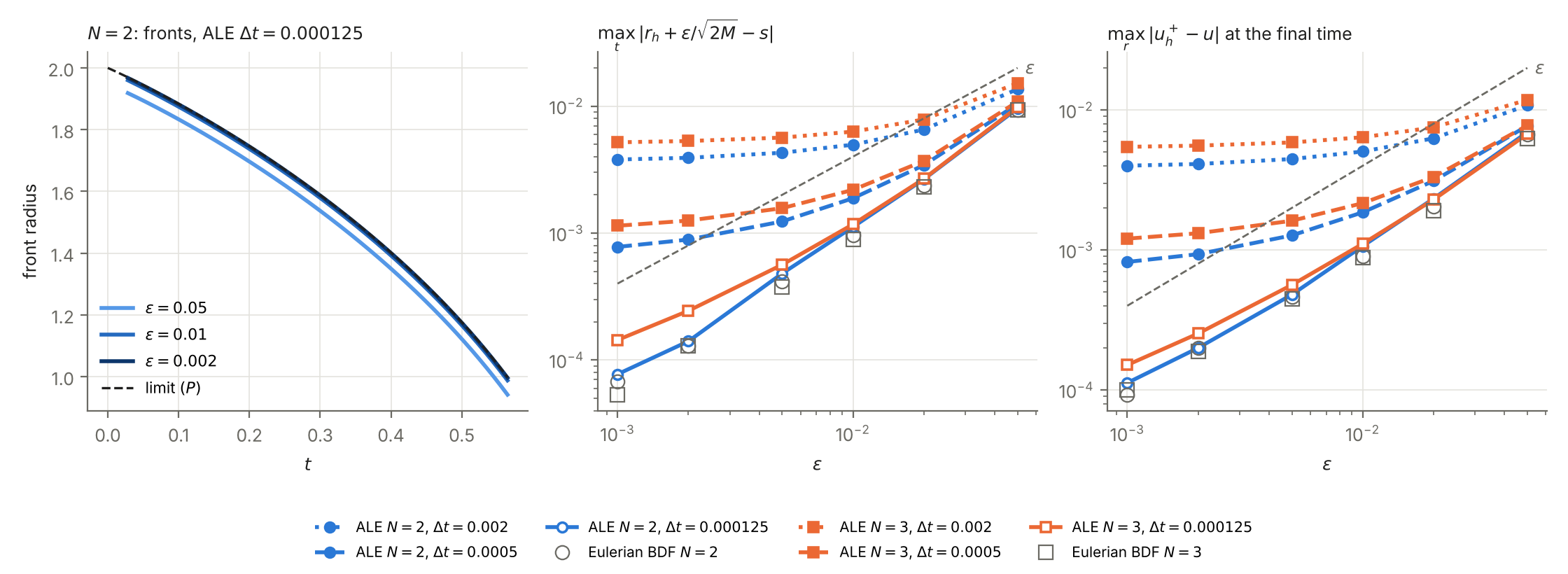}
\caption{\eqref{eq:Peps} against the computed limit problem for the radial data
$u_0=a(1-r^2/R_0^2)^+$ ($w_0=0.5$). Left: fronts ($N=2$). Middle: front$^*$. Right: sup
error at the final time. Colour: ALE with $\Delta t=0.002$, $0.0005$, $0.000125$; grey:
BDF.}\label{fig:C1}
\end{figure}

We take $w_0=0.5$ and data that satisfy the free boundary condition:
$u_0=A\tanh\big(\sqrt{2M}(x-2)/A\big)^+$ (planar, $A=1$, up to $t=6$) and
$u_0=a(1-r^2/R_0^2)^+$ with $a=\sqrt{2M}R_0/2$, $R_0=2$ (radial, $N=2,3$, until the radius
halves). Table~\ref{tab:limit} shows the errors of the BDF solver, and
Figure~\ref{fig:C1} those of ALE.

\begin{table}[ht]
\centering
\small
\begin{tabular}{r cc cc cc}
\toprule
& \multicolumn{2}{c}{planar} & \multicolumn{2}{c}{radial $N=2$} & \multicolumn{2}{c}{radial $N=3$}\\
$\ve$ & front$^*/\ve$ & sup$/\ve$ & front$^*/\ve$ & sup$/\ve$ & front$^*/\ve$ & sup$/\ve$\\
\midrule
0.05  & 0.391 & 0.369 & 0.187 & 0.134 & 0.186 & 0.125\\
0.02  & 0.325 & 0.327 & 0.118 & 0.102 & 0.115 & 0.096\\
0.01  & 0.302 & 0.311 & 0.096 & 0.090 & 0.089 & 0.089\\
0.005 & 0.300 & 0.312 & 0.083 & 0.092 & 0.075 & 0.089\\
0.002 & 0.321 & 0.333 & 0.065 & 0.100 & 0.065 & 0.095\\
0.001 & 0.296 & 0.302 & 0.067 & 0.093 & 0.053 & 0.101\\
\bottomrule
\end{tabular}
\caption{\eqref{eq:Peps} (BDF) against the computed solution of \eqref{eq:P}: error in the
position of the front after removing the shift $\ve/\sqrt{2M}$ (front$^*$), and sup error
of $(u^\ve)^+$ at the final time, both divided by $\ve$.}\label{tab:limit}
\end{table}

All the errors are $O(\ve)$. The radial data are pockets of fresh gas as in
Theorem~\ref{q:thm:rate}, and the errors agree with its rate; the planar front is not
covered by the theorem, whose setting requires a bounded fresh region, but it shows the
same order. The data have
$Y_0^\ve=0$ in the burnt region, and no logarithm is visible. Most of the raw error in the
front is the shift $\ve/\sqrt{2M}$ between the level set $\{u^\ve=\ve\}$ and the free
boundary. What remains is $0.3\,\ve$ in the planar case; in radial symmetry it decreases from
$0.19\,\ve$ to $0.05$--$0.07\,\ve$ as $\ve$ decreases, so it may be of lower order there. The planar value is explained by Proposition~\ref{prop:speed}: the
front moves essentially as a travelling wave with $A=1$, whose speed is
$c_0+0.044\,\ve+o(\ve)$ for $w_0=0.5$, and over $0\le t\le6$ this accumulates to
$\approx0.26\,\ve$. With ALE and $\Delta t=1.25\cdot10^{-4}$ the errors are comparable to those of
BDF ($7.7\cdot10^{-5}$ against $6.7\cdot10^{-5}$ for front$^*$ at $\ve=0.001$, $N=2$); for larger
time steps the time discretization dominates at small $\ve$ ($3.8\cdot10^{-3}$ and
$7.8\cdot10^{-4}$ for $\Delta t=2\cdot10^{-3}$ and $5\cdot10^{-4}$).

\subsection{A non-constant $w_0$}\label{sec:exp:varw}

\begin{figure}[ht]
\centering
\includegraphics[width=\textwidth]{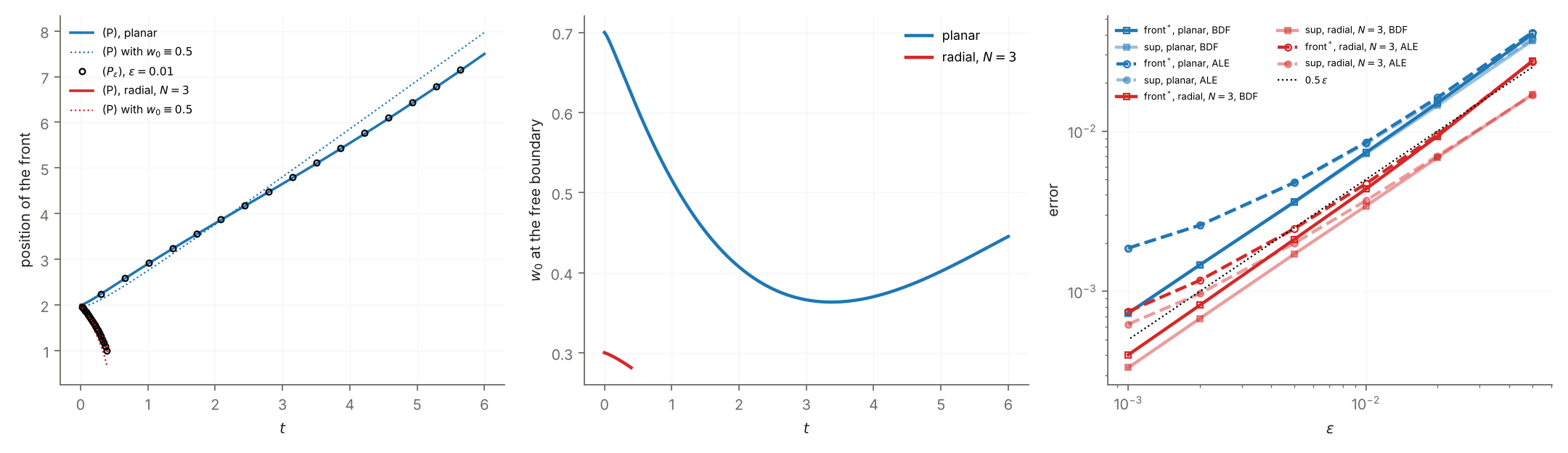}
\caption{Non-constant $w_0$. Left: fronts of the computed limit problem \eqref{eq:P}
(solid), of \eqref{eq:P} with $w_0\equiv0.5$ (dotted) and of \eqref{eq:Peps} with $\ve=0.01$
(circles). Middle: $w_0$ at the free boundary. Right: front$^*$ and sup error, as in
Table~\ref{tab:limit}, for BDF and for ALE ($\Delta t=2.5\cdot10^{-4}$ planar,
$1.25\cdot10^{-4}$ radial).}\label{fig:I1}
\end{figure}

In all the experiments above $w^\ve\equiv\ve w_0$ with $w_0$ constant. When $w^\ve(\cdot,0)=\ve
W_0$ with $W_0$ not constant, $w_0(x,t)$ is the solution of the heat equation with initial
datum $W_0$, and the free boundary condition of \eqref{eq:P} is
$|\nabla u|=\sqrt{2M(w_0(x,t))}$: the flame moves in a mixture that is stratified in
enthalpy. This is the feature that distinguishes \eqref{eq:P} from the problem of
\cite{CaffarelliVazquez1995}. The scheme needs no change, since $W_h=Y-U$ solves the
discrete heat equation, and in the front-fixing solver $M$ is evaluated at the current
front. We consider a planar front, with $L=24$,
\[
w_0(x,t)=0.5+0.4\,e^{-k^2t}\cos(kx),\quad k=\tfrac\pi6,\qquad
u_0=A\tanh\big(\sqrt{2M(w_0(2,0))}\,(x-2)/A\big)^+,\quad A=1,
\]
up to $t=6$, and a pocket of fresh gas in radial symmetry, with $N=3$, $L=6$,
\[
w_0(r,t)=0.5-0.3\,e^{-k^2t}\,\frac{\sin(kr)}{kr},\quad \tan(kL)=kL,\qquad
u_0=a\,(1-r^2/R_0^2)^+,\quad R_0=2,
\]
with $a=\sqrt{2M(w_0(R_0,0))}\,R_0/2$,
until the radius halves. In both cases $w_0$ has zero flux at the ends of the domain and
the data satisfy the free boundary condition at $t=0$; in the second one $w_0$
is monotone under parabolic dilations in a neighbourhood of the pocket, which is the
condition \eqref{q:eq:mono} of Theorem~\ref{q:thm:rate} (stated there in the whole space):
it satisfies (W1) and (W2) of
Section~\ref{cv:sec:main}, with $L_w=\pi/k\approx4.2$, so that it is of the kind considered in
Theorem~\ref{cv:thm:main}. Along the free
boundary $w_0$ goes from $0.70$ to $0.36$ and back to $0.45$ in the planar case, and
from $0.30$ to $0.28$ in the radial one (Figure~\ref{fig:I1}); at $t=6$ the planar front
is at $x=7.50$, while with $w_0$ replaced by its mean value $0.5$ it would be at $7.97$.

The errors of \eqref{eq:Peps}, computed with BDF, against the computed limit problem are
proportional to $\ve$, with no visible logarithm: for $\ve$ from $0.05$ to $0.001$, front$^*/\ve$
(with the shift $\ve/\sqrt{2M}$ evaluated at the current front) is $0.80$, $0.75$, $0.74$,
$0.72$, $0.73$, $0.73$ in the planar case and $0.55$, $0.47$, $0.44$, $0.42$, $0.41$, $0.40$
in the radial one, and the sup error divided by $\ve$ is between $0.72$ and $0.73$, and
between $0.33$ and $0.34$, respectively. With ALE the errors are the same up to the time
discretization, which is visible for the smallest values of $\ve$ (front$^*=1.9\,\ve$ and
$0.75\,\ve$ for $\ve=0.001$), as in Section~\ref{sec:limit}. The discrete $W_h/\ve$ differs
from the exact $w_0$ by at most $10^{-4}$ at the final time.

\subsection{Dependence of the limit on $w_0$}\label{sec:exp:w0}
\begin{figure}[ht]
\centering
\includegraphics[width=0.85\textwidth]{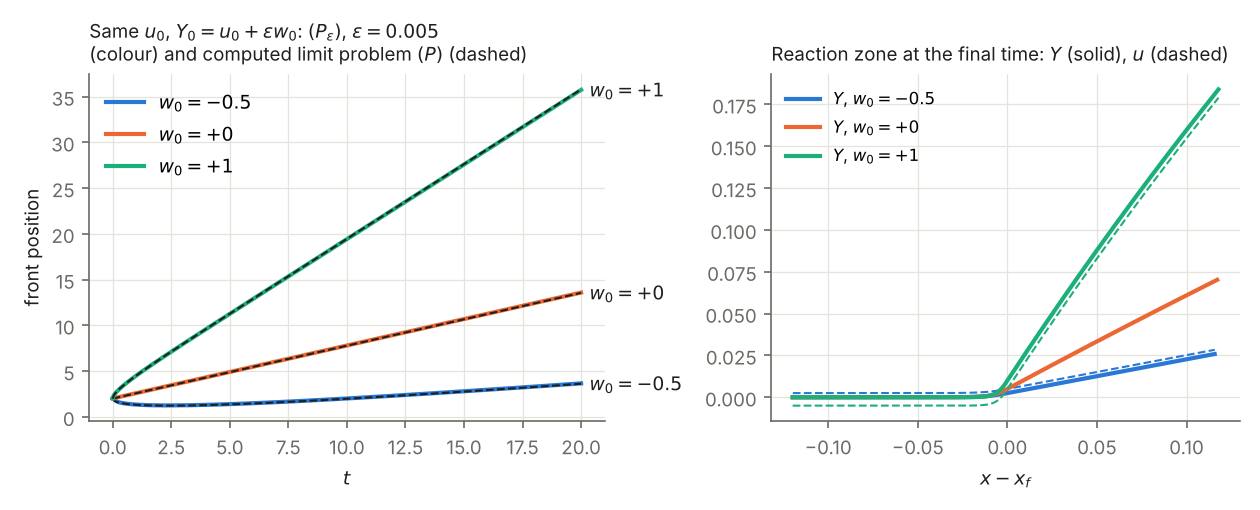}
\caption{Same $u_0$, three values of $w_0$, $\ve=0.005$ (BDF). Left: fronts of
\eqref{eq:Peps} (colour) and of the computed limit problem \eqref{eq:P} (dashed).
Right: $Y$ (solid) and $u$ (dashed) near the front at $t=20$.}\label{fig:instability}
\end{figure}

To display the dependence of the limit on the perturbation $w_0$, proved in
\cite{FernandezBonderWolanski2000}, we
take the \emph{same} initial temperature $u_0$ (the limit travelling wave for
$w_0=0$) and three values of $w_0$, that is, $Y_0^\ve=u_0^\ve+\ve w_0$ with
$u_0^\ve\to u_0$: three mixtures whose content of reactant differs by a quantity of order
$\ve$, and whose flame temperatures therefore differ by the same amount. The limit problems
differ through $M(w_0)$, and so do the
flames (Figure~\ref{fig:instability}, computed with BDF): after $t=20$ the fronts are at
$x\approx3.7$, $13.6$ and $35.7$. The dashed lines are the free boundaries of
\eqref{eq:P} with each $M(w_0)$, computed with the solver of Section~\ref{sec:limit}; the
solutions of \eqref{eq:Peps} follow them. For $w_0=-0.5$ the initial gradient exceeds
$\sqrt{2M(-0.5)}$ and the front first recedes before propagating.

\subsection{An accelerating front: the setting of Theorem~\ref{cv:thm:main}}\label{sec:exp:thm}

For the extinction of the self-similar ball (Section~\ref{sec:selfsimilar}, $w_0=0.5$,
$T=1$) the front accelerates, because it is concave towards the fresh gas, and the
solution is not steady in any frame; it satisfies (V1), (V2) on $[0,T^*]$ for every
$T^*<1$, and on the time intervals that we use the speed of the front is at most
$R/(2\sqrt{T-t})\le R$. We first run (S) in the setting of
Theorem~\ref{cv:thm:main}, without remeshing, for $N=2$, $0\le t\le0.5$,
$\ve\in\{0.02,0.01,0.005\}$ and $h_f=\ve/16$, and record
$E_v=\max_{n,j}|U^n_j-v(x^n_j,t_n)|$. On a uniform mesh that does not move, $E_v$ is a
function of $\Delta t/\ve$ only: $0.023$, $0.077$, $0.18$, $0.26$ and $0.30$ for
$\Delta t/\ve=1/64$, $1/16$, $1/4$, $1$ and $4$, the same for the three values of $\ve$ to two
digits. On a band mesh whose band moves with the prescribed velocity $v_n=\dot R(t_n)$,
so that $\mu=0$ (the cells outside the rigid zone are stretched, up to $0.044$),
$E_v=\ve w_0$ for every $\Delta t\le\ve$, which is the value of $v-u^\ve$ in the burnt region, and
$E_v\le0.67\,\ve$ for $\Delta t=4\ve$. This is the term $\mu\min\{1,\Delta t/\ve\}$ of the theorem,
against $\Delta t$ (the constants $c_0$, $C_b$ are not explicit, so we cannot say that
\eqref{cv:eq:hyp} holds in these runs).

\begin{figure}[ht]
\centering
\includegraphics[width=0.6\textwidth]{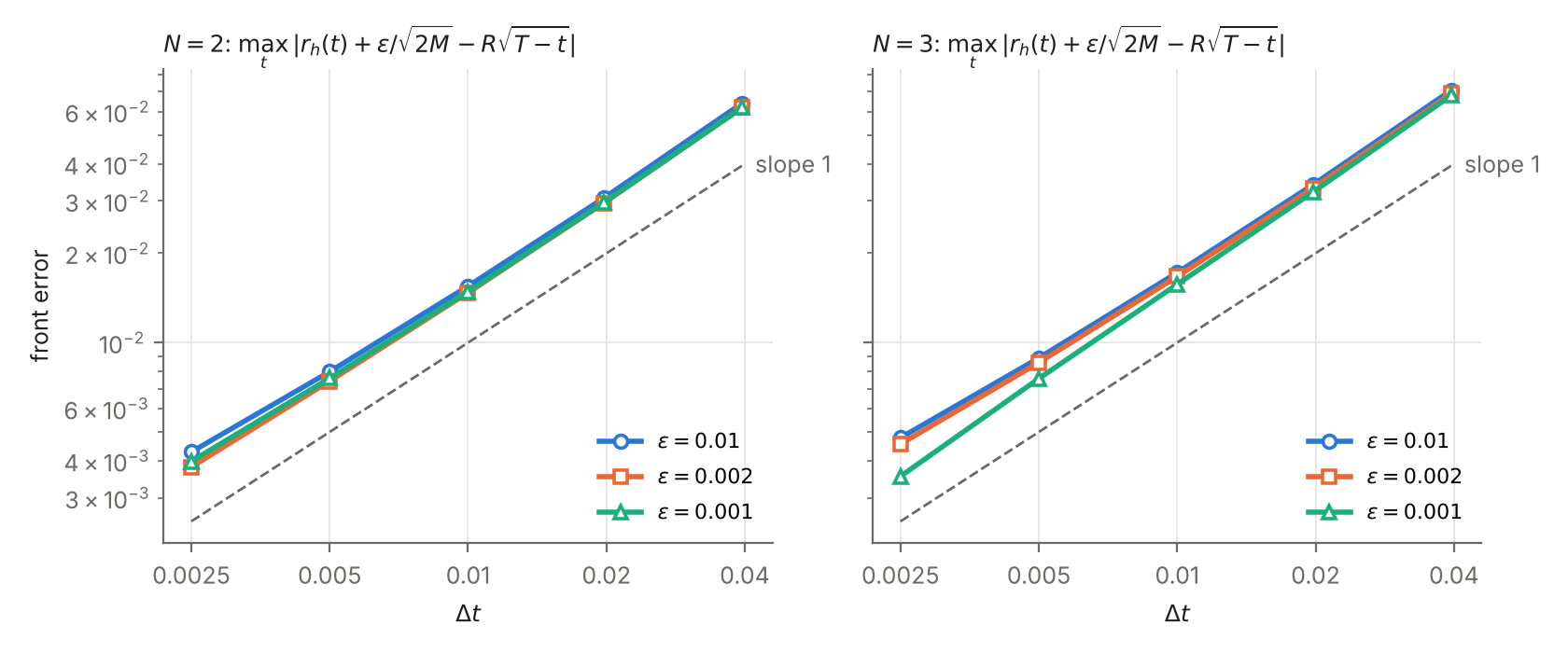}
\caption{ALE, extinction of a self-similar ball of fresh gas ($w_0=0.5$, $t\le0.75$): error in the
position of the front, after removing the shift $\ve/\sqrt{2M}$, against $\Delta t$, for
$N=2,3$ and $\ve\in\{0.01,0.002,0.001\}$.}\label{fig:B4}
\end{figure}

With the method itself (computed velocity, remeshing) the error in the position of the
front (Figure~\ref{fig:B4}), after
removing the shift $\ve/\sqrt{2M}$ of Remark~\ref{q:rem:rate}, is of first order in $\Delta t$:
between $1.4\,\Delta t$ and $1.9\,\Delta t$ for $\ve\in\{0.01,0.002,0.001\}$ and $N=2,3$, with no
systematic dependence on $\ve$. For general data the time step does enter, but as
$\Delta t$ and not as $\Delta t/\ve$, as in Theorem~\ref{cv:thm:main}; here $\Delta t/\ve$ ranges from $0.25$
to $40$, mostly beyond the time steps that the theorem covers.

\subsection{Travelling waves, and the cost of the method}\label{sec:exp:tw}

For travelling waves the exact solution $\Phi_\ve$ of \eqref{eq:Peps} is available, and the
discretization error can be separated from the error of the singular limit. We summarize
here the results; the details, with figures and tables, are in a supplementary
document available with the code.

\emph{Meshes that follow the front.} Starting from $U^0=\Phi_\ve$ with $w_0=0$, $A=1$ and
$h_f=\ve/16$, the method of Section~\ref{sec:ale} gives, up to $T=8$, a sup error
$E=4.2$--$4.8\cdot10^{-4}$ for every time step from $0.0025$ to $0.1$ and every $\ve$ from
$0.1$ to $0.001$, that is, for $\Delta t/\ve$ up to $100$: the error depends neither on the
time step nor on $\ve$, in agreement with Theorem~\ref{an:thm:ale}. It does depend on the
position of the initial front relative to the mesh (between $10^{-5}$ and $10^{-3}$ over
sixteen positions): it is a
\emph{lattice pinning} error: the front keeps its position $\vartheta h_f$ inside a cell of
the band and travels with a speed $c_\ve+e(\vartheta)$, where $e(\vartheta)$ oscillates with an
amplitude proportional to $(h_f/\ve)^2$ ($1.2\cdot10^{-4}$ for $h_f=\ve/16$). This is the
scaling $(h_f/\ve)^2$, not the $h_f/\ve$ of Theorem~\ref{cv:thm:main}; the coarse mesh
contributes $O(h_c^2)$. The computed mesh velocities differ from $c_\ve$ by at most
$1.6\cdot10^{-3}$, the largest differences occurring in the first steps. Inside the hypotheses of Theorem~\ref{an:thm:ale} (uniform mesh
translating with velocity $c_\ve$, $h\le\ve^2/8$) the error is below $10^{-5}$ and does not
depend on the time step either, and for a rate of class $C^{1,1}$ it is roughly
proportional to $h^2/\ve^3$; the constants of the bounds are far from sharp.

\emph{Meshes that do not move.} With the Eulerian version of (S) the front is slower than
the exact one, and the error of the speed is proportional to $\Delta t/\ve$: after a slow
transient it is between $-0.216$ and $-0.235$ times $\Delta t/\ve$, to be compared with the
value \eqref{an:eq:solv}. This is the term $\Delta t/\ve$ of Remark~\ref{cv:rem:main}(1). On a
rigid mesh that translates with a prescribed velocity $c_\ve+\delta v$ the front moves with
speed $c_\ve+\alpha\,\delta v$, where $\alpha$ grows from $0.04$ to $0.98$ as $\Delta t/\ve$ goes from
$0.03$ to $100$: this is the factor $\min\{1,\Delta t/\ve\}$ that multiplies $\mu$ in
Theorem~\ref{cv:thm:main}.

\emph{Cost.} For the same travelling wave up to $T=4$, the error in the position of the
front of the moving mesh method is fixed by the mesh of the band ($4.5\cdot10^{-4}$,
$10^{-4}$ and $10^{-5}$ for $h_f=\ve/16$, $\ve/32$ and $\ve/64$), with $240$ tridiagonal solves
on meshes of $3600$ to $6300$ nodes, for $\ve=0.01$ and $\ve=0.001$ alike. A method of lines on
the same band meshes without motion (BDF or LSODA) reaches errors between $10^{-5}$ and
$6\cdot10^{-5}$ with at least $700$ time steps for $\ve=0.01$ and more than $2\cdot10^4$ for
$\ve=0.001$, and the Eulerian version of (S) needs $\Delta t\approx\ve/1000$ for an error of
$10^{-3}$. The travelling wave is the most favourable case for a mesh that follows the
front, because the solution is steady in its frame; for other data the error is of
first order in $\Delta t$ (Section~\ref{sec:exp:thm}). We have not compared with moving mesh
methods based on equidistribution \cite{HuangRussell2011}.

\subsection{Approach to the self-similar extinction profile}

\begin{figure}[ht]
\centering
\includegraphics[width=0.85\textwidth]{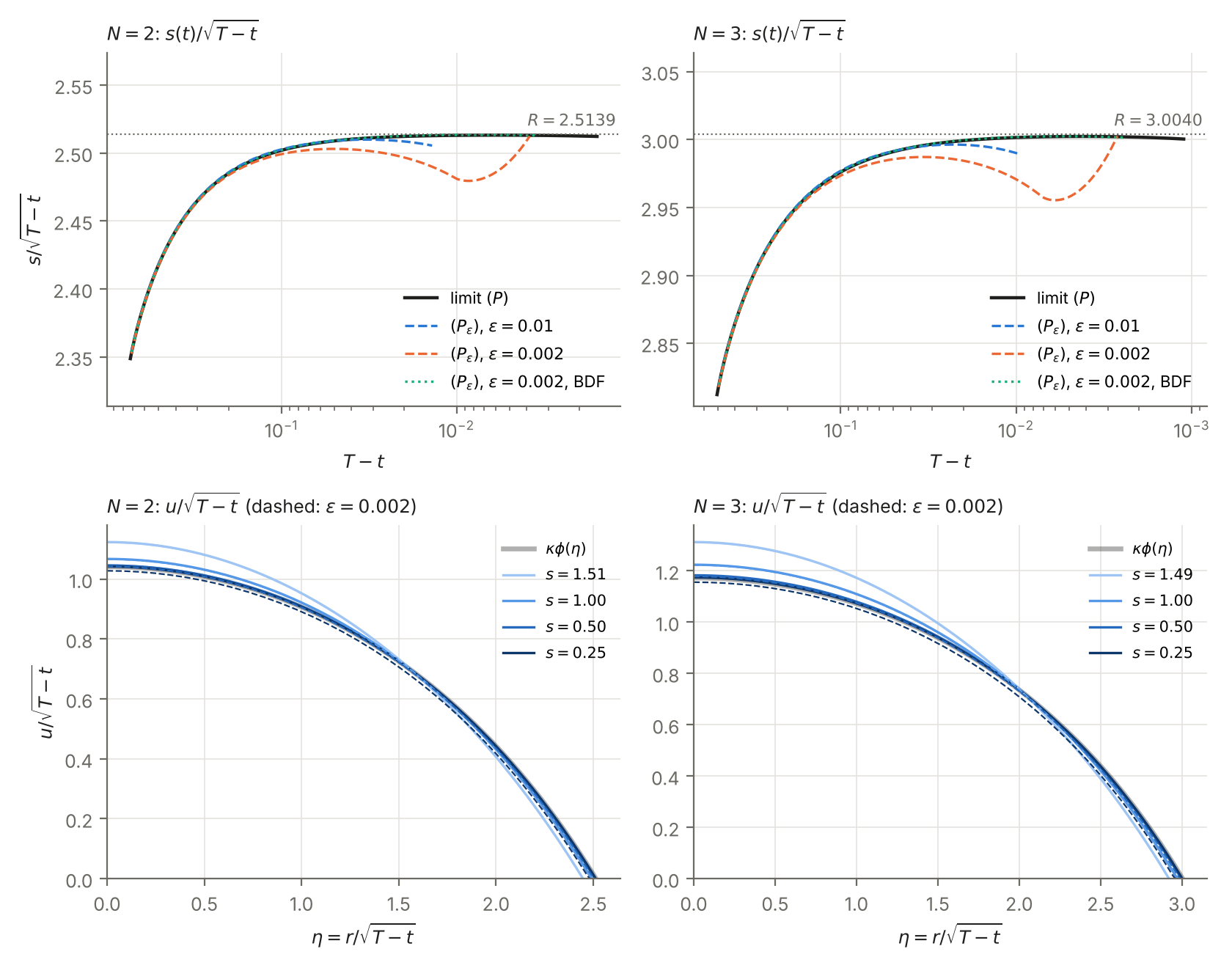}
\caption{Top: $s(t)/\sqrt{T-t}$ against $T-t$ for the limit problem (black) and for
\eqref{eq:Peps} (ALE, $\ve=0.01,0.002$; BDF, $\ve=0.002$), $N=2,3$. Bottom: rescaled
profiles $u/\sqrt{T-t}$ against $r/\sqrt{T-t}$ when the radius is $1.5$, $1$, $0.5$,
$0.25$, and the self-similar profile $h$ (grey).}\label{fig:C2}
\end{figure}

Galaktionov, Hulshof
and V\'azquez \cite{GalaktionovHulshofVazquez1997} proved that the extinction of
radial flames is asymptotically self-similar (see \cite{AbramsChoi2019} for non-radial
data). For the radial data above, followed until
extinction, Figure~\ref{fig:C2} shows $s(t)/\sqrt{T-t}$, which tends to the constant $R$
of Section~\ref{sec:selfsimilar} ($2.5139$ for $N=2$, $3.0040$ for $N=3$), and the rescaled
profiles, which tend to $h$. For the limit problem the ratio is $2.4474$, $2.4922$,
$2.5102$, $2.5130$ when the radius is $1.5$, $1$, $0.5$, $0.25$ ($N=2$). The extinction times
of \eqref{eq:Peps} agree with those of \eqref{eq:P} to three or four digits: for $N=2$,
$T=0.72509$ for \eqref{eq:P}, and $0.72502$ (BDF, $\ve=0.002$), $0.72464$ and $0.72533$
(ALE, $\ve=0.01$ and $0.002$). For $\ve=0.002$ the ALE curve dips by about $1.5\%$ below $R$
near extinction; this is the cut-off of the mesh velocity near $r=0$
(Section~\ref{sec:ale}), which makes the scheme partly Eulerian there.

\section{Conclusions}\label{sec:conclusions}

When the free boundary problem \eqref{eq:P} has a regular classical solution and the
data are strictly starshaped, the solutions of \eqref{eq:Peps} converge to it with
rate $O(\ve)$ from above and $O(\ve\log(1/\ve))$ from below
(Theorem~\ref{q:thm:rate}), up to the distance $\eta_\ve$ between $w^\ve/\ve$ and its limit;
this makes the convergence theorems of
\cite{FernandezBonderWolanski2004} quantitative, under stronger assumptions on the
regularity of the super- and subsolutions and on $f$. For constant $w_0$ the new
information with respect to the barriers of \cite{Petrosyan2002} is the absence of the
logarithm in the upper bound. For planar waves the first-order
correction of the speed is explicit (Proposition~\ref{prop:speed}).

On the numerical side, the linearly implicit scheme (S) satisfies a discrete
comparison principle for every time step and every admissible mesh, and the barriers
of the rate theorem, with a dilation that depends on time and a damping proportional
to the truncation errors of the layer, are discrete barriers. For data in planar or
radial symmetry, a regular solution whose radius satisfies $R-2t\dot R>0$, a perturbation
$w_0$ that is monotone under parabolic dilations, and band
meshes with prescribed velocities, the fully discrete solution
converges to the solution of the free boundary problem with the bound
$C(\ve\log\frac1\ve+\frac{h_f}\ve+h_c+\Delta t+\mu\min\{1,\frac{\Delta t}\ve\})$
(Theorem~\ref{cv:thm:main}), for time steps $\Delta t\lesssim\ve$ divided by the bound: the layer
has to be resolved, and this is enough; on meshes that do not move the time step is divided
by $\ve$, and on meshes that follow the front it is not. In the experiments the terms in
$\Delta t$ and $\mu$ are observed, also for $\Delta t\gg\ve$, while the mesh sizes enter as
$(h_f/\ve)^2$ and $h_c^2$.

The scope of these results is that of the model: constant density, no heat losses and
equal diffusivities. The last assumption is what makes $w^\ve$ caloric and reduces the
system to a scalar equation with a comparison principle, on which both the rate and the
scheme rest; with a Lewis number different from one, or with heat losses, this structure
is lost, and neither the barriers nor the monotonicity of the scheme extend directly. The
rate theorem and the convergence theorem concern pockets of fresh gas; planar and
expanding flames, for which the fresh region is unbounded, are covered only by the
travelling wave estimates and by the experiments. And the experiments are in planar or
radial symmetry, in the regime of very large activation energy.

Several questions remain open: (i) an a priori bound for the implicit mesh velocity,
the remeshing, time steps $\Delta t\gg\ve$ for general data, and the sharp dependence on the
mesh sizes (Section~\ref{an:sec:impl});
(ii) discrete barriers without symmetry, for which the barriers of
Section~\ref{sec:rate} are available in any dimension, and for travelling waves and
other configurations with unbounded fresh region; (iii)
whether the logarithm in Theorem~\ref{q:thm:rate} is necessary, and rates for limits that
are not classical, or without starshapedness, in particular near the non-uniqueness and
degenerate behaviours of \cite{PetrosyanYip2008,WeissZhang2009}; (iv) in the method: the
lattice pinning of the front, the cut-off of the mesh velocity near the origin, which
makes the scheme partly Eulerian near extinction, and several fronts; and (v) the
problem in two and three dimensions, with adaptive refinement and coarsening of a mesh
that follows the flame, or with a moving mesh method
\cite{HuangRussell2011}, for which the elliptic (Bernoulli) case
\cite{FlucherRumpf1997,ShenDaiHuang2024} is a natural first step.

\section*{Use of AI tools}

The problem and the main ideas and the original draft of this work are due to the author. The work was prepared
with the assistance of Claude Code (Anthropic). The tool wrote all the numerical code, ran the simulations and produced the figures and tables; it assisted in the development and writing of the proofs, in drafting and
revising the manuscript, and in searching the literature. The author reviewed and
verified all the statements, proofs, computations and references, and takes full
responsibility for the content of the article.

\bibliographystyle{amsplain}
\bibliography{references}

\end{document}